\documentclass[final]{siamltex}
\usepackage{amsfonts}
\usepackage{amsmath}
\usepackage{amssymb}
\usepackage{enumerate}
\usepackage{graphicx}
\usepackage{multirow}
\usepackage{url}
\usepackage{verbatim}
\usepackage[pdftex]{thumbpdf}
\usepackage[pdftex,
pdfstartview=FitH,bookmarks=true,bookmarksnumbered=true,bookmarksopen=true,hypertexnames=false,breaklinks=true,
colorlinks=true,  linkcolor=blue,anchorcolor=blue,
citecolor=blue,filecolor=blue,  menucolor=blue]{hyperref}
\providecommand{\U}[1]{\protect\rule{.1in}{.1in}}
\newtheorem{algorithm}[theorem]{Algorithm}

\newtheorem{remark}[theorem]{Remark}
\begin{document}

\title{Inverse subspace iteration for spectral \\stochastic finite element methods \thanks{This work is based upon work
supported by the U.\, S.\, Department of Energy Office of Advanced Scientific
Computing Research, Applied Mathematics program under Award Number
DE-SC0009301, and by the U.\, S.\, National Science Foundation under grants
DMS–1418754 and DMS1521563.}}
\author{Bed\v{r}ich Soused\'{\i}k\thanks{Department of Mathematics and Statistics,
University of Maryland, Baltimore County, 1000 Hilltop Circle, Baltimore,
MD~21250 (\texttt{sousedik@umbc.edu}). }
\and Howard C. Elman\thanks{Department of Computer Science and Institute for
Advanced Computer Studies, University of Maryland, College Park, MD 20742
(\texttt{elman@cs.umd.edu}) }}
\maketitle

\begin{abstract}
We study random eigenvalue problems in the context of spectral stochastic
finite elements. In particular, given a parameter-dependent, symmetric
positive-definite matrix operator, we explore the performance of algorithms
for computing its eigenvalues and eigenvectors represented using polynomial
chaos expansions.
We formulate a version of stochastic inverse subspace iteration, which is
based on the stochastic Galerkin finite element method, and we compare its
accuracy with that of Monte Carlo and stochastic collocation methods.
The coefficients of the eigenvalue expansions are computed from a stochastic
Rayleigh quotient. 
Our approach allows the computation of  interior eigenvalues by deflation methods, 
and we can also compute the coefficients of multiple eigenvectors 
using a stochastic variant of the modified Gram-Schmidt process.
The effectiveness of the methods is illustrated by numerical experiments on
benchmark problems arising from vibration analysis.
\end{abstract}

\section{Introduction}

\label{sec:introduction}Eigenvalue analysis plays an essential role in many
applications, for example, dynamic response of structures, stability of flows,
and nuclear reactor criticality calculations. In traditional approaches, the
physical characteristics of models are considered to be known and the
eigenvalue problem is deterministic. However, in many important cases there is
uncertainty, for example, due to material imperfections, boundary conditions
or external loading, and the exact values of physical parameters are not
known. If the parameters are treated as random processes, the associated
matrix operators have a random structure as well, and the uncertainty is
translated into eigenvalues and eigenvectors. Techniques used to solve this
class of problems include Monte Carlo
methods~\cite{Nightingale-2007-MCE,Pradlwarter-2002-REP}, which are known to
be robust but slow, and perturbation
methods~\cite{Kaminski-2013-SPM,Kleiber-1992-SFE,Shinozuka-1972-REP,vomScheidt-1983-REP}, which are limited to models with low variability of uncertainty.

In this study, we explore the use of spectral stochastic finite element
methods (SSFEM)~\cite{Ghanem-1991-SFE,LeMaitre-2010-SMU,Xiu-2010-NMS} for the
solution of eigenvalue problems. The methods are based on an assumption that
the stochastic process is described in terms of polynomials of random
variables, and they produce discrete solutions that, with respect to the
stochastic component, are also polynomials in these random variables.
This framework is known as the generalized polynomial chaos
(gPC)~\cite{Ghanem-1991-SFE,Xiu-2002-WAP}. There are two main ways to use
this approach: stochastic Galerkin finite elements and stochastic
collocation~(SC). The first method translates the stochastic problem by means
of Galerkin projection into one large coupled deterministic system; the second
method samples the model problem at a predetermined set of \emph{collocation
points}, which yields a set of uncoupled deterministic problems. Although
numerous algorithms for solving stochastic partial differential equations by
SSFEM have been proposed, the literature addressing eigenvalue problems is
limited. Verhoosel~et al.~\cite{Verhoosel-2006-ISR} proposed an algorithm for
inverse iteration in the context of stochastic Galerkin finite elements, and
Meidani and Ghanem~\cite{Meidani-2012-SMD,Meidani-2014-SPI} proposed
stochastic subspace iteration using a stochastic version of the QR algorithm.
In alternative approaches, Ghanem and
Ghosh~\cite{Ghanem-2007-ECR,Ghosh-2008-SCA} proposed two numerical schemes:
one based on the Newton-Raphson method, and another based on an optimization
problem (see also~\cite{Ghosh-2013-ARE,Sarrouy-2012-SAE}), Pascual and
Adhikari~\cite{Pascual-2012-HPP} introduced several hybrid
perturbation-Polynomial Chaos approaches, and
Williams~\cite{Williams-2010-MSSnote,Williams-2010-MSS,Williams-2013-MSS}
presented a method that avoids the nonlinear terms in the conventional method
of stochastic eigenvalue calculation but introduces an additional independent variable.

We formulate a version of stochastic inverse subspace iteration which is based
on the stochastic Galerkin finite element method. We assume that the symmetric
positive-definite matrix operator is given in the form of a polynomial chaos
expansion, and we compute the coefficients of the polynomial chaos expansions
of the eigenvectors and eigenvalues. We also compare this method with the
stochastic collocation method in the larger context of spectral stochastic
finite element methods. In particular, we use both these methods to explore
stochastic eigenvalues and give an assessment of their accuracy. Our starting
point for stochastic inverse subspace iteration is based
on~\cite{Meidani-2014-SPI,Verhoosel-2006-ISR}. In order to increase efficiency
of the algorithm, we first solve the underlying mean problem and we use the
solution as the initial guess for the stochastic inverse subspace iteration,
which computes a correction of the expected value of the eigenvector from the
mean and coefficients of the higher order terms in the gPC\ expansion. The
gPC\ coefficients of the eigenvalue expansions are computed from a stochastic
Rayleigh quotient. We also show that in fact the Rayleigh quotient itself
provides a fairly close estimate of the eigenvalue expansion using only the
mean coefficients of the corresponding eigenvector. In our approach, it is
also relatively easy to deal with badly separated eigenvalues because one can
apply deflation to the mean matrix in the same way as in the deterministic case.\ 

The paper is organized as follows. In Section~\ref{sec:problem} we introduce
the stochastic eigenvalue problem and outline the Monte Carlo, stochastic
collocation and stochastic Galerkin methods. In Section~\ref{sec:SISI} we
formulate the algorithm of stochastic inverse subspace iteration. In
Section~\ref{sec:numerical} we report the results of numerical experiments,
and in Section~\ref{sec:conclusion} we summarize and conclude our work.

\section{Stochastic eigenvalue problem}

\label{sec:problem}Let $\left(  \Omega,\mathcal{F},\mathcal{P}\right)  $ be a
complete probability space, that is, $\Omega$ is the sample space
with $\sigma$-algebra $\mathcal{F}$ and probability measure~$\mathcal{P}$,
and let $D\subset\mathbb{R}^{d}$\ be a bounded physical domain. 
We assume that the randomness in the mathematical model is induced by 
a vector~$\xi:\Omega\rightarrow\Gamma 
\subset\mathbb{R}^{m_{\xi}}$ of independent, identically distributed (i.i.d.) 
random variables $\xi_{1}(\omega), \dots ,\xi_{m_{\xi}}(\omega)$.  Let 
$\mathcal{B}(\Gamma)$ denote 
the Borel $\sigma$-algebra on~$\Gamma$ induced by~$\xi$ 
and $\mu$  the induced measure.  
Then, the expected
value of the product of measurable functions on~$\Gamma$ 
determines a Hilbert 
space$~L^{2}\left(  \Gamma,\mathcal{B}(\Gamma),\mu\right)$ 
with inner product 
\begin{equation}
\left\langle u,v\right\rangle =\mathbb{E}\left[ u v \right]  =\int_{\Gamma
} u(\xi) v(\xi)  \,d\mu (\xi)  ,
\label{eq:expectation}
\end{equation}
where the symbol$~\mathbb{E}$\ denotes the mathematical expectation.

In computations, we will work with a set $\left\{  \psi_{\ell} \right\}$ of 
orthonormal polynomials 
such that $\left\langle \psi_{j}\psi_{k}\right\rangle =\delta_{jk}$, where
$\delta_{jk}$ is the Kronecker delta and $\psi_0$  is constant.
This set, the gPC\ basis, spans a finite-dimensional subspace of 
$L^{2}\left(  \Gamma,\mathcal{B}(\Gamma),\mu\right)$.
We will also suppose we are given a
symmetric positive-definite matrix-valued random 
variable~$A\left( x,\xi \right)$ represented as
\begin{equation}
A\left( x, \xi\right)  =\sum_{\ell=0}^{M_{A}}A_{\ell} \left( x \right) \psi_{\ell}\left(
\xi\right)  ,
\label{eq:gPC-A}
\end{equation}
where each $A_{\ell}$ is a deterministic\ matrix of size$~M_{x}\times
M_{x}$, with $M_{x}$ determined by the discretization of the physical
domain, and $A_{0}$ is the matrix corresponding to the mean value 
of the matrix~$A\left(x,\xi\right)$, that is $A_0=\mathbb{E}\left[A\left(x,\xi\right)\right]$. 
The representation~(\ref{eq:gPC-A}) is typically obtained
from an expansion of a random process;
examples are given in Section~\ref{sec:numerical} on numerical experiments. 
We will also use the notation 
\begin{equation}
c_{\ell jk}=\mathbb{E}\left[  \psi_{\ell} \psi_{j} \psi_{k} \right] .
\label{eq:cijk}\end{equation}
We are interested in a solution of the following stochastic eigenvalue
problem: find a set of random eigenvalues$~\lambda^{s}$ and corresponding
eigenvectors$~u^{s}$, $s=1,\dots n_{s}$, which almost surely (a.s.) satisfy
the equation
\begin{equation}
A\left(  x,\xi\right)  u^{s}\left(  x,\xi\right)  =\lambda^{s}\left(
\xi\right)  u^{s}\left(  x,\xi\right)  ,\qquad\forall x\in D.
\label{eq:standard} 
\end{equation}
We will search for expansions of\ a set of$~n_{s}$ eigenvalues and
eigenvectors in the form
\begin{equation}
\lambda^{s}\left(  \xi\right)  =\sum_{k=0}^{M_{\lambda}}\lambda_{k}^{s} 
\psi_{k}\left(  \xi\right)  ,\qquad u^{s}\left(  x,\xi\right)  =
\sum_{k=0}^{M_{\xi}}u_{k}^{s}\left(  x\right)  \psi_{k}\left(  \xi\right)  ,
\label{eq:gPC-lambda-u} 
\end{equation}
where $\lambda_{k}^{s}$ and $u_{k}^{s}$ are\ coefficients defined by
projection on the basis $\left\{  \psi_{k} \right\}$, 
\begin{equation}
\lambda_{k}^{s}=\left\langle \lambda^{s},\psi_{k}\right\rangle ,
\quad k=0,\dots, M_{\lambda}, 
\qquad
u_{k}^{s}=\left\langle u^{s},\psi_{k}\right\rangle ,
\quad
 k=0,\dots, M_{\xi}.
\label{eq:coeffs-lambda-u} 
\end{equation}
We will consider several ways to approximate these quantities.

\subsection{Monte Carlo and stochastic collocation methods}

Both the Monte Carlo and the stochastic collocation methods are based on
sampling. This entails the solution of independent deterministic eigenvalue
problems at a set of sample points$~\left\{  \xi^{\left(  q\right)  }\right\}$,
\begin{equation}
A\left(  \xi^{\left(  q\right)  }\right)  u^{s}\left(  \xi^{\left(  q\right)
}\right)  =\lambda^{s}\left(  \xi^{\left(  q\right)  }\right)  u^{s}\left(
\xi^{\left(  q\right)  }\right)  ,\qquad s=1,\dots,n_{s}.
\label{eq:standard-sample} 
\end{equation}
In the Monte Carlo method,  
the sample points $\xi^{\left( q \right)}$, $q=1,\dots,N_{MC}$, 
are generated randomly, following the distribution 
of the random variables$~\xi_i$, $i=1,\dots,m_\xi$ 
and moments of the solution are obtained from ensemble averaging. 
For stochastic collocation, the sample points~$\xi^{\left(  q\right)}$, 
$q=1,\dots,N_{q}$, consist of a predetermined set of \emph{collocation points}. 
This approach derives from a methodology for performing quadrature or
interpolation in multidimensional space using a small number of points, a
so-called sparse
grid~\cite{Gerstner-1998-NIU,Novak-1996-HDI,Trefethen-2008-IGQ}.

There are two ways to implement stochastic collocation~\cite{Xiu-2010-NMS}. We
can either construct a Lagrange interpolating polynomial that interpolates at
the collocation points, or we can use a discrete projection in the so-called
pseudospectral approach, to obtain coefficients of expansion in an a priori
selected basis of stochastic functions. In this study, we use the second
approach because it facilitates a direct comparison with the stochastic
Galerkin method. In particular, the coefficients 
in the expansions~(\ref{eq:gPC-lambda-u}) are determined by 
evaluating~(\ref{eq:coeffs-lambda-u}) in the sense of~(\ref{eq:expectation}) 
using numerical quadrature as 
\begin{equation}
\lambda_{k}^{s}=\sum_{q=1}^{N_{q}}\lambda^{s}\,\psi_{k}\left(  \xi^{\left(
q\right)  }\right)  \,w^{\left(  q\right)  },\quad u_{ik}^{s}=\sum
_{q=1}^{N_{q}}u^{s}\left(  x_{i}\right)  \,\psi_{k}\left(  \xi^{\left(
q\right)  }\right)  \,w^{\left(  q\right)  }, 
\label{eq:Q-lambda-u} 
\end{equation}
where $\xi^{\left(  q\right)  }$ are the collocation (quadrature) points,\ and
$w^{\left(  q\right)  }$ are quadrature weights. 
We refer to \cite{LeMaitre-2010-SMU} for a discussion of quadrature rules.
Details of the rule we use in experiments are 
discussed in Section~4 (prior to Section~4.1).

\subsection{Stochastic Galerkin method}

The stochastic Galerkin method\ is based on the projection
\begin{equation}
\left\langle Au^{s},\psi_{k}\right\rangle =
\left\langle \lambda^{s}u^{s},\psi_{k}\right\rangle ,\qquad k=0,\dots,M_{\xi},
\quad s=1,\dots,n_{s}.
\label{eq:variational}
\end{equation}
Substituting the expansions~(\ref{eq:gPC-A}) and (\ref{eq:gPC-lambda-u}) into
(\ref{eq:variational}) yields a nonlinear algebraic system
\begin{equation}
\sum_{j=0}^{M_{\xi}}\sum_{\ell=0}^{M_{A}}c_{\ell jk}A_{\ell}u_{j}^{s}=
\sum_{j=0}^{M_{\xi}}\sum_{i=0}^{M_{\lambda}}c_{ijk}\lambda_{i}^{s}u_{j}^{s},
\qquad k=0,\dots,M_{\xi},\quad s=1,\dots,n_{s}, 
\label{eq:global-system}
\end{equation}
to solve for the coefficients $\lambda_{i}^{s}$ and $u_{j}^{s}$.
The Galerkin solution is then given by~(\ref{eq:gPC-lambda-u}).

We will also consider a shifted variant of this method with a deterministic
shift$~\rho$, introduced in~\cite{Verhoosel-2006-ISR}, which can be used to
find a single
interior eigenvalue, Thus we drop the superscript $s$, and the shifted
counterpart of~(\ref{eq:variational}) is
\begin{equation}
\left\langle \left(  A-\rho I\right)  u,\psi_{k}\right\rangle =\left\langle
\left(  \lambda-\rho\right)  u,\psi_{k}\right\rangle ,\qquad k=0,\dots,M_{\xi}.
\label{eq:variational-shift}
\end{equation}
Writing this equation out using the gPC\ expansions gives
\[
\left\langle \left(  \sum_{\ell=0}^{M_{A}}A_{\ell}\psi_{\ell}-\rho I\psi
_{0}\right)  \sum_{j=0}^{M_{\xi}}u_{j}\psi_{j},\psi_{k}\right\rangle =\left\langle
\left(  \sum_{i=0}^{M_{\lambda}}\lambda_{i}\psi_{i}-\rho\psi_{0}\right)
\sum_{j=0}^{M_{\xi}}u_{j}\psi_{j},\psi_{k}\right\rangle ,
\]
which leads to a modified system 
\begin{equation}
\sum_{j=0}^{M_{\xi}}\sum_{\ell=0}^{M_{A}}c_{\ell jk}\widetilde{A}_{\ell}u_{j}
=\sum_{j=0}^{M_{\xi}}\sum_{i=0}^{M_{\lambda}}c_{ijk}\widetilde{\lambda}_{i}%
u_{j},\qquad k=0,\dots,M_{\xi}, 
\label{eq:global-system-shift}
\end{equation}
instead of~(\ref{eq:global-system}), where
\begin{align}
\widetilde{A}_{0}  &  =A_{0}-\rho I,\qquad\widetilde{A}_{\ell}=A_{\ell}%
,\quad\ell=1,\dots,M_{A},\label{eq:A-shift}\\
\widetilde{\lambda}_{0}  &  =\lambda_{0}-\rho,\qquad\widetilde{\lambda}%
_{i}=\lambda_{i},\quad i=1,\dots,M_{\lambda}, 
\label{eq:lambda-shift} 
\end{align}
and $\left\{  \lambda_{i}\right\}  _{i=0}^{M_{\lambda}}$\ are the quantities
that would be obtained with$~\rho=0$.

As will be seen in numerical experiments, deflation
of the mean matrix $A_{0}$ is more robust than using a shift
for identification of interior eigenvalues.
For inverse iteration, deflation can be done via
\begin{equation}
\widetilde{A}_{0}=A_{0}+\sum_{d=1}^{n_{d}}\left[  C_{\lambda}-\lambda_{0}^{d}\right]
\left(  u_{0}^{d}\right)  \left(  u_{0}^{d}\right)  \,\!\!^{T},\qquad
\widetilde{A}_{\ell}=A_{\ell},\quad\ell=1,\dots,M_{A}, 
\label{eq:deflation-2} 
\end{equation}
where $\left(  \lambda_{0}^{d},u_{0}^{d}\right)  ,$ $d=1,\dots,n_{d}$, are the
zeroth order coefficients of the eigenpairs to be deflated, and $C_{\lambda}$ is a
constant such that $C_{\lambda}\gg\lambda_{0}^{d}$ for $d=1,\dots,n_{d}$, for example
$C_{\lambda}=\max_{s}\left(  \lambda^{s}\right)  $. Note that there are other types of
deflation, where, for example, the computation proceeds with a smaller
transformed submatrix from which the deflated eigenvalue is explicitly
removed. Since this is complicated for matrix operators in the form of the
expansion (\ref{eq:gPC-A}), we do not consider this approach here.

\section{Stochastic inverse subspace iteration}

\label{sec:SISI} The stochastic inverse subspace iteration algorithm is based
on a stochastic Galerkin projection. In order to motivate the
stochastic algorithm, we first recall the deterministic inverse subspace
iteration that allows to find several smallest eigenvalues and corresponding
eigenvectors of a given matrix. Then, we give a formal statement of the
stochastic variant and relate it to stochastic inverse
iteration~\cite{Verhoosel-2006-ISR}. Finally, we describe the components of
the algorithm in detail and relate it to stochastic subspace
iteration~\cite{Meidani-2014-SPI}. Our strategy is motivated by the
deterministic inverse subspace iteration, when the aim is to find small
eigenvalues by finding large eigenvalues of the inverse problem. 

\begin{algorithm}
[Deterministic inverse subspace iteration (DISI)]\label{alg:DISI}
Let $u^{1},\dots,u^{n_{s}}$ be a set of $n_{s}$ orthonormal
vectors, and let$~\mathcal{A}$ be a symmetric positive-definite matrix. 

\textbf{for} $it=0,1,2,\dots$

\begin{enumerate}
\item Solve the system for $v^{s,\left(  it\right)  }$: 
\begin{equation}
\mathcal{A}v^{s,\left(  it\right)  }=u^{s,\left(  it\right)  },\qquad
s=1,\dots,n_{s}. \label{eq:alg-DISI-solve} 
\end{equation}

\item If $n_{s}=1$, normalize as\ $u^{s,\left(  it+1\right)  }=v^{s,\left(
it\right)  }/\left\Vert v^{s,\left(  it\right)  }\right\Vert $, or else if
$n_{s}>1$ use the Gram-Schmidt process and transform the set of vectors
$v^{s,\left(  it\right)  }$ into an orthonormal set $u^{s,\left(  it+1\right)
},$ $s=1,\dots,n_{s}$.

\item Check convergence.
\end{enumerate}

\textbf{end}

Use the Rayleigh quotient to compute the eigenvalues $\lambda^{s}=\left(
u^{s,\left(  it+1\right)  }\right)  ^{T}\mathcal{A}u^{s,\left(  it+1\right)
}$.
\end{algorithm}

The proposed stochastic variant of this algorithm follows:

\begin{algorithm}
[Stochastic inverse subspace iteration (SISI)]\label{alg:SISI}Find~$n_{e}$
eigenpairs of the deterministic problem
\begin{equation}
A_{0}\,\overline{U}=\overline{U}\,\overline{\Lambda},\qquad\overline
{U}=\left[
\begin{array}
[c]{ccc}%
\overline{u}^{1}, & \dots, & \overline{u}^{n_{e}}%
\end{array}
\right]  ,\qquad\overline{\Lambda}=diag\left(  \overline{\lambda}^{1}%
,\dots,\overline{\lambda}^{n_{e}}\right)  , \label{eq:SISI_mean}%
\end{equation}
choose$~n_{s}$ eigenvalues of interest with indices $e=\left\{  e_{1}%
,\dots,e_{n_{s}}\right\}  \subset\left\{  1,\dots,n_{e}\right\}  $, set up 
matrices$~\mathcal{A}_{\ell}$, $\ell=0,\dots,M_A$, either as
$\mathcal{A}_{\ell}=A_{\ell}$, or $\mathcal{A}_{\ell}=\widetilde{A}_{\ell}$
using deflation such as~(\ref{eq:deflation-2}) and initialize
\begin{align}
u_{0}^{1,\left(  0\right)  }  &  =\overline{u}^{e_{1}},\quad u_{0}^{2,\left(
0\right)  }=\overline{u}^{e_{2}},\quad\dots\quad u_{0}^{n_{s},\left(
0\right)  }=\overline{u}^{e_{n_{s}}},\label{eq:SISI-u_i_1}\\
u_{i}^{s,\left(  0\right)  }  &  =0,\qquad s=1,\dots,n_{s},\quad i=1,\dots,M_{\xi}.
\label{eq:SISI-u_i_2}%
\end{align}

\textbf{for} $it=0,1,2,\dots$

\begin{enumerate}
\item Solve the stochastic Galerkin system for $v_{j}^{s,\left(  it\right)  }%
$, $j=0,\dots,M_{\xi}$:%
\begin{equation}
\sum_{j=0}^{M_{\xi}}\sum_{\ell=0}^{M_{A}}c_{\ell jk}\mathcal{A}_{\ell}%
v_{j}^{s,\left(  it\right)  }=u_{k}^{s,\left(  it\right)  },\qquad
k=0,\dots,M_{\xi},\quad s=1,\dots,n_{s}. \label{eq:alg-SISI-solve}%
\end{equation}

\item If $n_{s}=1$, normalize using the quadrature
rule~(\ref{eq:vector-normalize}), or else if $n_{s}>1$ orthogonalize using the
stochastic modified Gram-Schmidt process: transform the set of coefficients
$v_{j}^{s,\left(  it\right)  }$ into a set $u_{j}^{s,\left(  it+1\right)  },$
$j=0,\dots,M_{\xi}$, $s=1,\dots,n_{s}$.

\item Check convergence.
\end{enumerate}

\textbf{end}

Use the stochastic Rayleigh quotient~(\ref{eq:RQ}) to compute the eigenvalue expansions.
\end{algorithm}

\medskip

Stochastic inverse iteration \cite[Algorithm~2]{Verhoosel-2006-ISR},
corresponds to the case where a stochastic expansion of a single eigenvalue
(in~\cite{Verhoosel-2006-ISR} with $M_{\lambda}=M_{\xi}$) is
sought; in this case, we can select a shift~$\rho$ using the solution of the
mean problem (\ref{eq:SISI_mean}) and modify the mean matrix
using~(\ref{eq:A-shift}). Step~$1$ of Algorithm~\ref{alg:SISI} then consists
of two parts:

\begin{enumerate}
\item[1(a).] Use the stochastic Rayleigh quotient~(\ref{eq:RQ}) to compute the
coefficients $\lambda_{i}^{\left(  it\right)  }$, $i=0,\dots,M_{\lambda}$ of
the eigenvalue expansion~(\ref{eq:gPC-lambda-u}), and set up the right-hand
side components as 
\[
\Lambda_{0}^{\left(  it\right)  }=\left(  \lambda_{0}^{\left(  it\right)
}-\rho\right)  I_{M_{x}},\qquad\Lambda_{i}^{\left(  it\right)  }=\lambda
_{i}^{\left(  it\right)  }I_{M_{x}},\quad i=1,\dots,M_{\lambda},
\]
where $I_{M_{x}}$ is the identity matrix of size$~M_{x}$.

\item[1(b).] Solve the stochastic Galerkin system for $v_{j}^{\left(
it\right)  }$, $j=0,\dots,M_{\xi}$:%
\begin{equation}
\sum_{j=0}^{M_{\xi}}\sum_{\ell=0}^{M_{A}}c_{\ell jk}\widetilde{A}_{\ell}%
v_{j}^{\left(  it\right)  }=\sum_{j=0}^{M_{\xi}}\sum_{i=0}^{M_{\lambda}}%
c_{ijk}\Lambda_{i}^{\left(  it\right)  }u_{j}^{\left(  it\right)  },\qquad
k=0,\dots,M_{\xi}. \label{eq:alg-SII-solve}%
\end{equation}

\end{enumerate}

\begin{remark}
In the deterministic version of inverse iteration, the shift$~\left(
\lambda-\rho\right)$ applied to the vector on the right-hand side is
dropped, and each step entails solution of the system
\begin{equation}
\left(  A_{0}-\rho I\right)  u_{0}^{\left(  it+1\right)  }=u_{0}^{\left(
it\right)  }. \label{eq:inv-it-deterministic}%
\end{equation}
Moreover, in the deterministic version of Rayleigh quotient iteration, the
eigenvalue estimate$~\lambda^{\left(  it\right)  }$ is used instead of$~\rho$
in (\ref{eq:inv-it-deterministic}).
Here we retain the shift in the iteration due to the presence of the
stochastic Galerkin projection, see~(\ref{eq:variational-shift}). In
particular, the shift in the left-hand side is fixed to$~\rho$ and the
estimate of the eigenvalue expansion is used in the setup of the right-hand
side in~(\ref{eq:alg-SII-solve}). Thus, stochastic inverse iteration 
is not an exact counterpart of either deterministic inverse iteration or
Rayleigh quotient\ iteration.
\end{remark}

We now describe components of Algorithm~\ref{alg:SISI} in detail.

\paragraph{Matrix-vector multiplication}

Computation of the stochastic Rayleigh quotient requires a stochastic version
of a matrix-vector product, which corresponds to evaluation of the projection
\[
v_{k}=\left\langle v,\psi_{k}\right\rangle =\left\langle Au,\psi
_{k}\right\rangle ,\qquad k=0,\dots,M_{\xi}.
\]
In more detail, this is
\[
\left\langle \sum_{i=0}^{M_{\xi}}v_{i}\psi_{i}\left(  \xi\right)  ,\psi
_{k}\right\rangle =\left\langle \left(  \sum_{\ell=0}^{M_{A}}A_{\ell}%
\psi_{\ell}\left(  \xi\right)  \right)  \left(  \sum_{j=0}^{M_{\xi}}u_{j}\psi
_{j}\left(  \xi\right)  \right)  ,\psi_{k}\right\rangle ,\qquad k=0,\dots,M_{\xi},
\]
so the coefficients of the expansion of the vector$~v$ are
\begin{equation}
v_{k}=\sum_{j=0}^{M_{\xi}}\sum_{\ell=0}^{M_{A}}c_{\ell jk}A_{\ell}u_{j},\qquad
k=0,\dots,M_{\xi}. \label{eq:mat-vec}%
\end{equation}
The use of this computation for the Rayleigh quotient is described below.
Algorithm~\ref{alg:SISI} can also be modified to perform subspace
iteration~\cite[Algorithm~$4$]{Meidani-2014-SPI} for identifying the largest
eigenvalue of$~A$. In this case, the solve in Step$~1$ of
Algorithm~\ref{alg:SISI} is replaced by a matrix-vector 
product~(\ref{eq:mat-vec}) in which $u_{j}=u_{j}^{s}$, $v_{k}=v_{k}^{s}$.

\paragraph{Stochastic Rayleigh Quotient}

In the deterministic case, the Rayleigh quotient\ is used to compute the
eigenvalue corresponding to a normalized eigenvector$~u$ as $\lambda=u^{T}v$,
where $v=Au$. For the stochastic Galerkin method, the Rayleigh quotient
defines the coefficients of a stochastic expansion of the eigenvalue defined
via a projection 
\[
\lambda_{k}=\left\langle \lambda,\psi_{k}\right\rangle =\left\langle
u^{T}v,\psi_{k}\right\rangle ,
\qquad k=0,\dots, M_{\lambda}.
\]
In our implementation we used $M_{\lambda}=M_{\xi}$. 
The coefficients of$~v$ are computed simply
using the matrix-vector product~(\ref{eq:mat-vec}). In more detail, this is
\[
\left\langle \sum_{i=0}^{M_{\xi}}\lambda_{i}\psi_{i}\left(  \xi\right)  ,\psi
_{k}\right\rangle =\left\langle \left(  \sum_{i=0}^{M_{\xi}}u_{i}\psi_{i}\left(
\xi\right)  \right)  ^{T}\left(  \sum_{j=0}^{M_{\xi}}v_{j}\psi_{j}\left(
\xi\right)  \right)  ,\psi_{k}\right\rangle ,\qquad k=0,\dots,M_{\xi}.
\]
so the coefficients $\lambda_{k}$\ are obtained as
\begin{equation}
\lambda_{k}=\sum_{j=0}^{M_{\xi}}\sum_{i=0}^{M_{\xi}}c_{ijk}\left\langle u_{i}%
,v_{j}\right\rangle _{\mathbb{R}},\qquad k=0,\dots,M_{\xi}, 
\label{eq:RQ} 
\end{equation}
where the notation$~\left\langle \cdot,\cdot\right\rangle _{\mathbb{R}}$
refers to the inner product of two vectors on Euclidean $M_{x}$-dimensional
space. It is interesting to note that~(\ref{eq:RQ}) is a Hadamard product,
see, e.g.,~\cite[Chapter~5]{Horn-1991-TMA}.

\begin{remark}
We used $M_{\lambda}=M_{\xi}$ in~(\ref{eq:global-system}), which is
determined by the definitions of eigenvalues and eigenvectors in~(\ref{eq:standard}),
and we used the same convention to compute the Rayleigh quotient~(\ref{eq:RQ}). 
It would be possible to compute~$\lambda_k$ 
for $k=M_{\xi}+1,\dots, M_{A}$ as well, since the inner product~$u^Tv$ 
of two eigenvectors which are expanded using chaos polynomials up to degree~$p$ 
has nonzero chaos coefficients up to degree~$2p$. 
Because $M_{\xi}<M_{A}$, this means that some terms are missing in the sum used to
construct the right-hand side of~(\ref{eq:RQ}).
An alternative to using this truncated sum is to use a {\em full} representation of
the Rayleigh quotient using the projection 
\[
\lambda_{k}
=\left\langle
u^{T}Au,\psi_{k}\right\rangle ,\qquad k=0,\dots,M_{\lambda}.
\]
In more detail, this
uses $M_{\lambda}=M_{A}$ and is given by
\[
\left\langle \sum_{i=0}^{M_{\lambda}}\lambda_{i}\psi_{i}\left(  \xi\right)
,\psi_{k}\right\rangle =\left\langle \left(  \sum_{i=0}^{M_{\xi}}u_{i}\psi
_{i}\left(  \xi\right)  \right)  ^{T}\left(  \sum_{\ell=0}^{M_{A}}A_{\ell}\psi_{\ell}\left(  \xi\right)  \right)  \left(  \sum_{j=0}^{M_{\xi}}v_{j}\psi_{j}\left(  \xi\right)  \right)  ,\psi_{k}\right\rangle ,
\]
where $k=0,\dots,M_{\lambda}$. So the coefficients $\lambda_{k}$\ are obtained
as
\begin{equation}
\lambda_{k}=\sum_{j=0}^{M_{\xi}}\sum_{i=0}^{M_{\xi}}\sum_{\ell=0}^{M_{A}}c_{\ell ijk
}\left(  u_{i}^{T}A_{\ell}u_{j}\right)  ,\qquad k=0,\dots,M_{\lambda
},\label{eq:RQ-full}\end{equation}
where
\begin{equation}
\label{eq:cijkl}
c_{\ell ijk}=\mathbb{E}\left[  \psi_{\ell} \psi_{i} \psi_{j} \psi_{k
}\right]  .
\end{equation}
We implemented and tested in numerical experiments both
computations~(\ref{eq:RQ}) and (\ref{eq:RQ-full}) 
and found the results to be virtually identical.
Note that~(\ref{eq:RQ-full}) is significantly more costly than~(\ref{eq:RQ}), so it
appears that there is no advantage to using~(\ref{eq:RQ-full}).
The construction~(\ref{eq:RQ}) appears to be new, 
but the truncated representation of $\lambda$ 
with $M_{\lambda}=M_{\xi}$ was also used 
in~\cite{Meidani-2014-SPI,Verhoosel-2006-ISR}.
\end{remark}

\paragraph{Normalization and the Gram-Schmidt process}

Let $~\left\Vert \cdot\right\Vert _{2}$ denote the norm induced by the inner
product $\left\langle \cdot,\cdot\right\rangle _{\mathbb{R}}$. That is, for a
vector$~u$ evaluated at a point$~\xi$, 
\begin{equation}
\left\Vert u\left(  \xi\right)  \right\Vert _{2}=\sqrt{\sum_{n=1}^{M_{x}%
}\left(  \left[  u\left(  \xi\right)  \right]  _{n}\right)  ^{2}}.
\label{eq:vector-norm} 
\end{equation}
We adopt the strategy used in$~$\cite{Meidani-2014-SPI}, 
whereby at each step of the stochastic iteration, the coefficients of the gPC
expansions of a given set of vectors$~\left\{  v^{s}\right\}  _{s=1}^{n_{s}}$
are transformed into an orthonormal set$~\left\{  u^{s}\right\}  _{s=1}%
^{n_{s}}$ such that
\begin{equation}
\left\langle u^{s}\left(  \xi\right)  ,u^{t}\left(  \xi\right)  \right\rangle
_{\mathbb{R}}=\delta_{st},\qquad\text{a.s}. 
\label{eq:orthonormal} 
\end{equation}
The condition~(\ref{eq:orthonormal}) is quite strict. However, because we
assume the eigenvectors have the form of stochastic polynomials that can be
easily sampled, the coefficients of the orthonormal eigenvectors can be
calculated relatively inexpensively using a discrete projection and a
quadrature rule as in~(\ref{eq:Q-lambda-u}). Note that each step of the
stochastic iteration entails construction of the eigenvector approximations at
the set of collocation points and, in contrast to the stochastic collocation
method, no deterministic eigenvalue problems are solved.
We also note that an alternative approach to normalization, based on solution 
of a certain nonlinear system was recently proposed 
by Hakula et al.~\cite{Hakula-2015-AMS}.

First, let us consider {\em normalization} of a vector, so $s=1$. 
The coefficients of a normalized
vector$~u_{k}^{1}$, for $k=0,\dots,M_{\xi}$, are computed from the
coefficients$~v_{k}^{1}$ as
\begin{equation}
u_{k}^{1}=\sum_{q=1}^{N_{q}}\frac{v^{1}\left(  \xi^{\left(  q\right)
}\right)  }{\left\Vert v^{1}\left(  \xi^{\left(  q\right)  }\right)
\right\Vert _{2}}\,\psi_{k}\left(  \xi^{\left(  q\right)  }\right)
\,w^{\left(  q\right)  }. 
\label{eq:vector-normalize}
\end{equation}
Then for general $s$, 
the {\em orthonormalization}~(\ref{eq:orthonormal}) is achieved 
by a stochastic version of the modified Gram-Schmidt algorithm 
proposed by Meidani and
Ghanem~\cite{Meidani-2014-SPI}.  It is based on the standard deterministic
formula, see, e.g.~\cite[Algorithm~8.1]{Trefethen-1997-NLA}, 
\[
u^{s}=v^{s}-\sum_{t=1}^{s-1}\frac{\left\langle v^{s},u^{t}\right\rangle _{\mathbb{R}}}{\left\langle u^{t},u^{t}\right\rangle _{\mathbb{R}}}\,u^{t},\quad s=2,\dots,n_{s}.
\]
For brevity, let us write $\chi^{ts}=\left\langle v^{s},u^{t}\right\rangle _{\mathbb{R}}/\left\langle u^{t},u^{t}\right\rangle _{\mathbb{R}}u^{t}$, 
so the expression above becomes
\begin{equation}
u^{s}=v^{s}-\sum_{t=1}^{s-1}\chi^{ts},\quad s=2,\dots,n_{s}.
\label{eq:vector-orthogonalize} 
\end{equation}
The stochastic counterpart of~(\ref{eq:vector-orthogonalize}) is obtained by
the stochastic Galerkin projection 
\[
\left\langle u^{s},\psi_{k}\right\rangle =\left\langle v^{s},\psi
_{k}\right\rangle -\sum_{t=1}^{s-1}\left\langle \chi^{ts},\psi_{k}%
\right\rangle ,\qquad k=0,\dots,M_{\xi},\quad s=2,\dots,n_{s}.
\]
Then the coefficients$~u_{k}^{s}$ are
\[
u_{k}^{s}=v_{k}^{s}-\sum_{t=1}^{s-1}\chi_{k}^{ts},\qquad k=0,\dots,M_{\xi},\quad
s=2,\dots,n_{s},
\]
where $\chi_{k}^{ts}$ are computed using a discrete projection and a
quadrature rule as in (\ref{eq:Q-lambda-u}),
\[
\chi_{k}^{ts}=\sum_{q=1}^{N_{q}}\chi^{ts}\left(  \xi^{\left(  q\right)
}\right)  \,\psi_{k}\left(  \xi^{\left(  q\right)  }\right)  \,w^{\left(
q\right)  }.
\]

\paragraph{Error assessment}

Ideally, we would like to minimize 
\[
\left\Vert \int_{\Gamma}r^{s}\,d\mu\left(  \xi\right)  \right\Vert _{2},\qquad
s=1,\dots,n_{s},
\]
where
\[
r^{s}=A\left(  \xi\right)  u^{s}\left(  \xi\right)  -\lambda^{s}\left(
\xi\right)  u^{s}\left(  \xi\right)  ,\qquad s=1,\dots,n_{s}
\]
is the true residual.
However, we are limited by the gPC\ framework.
In particular, the algorithm only provides the coefficients of expansion 
\[
\widetilde{r}_{k}^{s}=\left\langle Au^{s}-\lambda^{s}u^{s},\psi_{k}\right\rangle ,\qquad
k=0,\dots,M_{\xi},\quad s=1,\dots,n_{s}
\]
of the residual, 
i.e., the vector corresponding to the difference of the left and right-hand sides of 
(\ref{eq:global-system}).
One could assess accuracy using Monte Carlo sampling of this residual by computing
\[
r^{s}\left(  \xi^{i}\right)  =A\left(  \xi^{i}\right)  u^{s}\left(  \xi
^{i}\right)  -\lambda^{s}\left(  \xi^{i}\right)  u^{s}\left(  \xi^{i}\right)
,\qquad i=1,\dots, N_{MC},\quad s=1,\dots,n_{s},
\]
possibly at each step of the stochastic iteration.  
A much less expensive computation %, which we perform in our experiments, 
is to use the expansion
coefficients directly as an error indicator.
In particular, we can monitor the norms of the terms 
of~$\widetilde{r}_{k}^{s}$ corresponding to expected value 
and variance of~$r^{s}$, 
\begin{equation}
\varepsilon_{0}^{s,\left(  it\right)  }=\left\Vert \widetilde{r}%
_{0}^{s,\left(  it\right)  }\right\Vert _{2},\qquad\varepsilon_{\sigma^{2}%
}^{s,\left(  it\right)  }=\left\Vert \sum_{k=1}^{M_{\xi}}\left(  \widetilde{r}%
_{k}^{s,\left(  it\right)  }\right)  ^{2}\right\Vert _{2},\qquad
s=1,\dots,n_{s}. 
\label{eq:eps}
\end{equation}
We can also monitor the difference of the coefficients in two consecutive
iterations
\begin{equation}
u_{\Delta}^{s,\left(  it\right)  }=\left\Vert \left[
\begin{array}
[c]{c}%
%TCIMACRO{\unit{u}}%
%BeginExpansion
\operatorname{u}%
%EndExpansion
_{0}^{s,\left(  it\right)  }\\
\vdots\\%
%TCIMACRO{\unit{u}}%
%BeginExpansion
\operatorname{u}%
%EndExpansion
_{M_{\xi}}^{s,\left(  it\right)  }%
\end{array}
\right]  -\left[
\begin{array}
[c]{c}%
%TCIMACRO{\unit{u}}%
%BeginExpansion
\operatorname{u}%
%EndExpansion
_{0}^{s,\left(  it-1\right)  }\\
\vdots\\%
%TCIMACRO{\unit{u}}%
%BeginExpansion
\operatorname{u}%
%EndExpansion
_{M_{\xi}}^{s,\left(  it-1\right)  }%
\end{array}
\right]  \right\Vert _{2},\qquad s=1,\dots,n_{s}. \label{eq:u_Delta}%
\end{equation}

\section{Numerical experiments}

\label{sec:numerical}In this section, we report on computations of estimates
of the probability density functions (pdf) of certain distributions. The plots
presented below that illustrate these were obtained using the \texttt{Matlab}
function \texttt{ksdensity}, which computes a distribution estimate from
samples. These samples were computed either directly by\ the Monte Carlo
method or by sampling the gPC expansions~(\ref{eq:gPC-lambda-u}) obtained from
stochastic inverse subspace iteration or stochastic collocation. In
particular, we report pdf estimates of eigenvalue distributions, and of the
$\ell^{2}$-norm of the eigenvector approximation errors
\begin{equation}
\varepsilon_{u}^{s}\left(  \xi^{\left(  i\right)  }\right)  =\frac{\left\Vert
u^{s}\left(  \xi^{\left(  i\right)  }\right)  -u_{MC}^{s}\left(  \xi^{\left(
i\right)  }\right)  \right\Vert _{2}}{\left\Vert u_{MC}^{s}\left(
\xi^{\left(  i\right)  }\right)  \right\Vert _{2}},\qquad i=1,\dots
,N_{MC},\quad s=1,\dots,n_{s}, \label{eq:eps_u}%
\end{equation}
where $u^{s}\left(  \xi^{\left(  i\right)  }\right)  $ are\ samples of
eigenvectors obtained from either stochastic inverse (subspace) iteration or
stochastic collocation. We also report the pdf estimates of the normalized
$\ell^{2}$-norms of the true residual distribution
\begin{equation}
\varepsilon_{r}^{s}\left(  \xi^{i}\right)  =\frac{\left\Vert r^{s}\left(
\xi^{i}\right)  \right\Vert _{2}}{\left\Vert A\left(  \xi^{i}\right)
\right\Vert _{2}},\qquad i=1,\dots N_{MC},\quad s=1,\dots,n_{s}.
\label{eq:true_res}%
\end{equation}

We have implemented the methods in \texttt{Matlab} and applied it to vibration
analysis of undamped structures, using the code
from~\cite{Ferreira-2009-MCF-1Ed}. For these models, the associated mean
problem gives rise to symmetric positive-definite matrices. For the
parametrized uncertain term in the problem definition, we take Young's
modulus, which is a proportionality constant relating strains and stresses in
Hooke's law, as 
\begin{equation}
E\left(  x,\xi\right)  =\sum_{\ell=0}^{M_{A}}E_{\ell}\left(  x\right)
\psi_{\ell}\left(  \xi\right)  \label{eq:E_gPC}%
\end{equation}
to be a truncated lognormal process transformed from an underlying Gaussian
random process using a procedure described in~\cite{Ghanem-1999-NGS}. That is,
$\psi_{\ell}\left(  \xi\right)  $, $\ell=0,\dots,M_{A}$,\ is a set of $N_{\xi
}$-dimensional products of univariate Hermite polynomials and, denoting the
coefficients of the Karhunen-Lo\`{e}ve expansion of the Gaussian process by
$g_{j}\left(  x\right)  $ and $\eta_{j}=\xi_{j}-g_{j}$, $j=1,\dots,m_{\xi}$, the
coefficients in expansion~(\ref{eq:E_gPC}) are computed as
\[
E_{\ell}\left(  x\right)  =\frac{\mathbb{E}\left[  \psi_{\ell}\left(
\eta\right)  \right]  }{\mathbb{E}\left[  \psi_{\ell}^{2}\left(  \xi\right)
\right]  }\exp\left[  g_{0}\left(  x\right)  +\frac{1}{2}\sum_{j=1}^{m_{\xi}%
}\left(  g_{j}\left(  x\right)  \right)  ^{2}\right]  .
\]
The covariance function of the Gaussian field was chosen to be%
\[
C\left(  x_{1},x_{2}\right)  =\sigma_{g}^{2}\exp\left(  -\frac{\left\Vert
x_{1}-x_{2}\right\Vert _{2}}{L_{corr}}\right)  ,
\]
where $L_{corr}$ is the correlation length of the random variables $\xi_{i}$,
$i=1,\dots,m_{\xi}$, and $\sigma_{g}$ is the standard deviation of the
Gaussian random field. Other parameters in the models were deterministic (see
below). Note that, according to~\cite{Matthies-2005-GML}, in order to
guarantee a complete representation of the lognormal process
by~(\ref{eq:E_gPC}), the degree of polynomial expansion of$~E\left(
x,\xi\right)  $ should be twice the degree of the expansion of the solution.
We follow the same strategy here. Denoting by$~p$ the degree of polynomial
expansions of $u\left(  x,\xi\right)  $\ and $\lambda\left(  x,\xi\right)  $,
the total numbers of the gPC\ polynomials are
see, e.g.,~\cite[p. 87]{Ghanem-1991-SFE} and~\cite[Section~5.2]{Xiu-2010-NMS}, 
\begin{equation}
M_{\xi}+1=\frac{\left(  m_{\xi}+p\right)  !}{m_{\xi}!p!},\quad M_A%
+1=\frac{\left(  m_{\xi}+2p\right)  !}{m_{\xi}!\left(  2p\right)  !}.
\label{eq:M} 
\end{equation}

Finite element spatial discretization leads to a generalized eigenvalue
problem of the form
\begin{equation}
\mathcal{K}\!\left(  \xi\right)  u=\lambda\mathcal{M}u, 
\label{eq:KM}
\end{equation}
where $\mathcal{K}\!\left(  \xi\right)  =
\sum_{\ell=0}^{M_{A}}K_{\ell} \psi_{\ell}\left(  \xi\right)  $ is 
the stochastic stiffness matrix given by the gPC\ expansion, 
and $\mathcal{M}$ is the deterministic mass matrix.
Although we can transform~(\ref{eq:KM}) into a standard eigenvalue problem
$\mathcal{M}^{-1}\mathcal{K}\!\left(  \xi\right)  u=\lambda u,$ we found that
the stochastic Rayleigh quotient is sensitive to the nonsymmetry of this
matrix operator. 
We note that this is well known in the deterministic case and instead, 
two-sided Rayleigh quotients are often used~\cite{Hochstenbach-2003-TSA}. 
Here, we used for simplicity the Cholesky factorization $\mathcal{M}=LL^{T}$ 
and transformed~(\ref{eq:KM}) into
\begin{equation}
L^{-1}\mathcal{K}\!\left(  \xi\right)  L^{-T}w=\lambda w, 
\label{eq:KM2}
\end{equation}
where $u=L^{-T}w$. So, the expansion of $A$\ corresponding 
to~(\ref{eq:gPC-A}) is
\begin{equation}
A=\sum_{\ell=0}^{M_{A}}A_{\ell}\psi_{\ell}\left(  \xi\right)  =\sum_{\ell
=0}^{M_{A}}\left[  L^{-1}K_{\ell}\left(  \mathcal{\xi}\right)  L^{-T}\right]
\psi_{\ell}\left(  \xi\right)  . 
\label{eq:A_ell}
\end{equation}

We used the \texttt{Matlab} function \texttt{eig} to solve 
the deterministic eigenvalue problems: the mean value problem 
in Algorithm~\ref{alg:SISI} and at all sample points~$\xi^{\left( q \right)}$. 
We compared the results for
the stochastic Galerkin methods with ones obtained using Monte Carlo
simulation and stochastic collocation. The stochastic Galerkin methods include
stochastic inverse subspace iteration from Algorithm \ref{alg:SISI}, and
direct use of stochastic Rayleigh quotient~(\ref{eq:RQ}). The latter entails
solving the deterministic mean problem~(\ref{eq:SISI_mean}) by \texttt{eig}
and using (\ref{eq:SISI-u_i_1})--(\ref{eq:SISI-u_i_2}) for $u$ 
in~(\ref{eq:RQ}), i.e., the coefficients from~$\overline{u}$ are used for the
zero-order terms of the polynomial chaos basis and the coefficients of
higher-order terms are set to zero. The coefficients of~$v$ were obtained from
the matrix-vector product~(\ref{eq:mat-vec}). This construction of eigenvalues
will be denoted by~RQ$^{\left(  0\right)  }$ to indicate that no stochastic
iteration was performed. The stochastic dimension was $m_{\xi}=3$, degree of
the gPC expansion of the solution$~p=3$, and degree of the gPC expansion of
the lognormal process$~2p$. Unless stated otherwise, we used $5\times10^{4}$
samples for the Monte Carlo method, and a Smolyak\ sparse grid 
with Gauss-Hermite quadrature points and grid level $4$ for collocation. 
With these settings, the size of$~c_{\ell jk}$ in~(\ref{eq:cijk}) was
$84\times20\times20$ with $806$ nonzeros, 
the size of $c_{\ell ijk}$ in~(\ref{eq:cijkl}) was $84\times20\times20\times84$ with $103,084$ nonzeros, 
and there were $N_{q}=69$ points on the sparse grid.

\subsection{Example 1:\ Timoshenko beam}

\begin{figure}[ptbh]
\begin{center}
\includegraphics[width=8.5cm]{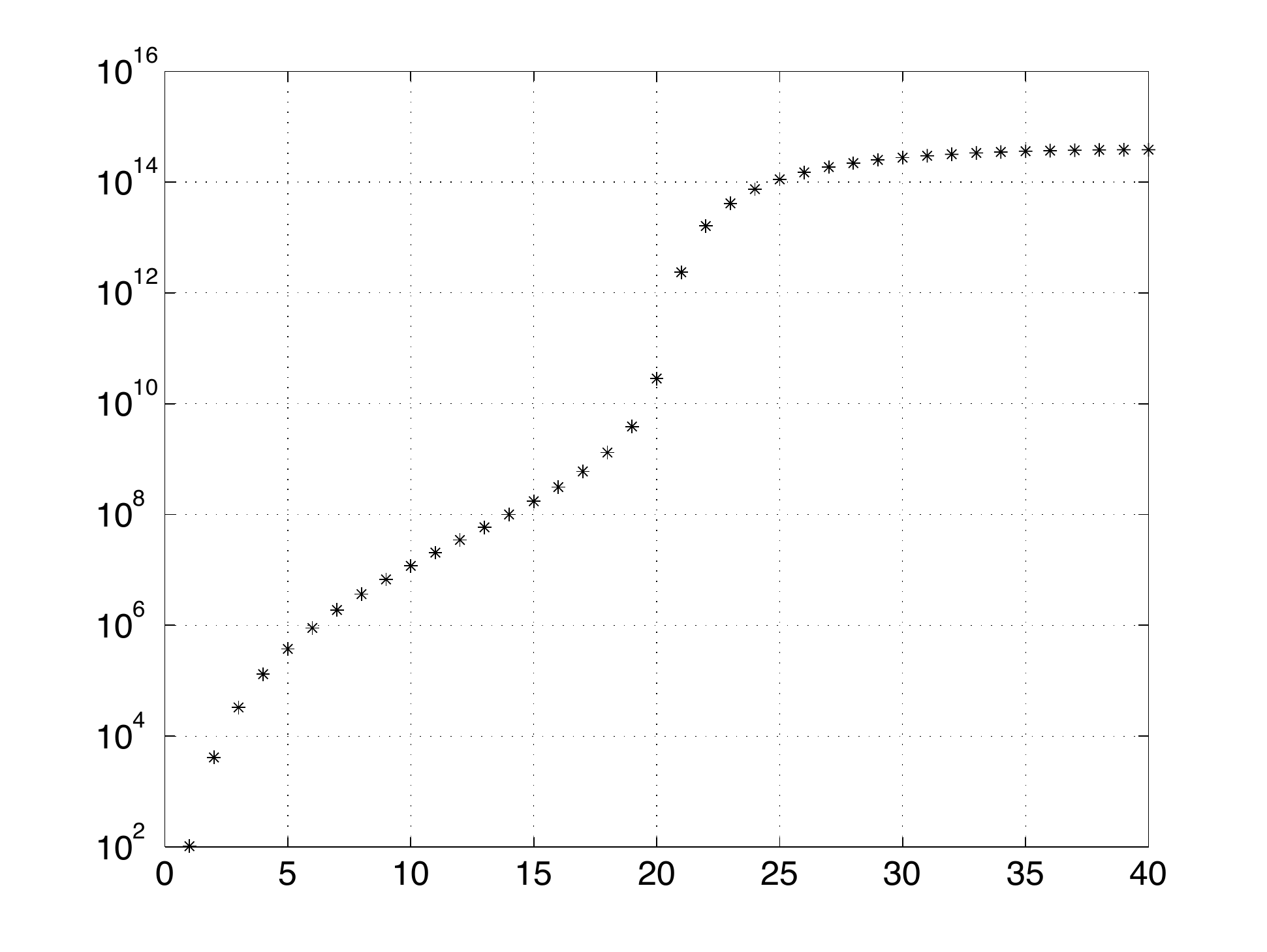}
\end{center}
\caption{Eigenvalues of the matrix $A_{0}$ corresponding to the Timoshenko
beam.} 
\label{fig:TB-eig-A_0} 
\end{figure}

For the first test problem, we analyzed free vibrations of a Timoshenko beam.
The kinetic energy of vibrations consists of two parts, one associated with
translations and one with rotations. The physical parameters of the cantilever
beam were set according to~\cite[Section 10.3]{Ferreira-2009-MCF-1Ed} 
as follows: the mean Young's modulus of the lognormal random field was
$E_{0}=10^{8}$, Poisson's ratio $\nu=0.30$, length$~L_{\text{beam}}=1$,
thickness$~0.001$, $\kappa=5/6$, and density $\rho=1$. The beam was
discretized using $20$ linear finite elements, i.e., with $40$ physical
degrees of freedom. The condition number of the mean matrix$~A_{0}$ from
(\ref{eq:A_ell})\ is$~3.7296\times10^{12}$, the norm $\left\Vert
A_{0}\right\Vert _{2}$ is $3.8442\times10^{14}$. 
The eigenvalues of$~A_{0}$\ are displayed in Figure~\ref{fig:TB-eig-A_0}. 
The correlation length was $L_{corr}=L_{\text{beam}}/4$, 
and the coefficient of variation $CoV$ of the
stochastic Young's modulus was set either to $0.1$ ($10\%$) or $0.25$
($25\%$), where $CoV=\sigma/E_{0}$, the ratio of the standard deviation and
the mean Young's modulus.

\begin{figure}[ptbh]
\begin{center}
\includegraphics[width=6.45cm]{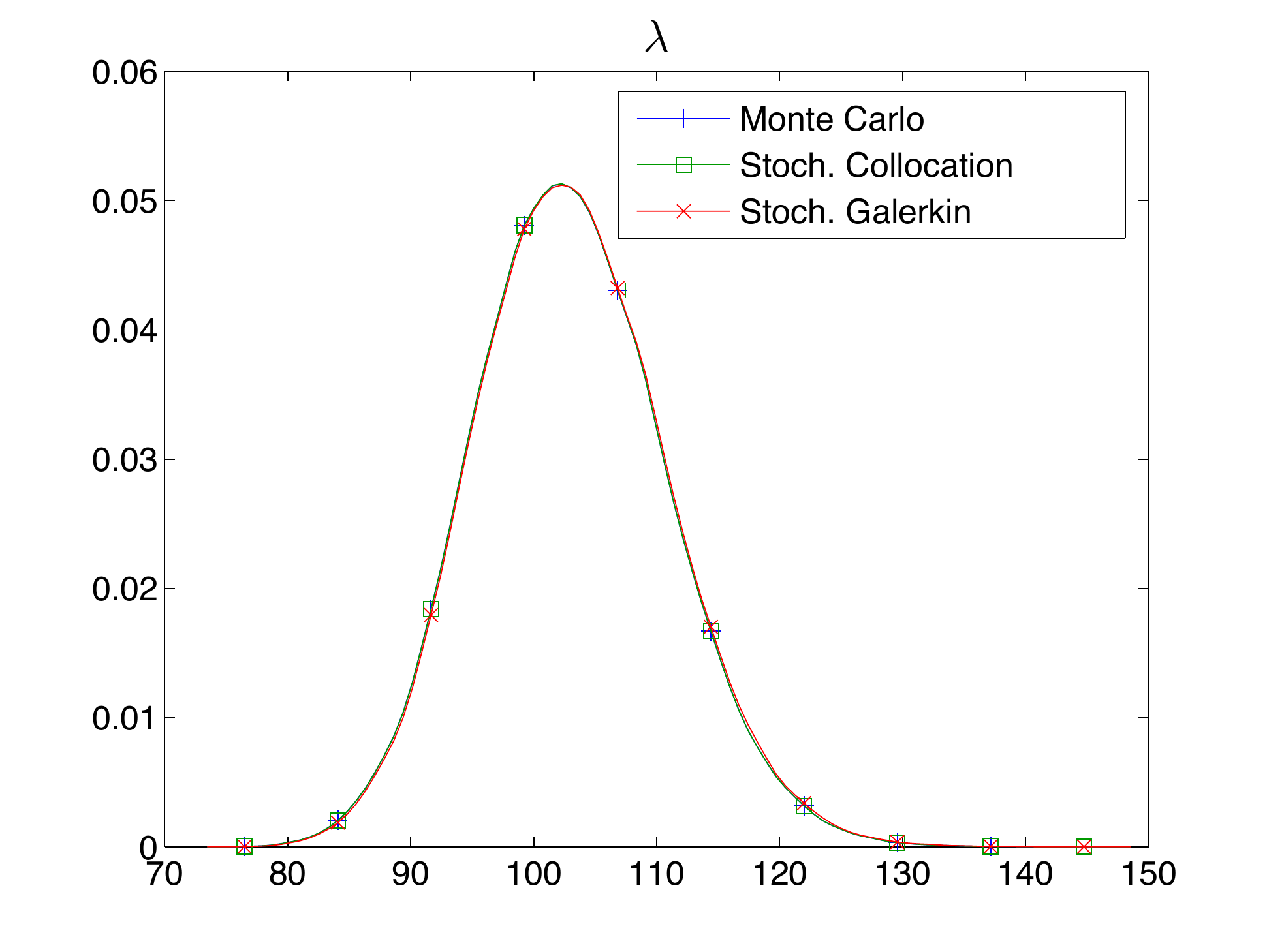}
\includegraphics[width=6.45cm]{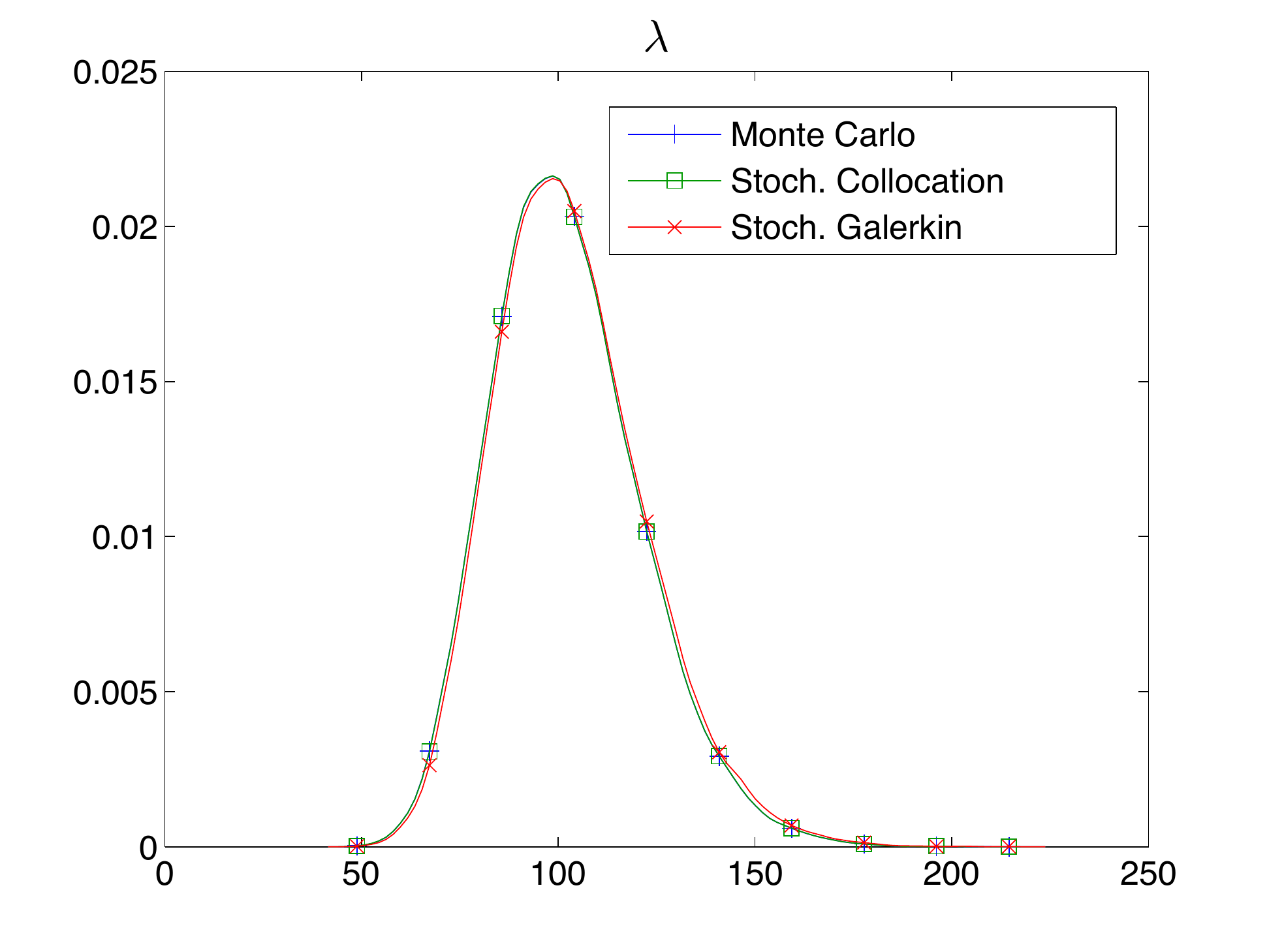}
\end{center}
\caption{Pdf estimates obtained from RQ$^{\left(  0\right)  }$, for the minimal
eigenvalue of the Timoshenko bean with $CoV=10\%$ (left) and $25\%$ (right).} 
\label{fig:TB-RQ} 
\end{figure}

First, we examine the performance of stochastic inverse iteration~(SII) and
compare it with stochastic collocation~(SC). We ran stochastic inverse
iteration with a fixed number of iterations, so plots of convergence
indicators~(\ref{eq:eps})--(\ref{eq:true_res}) shown below just illustrate the
performance of the algorithms. We computed estimates of pdfs for the
distributions of the eigenvalues and of the $\ell^{2}$-norm of the relative
eigenvector error (\ref{eq:eps_u}) corresponding to the minimal eigenvalue of
the Timoshenko beam, with $CoV=10\%$ and 25\%. Figure \ref{fig:TB-RQ} shows
the estimated eigenvalue distributions obtained using the \textquotedblleft
zero-step\textquotedblright\ computation~(RQ$^{\left(  0\right)  }$), which
uses only the mean solution~(\ref{eq:SISI-u_i_1})--(\ref{eq:SISI-u_i_2}). The
figure compares these distributions with those obtained using Monte Carlo and
stochastic collocation, and it is evident that the visible displays of the
three distributions are virtually indistinguishable. (Analogous plots, not
shown, obtained after one complete stochastic iteration produced essentially
identical plots.) As expected, the pdf estimates are narrower for $CoV=10\%$.
This computation is explored further in Tables~\ref{tab:TB-gPC-CoV10}
and~\ref{tab:TB-gPC-CoV25}, which show the first ten coefficients of the gPC
expansion of the smallest eigenvalue obtained using RQ$^{\left(  0\right)  }$,
one step and 20~steps of stochastic inverse iteration, and stochastic
collocation. It can be seen that RQ$^{(0)}$ provides good estimates of the
four coefficients corresponding to the mean ($d=0$) and linear terms ($d=1$)
of the expansion (\ref{eq:gPC-lambda-u}), and a single SII step significantly
improves the quality of the quadratic terms ($d=2$).\footnote{To test
robustness of the algorithms with respect to possible use of an inexact solver
of the deterministic mean value problem, we also examined perturbed initial
approximations $u_{0}^{s,\left(  0\right)  }=\overline{u}^{s}+\delta u^{s}$
for the stochastic iteration (\ref{eq:SISI-u_i_1}), where $\overline{u}^{s}$
is an eigenvector of the mean problem computed by \texttt{eig} and $\delta
u^{s}$ is a random perturbation with norm $10^{-6}$. We found this to have no
impact on performance in the sense that the columns for SII$^{(1)}$ and
SII$^{(20)}$ in Tables~\ref{tab:TB-gPC-CoV10}--\ref{tab:TB-gPC-CoV25} are
unchanged.}

\begin{table}[ptbh]
\caption{The first ten coefficients of the gPC expansion of the smallest
eigenvalue of the Timoshenko beam with $CoV=10\%$ using $0$, $1$ or $20$ steps
of stochastic inverse iteration, or using stochastic collocation. Here $d$ is
the polynomial degree and $k$ is the index of basis function 
in expansion~(\ref{eq:gPC-lambda-u}). } 
\label{tab:TB-gPC-CoV10} 
\begin{center}
\begin{tabular}
[c]{|c|c|r|r|r|r|}\hline
$d$ & $k$ & \multicolumn{1}{|c}{RQ$^{(0)}$} & \multicolumn{1}{|c|}{SII$^{(1)}%
$} & \multicolumn{1}{|c|}{SII$^{(20)}$} & \multicolumn{1}{c|}{SC}\\\hline
0 & 0 & 103.0823 & 102.9308 & 102.9307 & 102.9319\\\hline
\multirow{3}{*}{1} & 1 & 5.7301 & 5.7231 & 5.7231 & 5.7220\\
& 2 & -4.7970 & -4.7854 & -4.7854 & -4.7848\\
& 3 & 2.1156 & 2.1075 & 2.1075 & 2.1072\\\hline
\multirow{6}{*}{2} & 4 & 0.2361 & 0.2144 & 0.2144 & 0.2142\\
& 5 & -0.2540 & -0.2803 & -0.2804 & -0.2807\\
& 6 & 0.0841 & 0.1523 & 0.1523 & 0.1523\\
& 7 & 0.1873 & 0.1272 & 0.1271 & 0.1250\\
& 8 & -0.1437 & -0.0507 & -0.0506 & -0.0507\\
& 9 & 0.0961 & -0.0372 & -0.0373 & -0.0382\\\hline
\end{tabular}
\end{center}
\end{table}

\begin{table}[ptbh]
\caption{The first ten coefficients of the gPC expansion of the smallest
eigenvalue of the Timoshenko beam with $CoV=25\%$ using $0$, $1$ or $20$ steps
of stochastic inverse iteration, or using stochastic collocation. Here $d$ is
the polynomial degree and $k$ is the index of basis function 
in expansion~(\ref{eq:gPC-lambda-u}). } 
\label{tab:TB-gPC-CoV25}
\begin{center}
\begin{tabular}
[c]{|c|c|r|r|r|r|}\hline
$d$ & $k$ & \multicolumn{1}{|c}{RQ$^{(0)}$} & \multicolumn{1}{|c|}{SII$^{(1)}%
$} & \multicolumn{1}{|c|}{SII$^{(20)}$} & \multicolumn{1}{c|}{SC}\\\hline
0 & 0 & 103.0823 & 102.1705 & 102.1670 & 102.1713\\\hline
\multirow{3}{*}{1} & 1 & 14.0453 & 13.9402 & 13.9402 & 13.9408\\
& 2 & -11.7568 & -11.5862 & -11.5859 & -11.5848\\
& 3 & 5.1830 & 5.0654 & 5.0651 & 5.0669\\\hline
\multirow{6}{*}{2} & 4 & 1.4284 & 1.2919 & 1.2918 & 1.2909\\
& 5 & -1.5368 & -1.6766 & -1.6767 & -1.6764\\
& 6 & 0.5090 & 0.9030 & 0.9032 & 0.9035\\
& 7 & 1.1331 & 0.7533 & 0.7530 & 0.7529\\
& 8 & -0.8696 & -0.2965 & -0.2960 & -0.2955\\
& 9 & 0.5812 & -0.2215 & -0.2220 & -0.2222\\\hline
\end{tabular}
\end{center}
\end{table}

\begin{figure}[ptbh]
\begin{center}
\includegraphics[width=6.45cm]{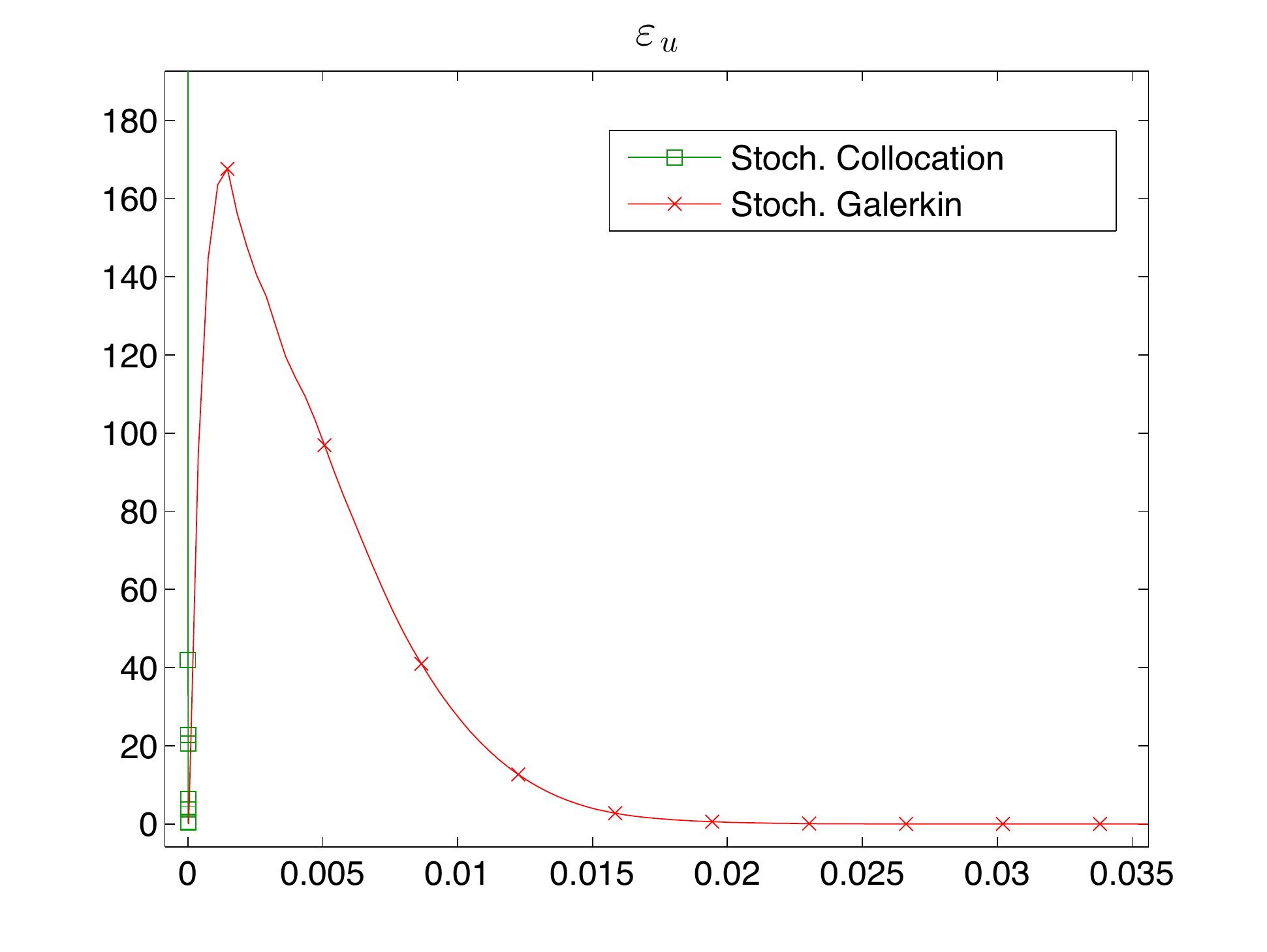}
\includegraphics[width=6.45cm]{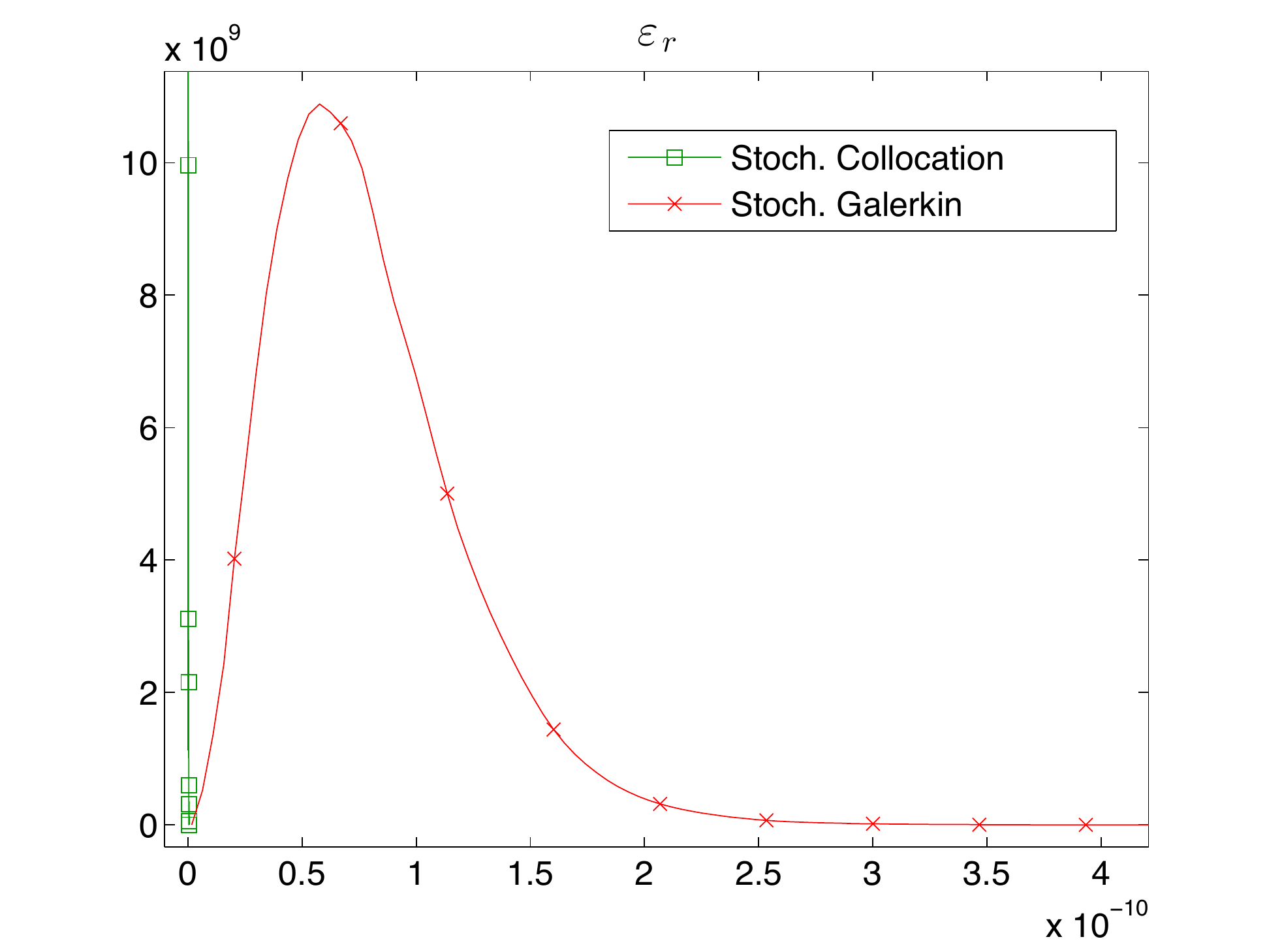}
\includegraphics[width=6.45cm]{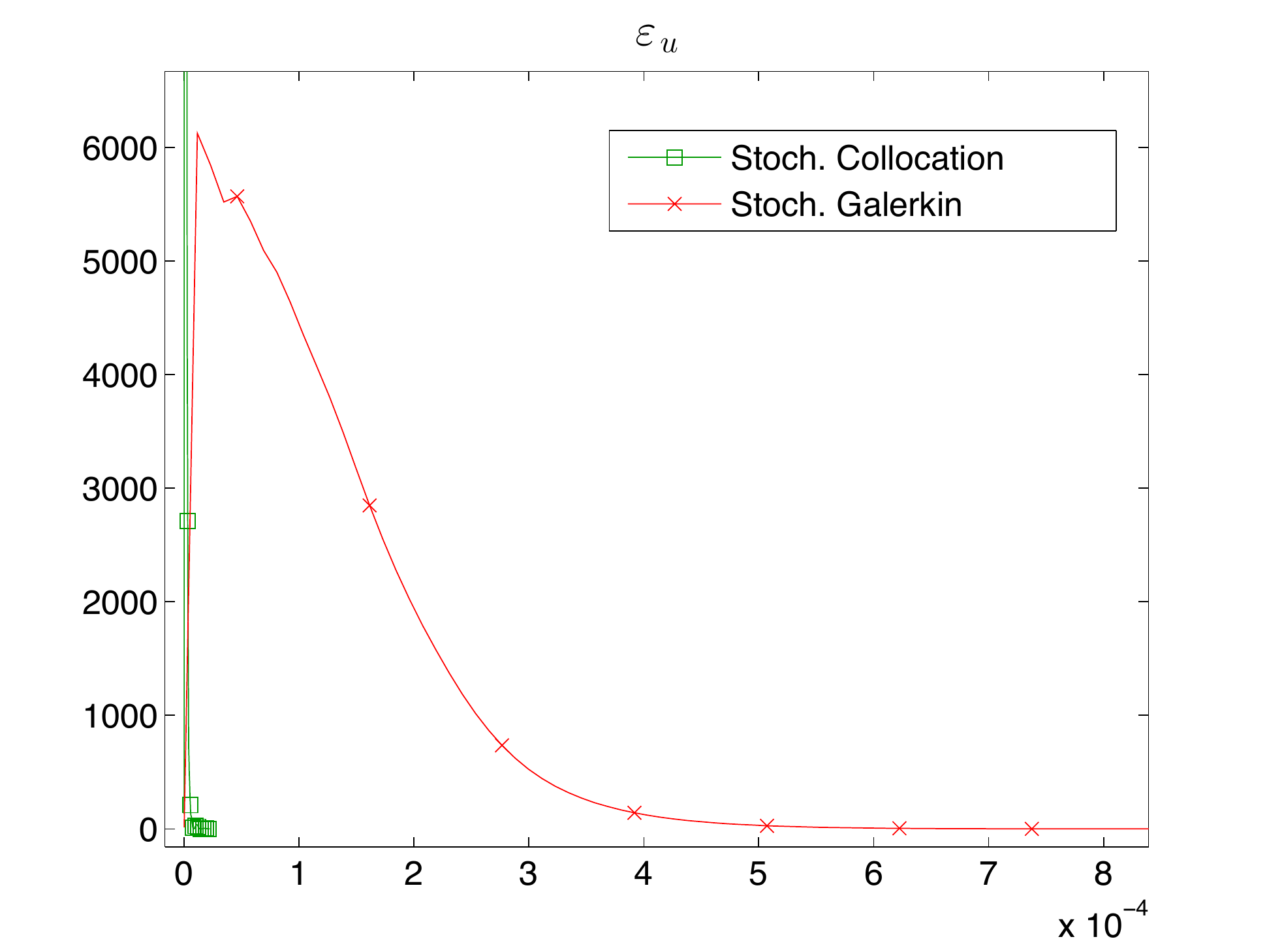}
\includegraphics[width=6.45cm]{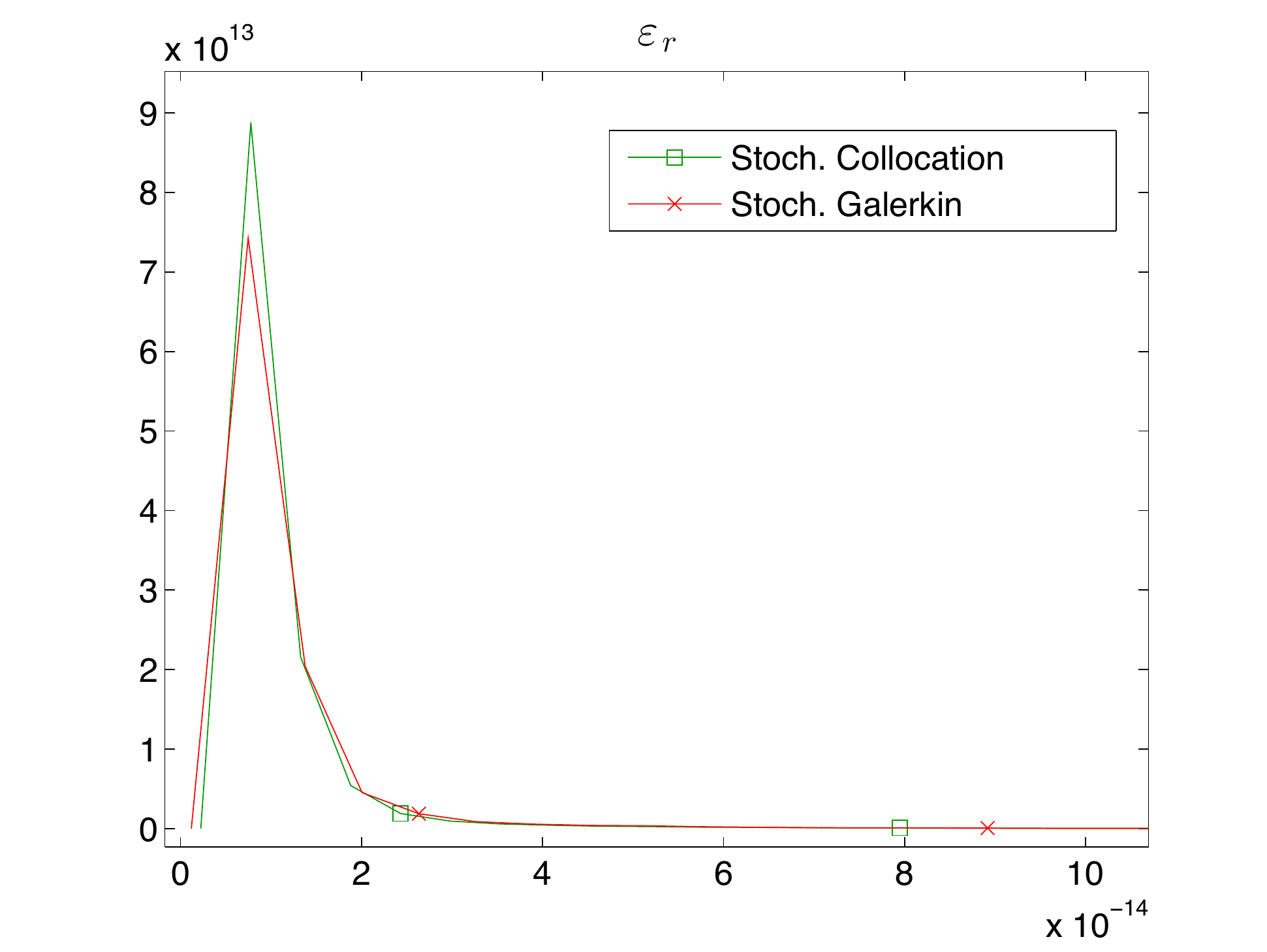}
\includegraphics[width=6.45cm]{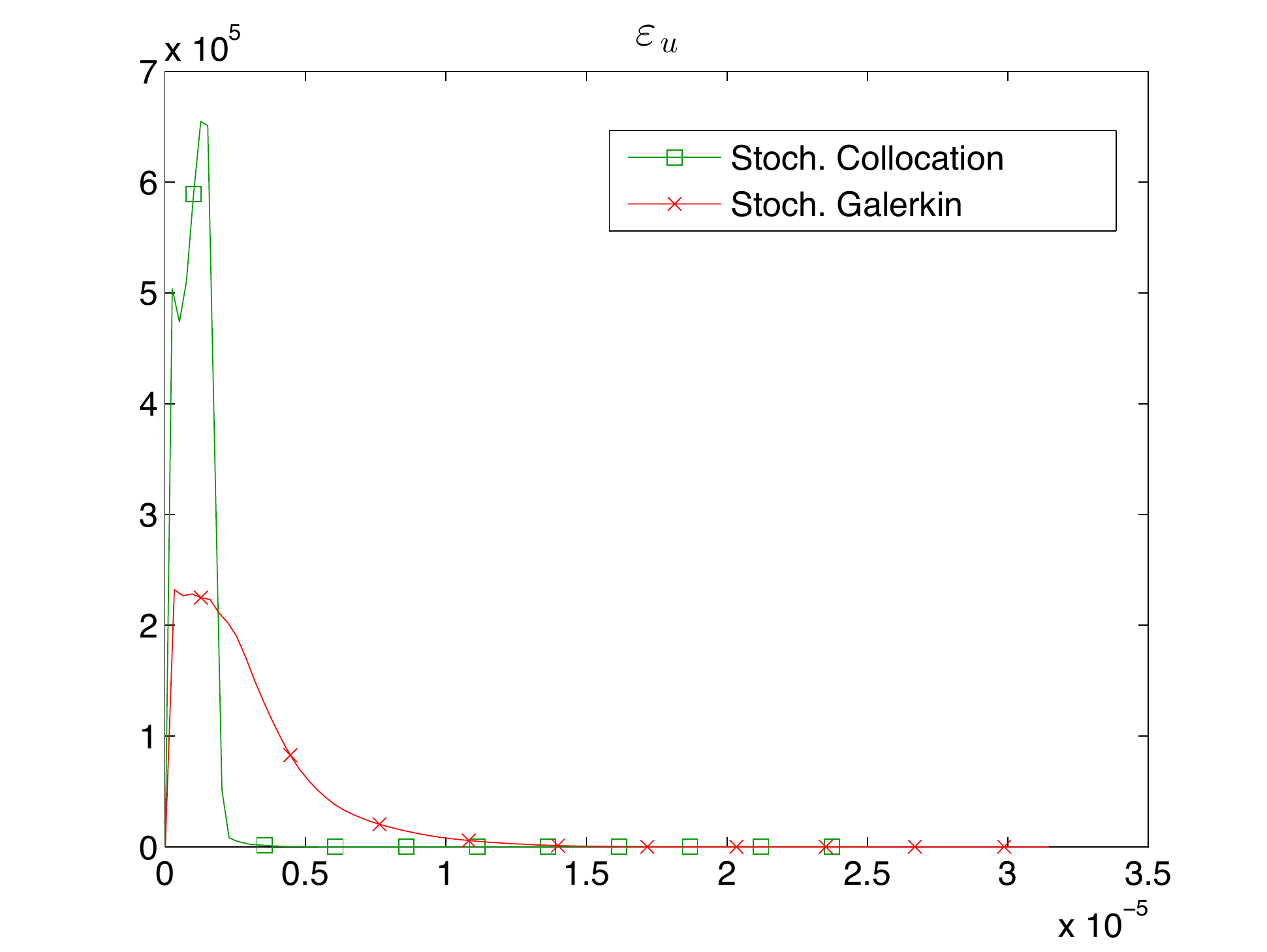}
\includegraphics[width=6.45cm]{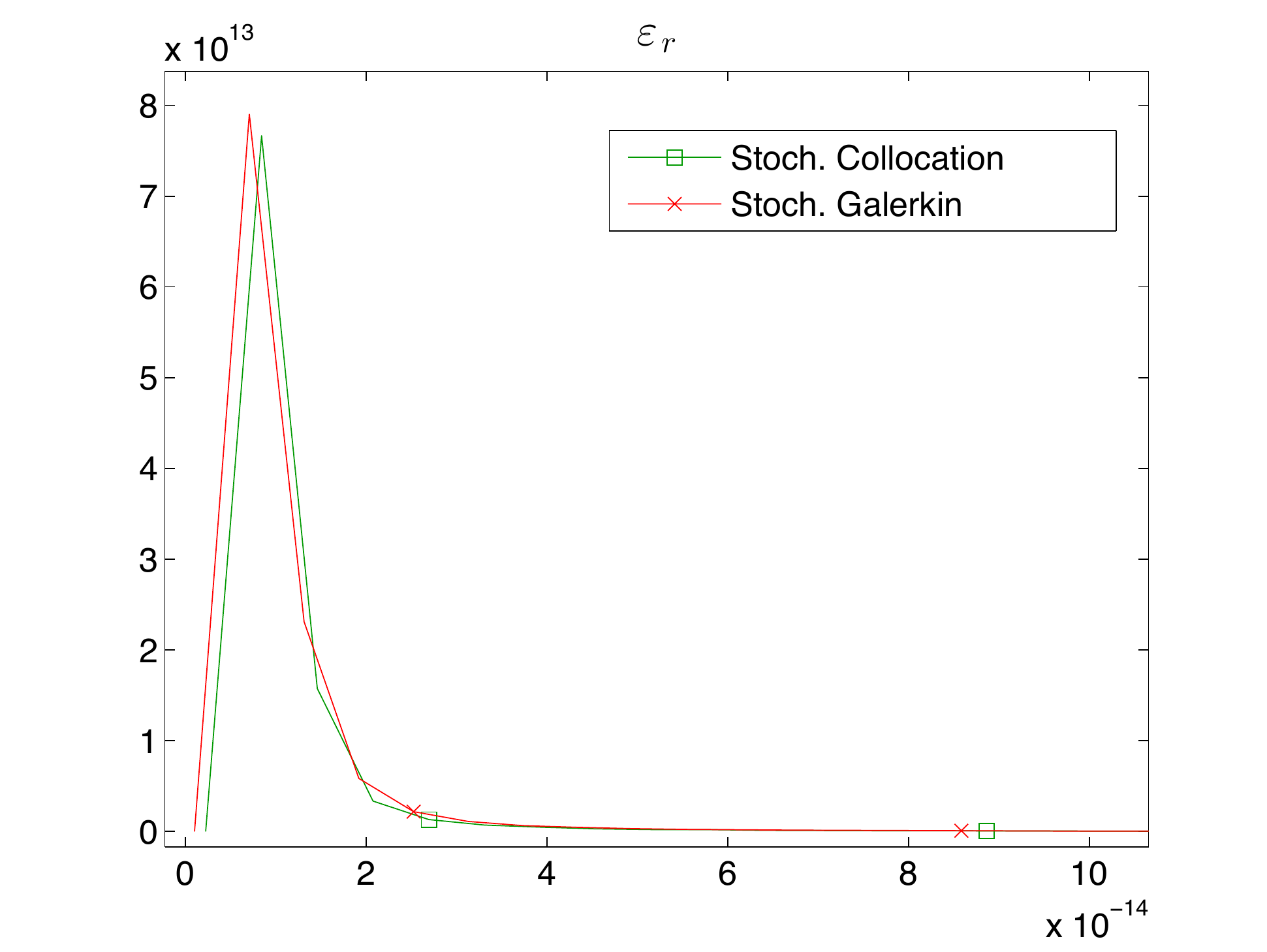}
\end{center}
\caption{Plots of the pdf estimate of the $\ell^{2}$-norms of the relative
eigenvector error (\ref{eq:eps_u}) (left) and the residual (\ref{eq:true_res})
(right) corresponding to the smallest eigenvalue of the Timoshenko beam with
$CoV=10\%$ obtained using stochastic Rayleigh quotient RQ$^{(0)}$ (top), and
after stochastic inverse iteration $1$ (middle) and $2$ (bottom).} 
\label{fig:TB-it-CoV10} 
\end{figure}

\begin{figure}[ptbh]
\begin{center}
\includegraphics[width=6.45cm]{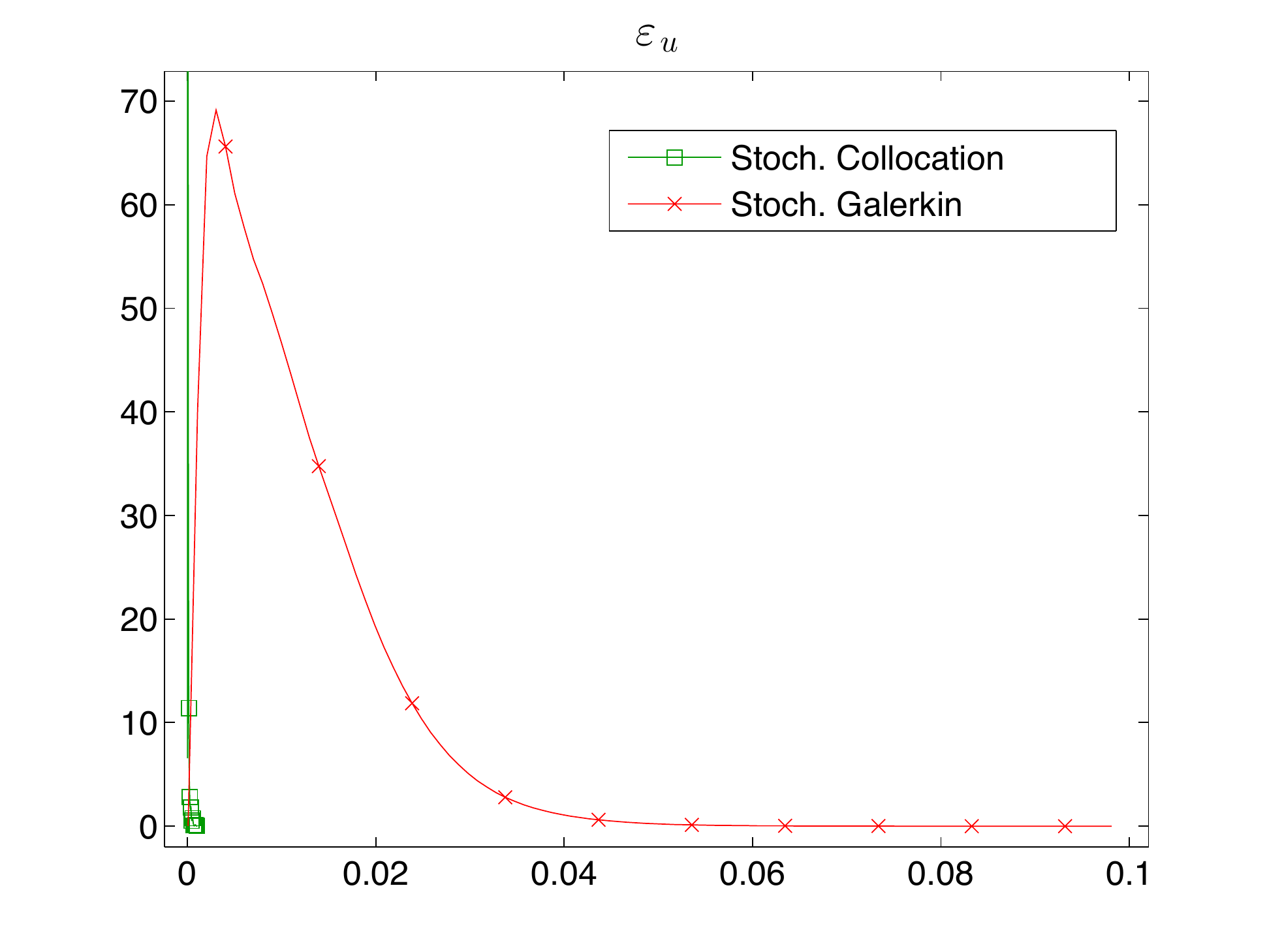}
\includegraphics[width=6.45cm]{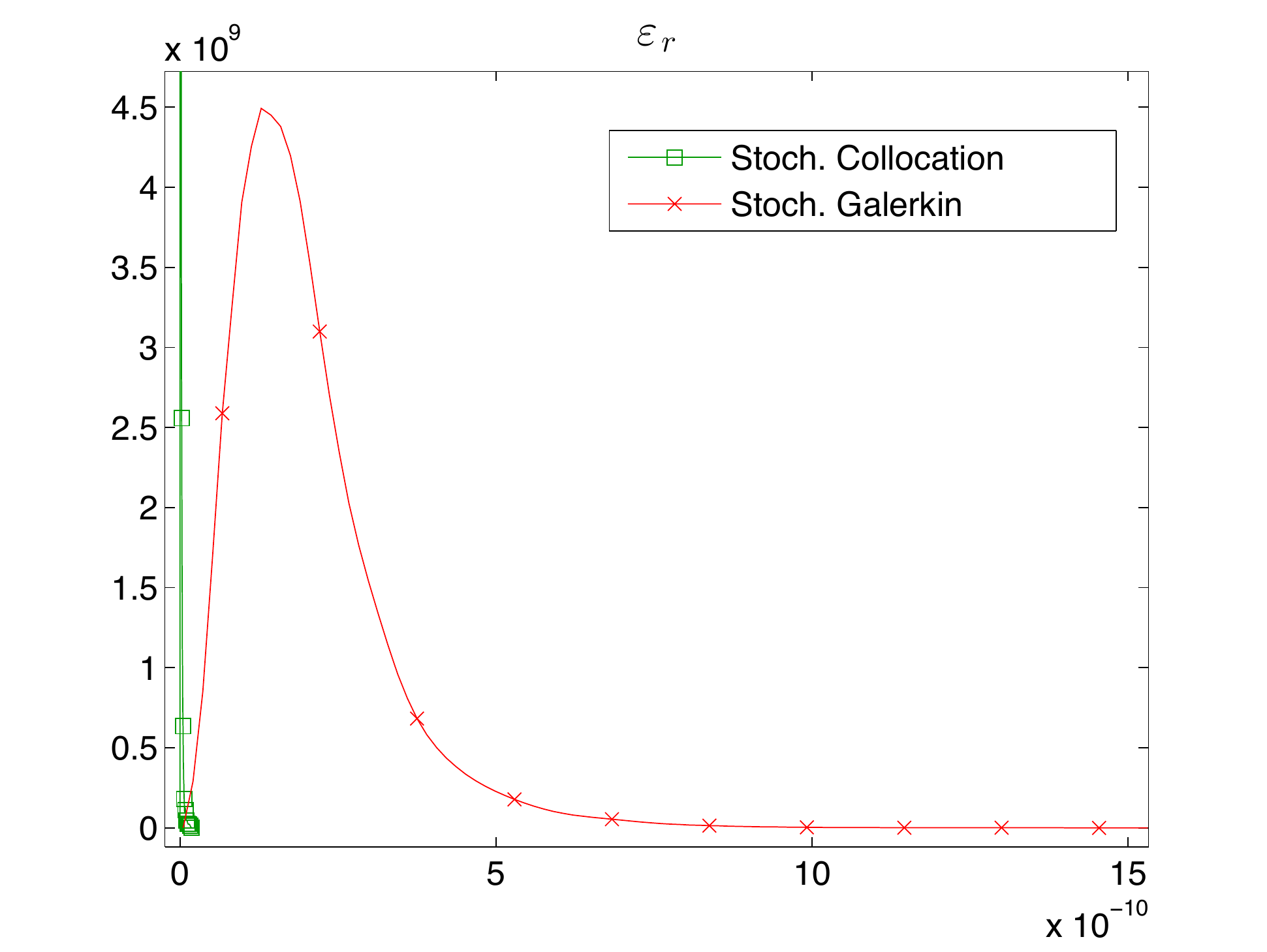}
\includegraphics[width=6.45cm]{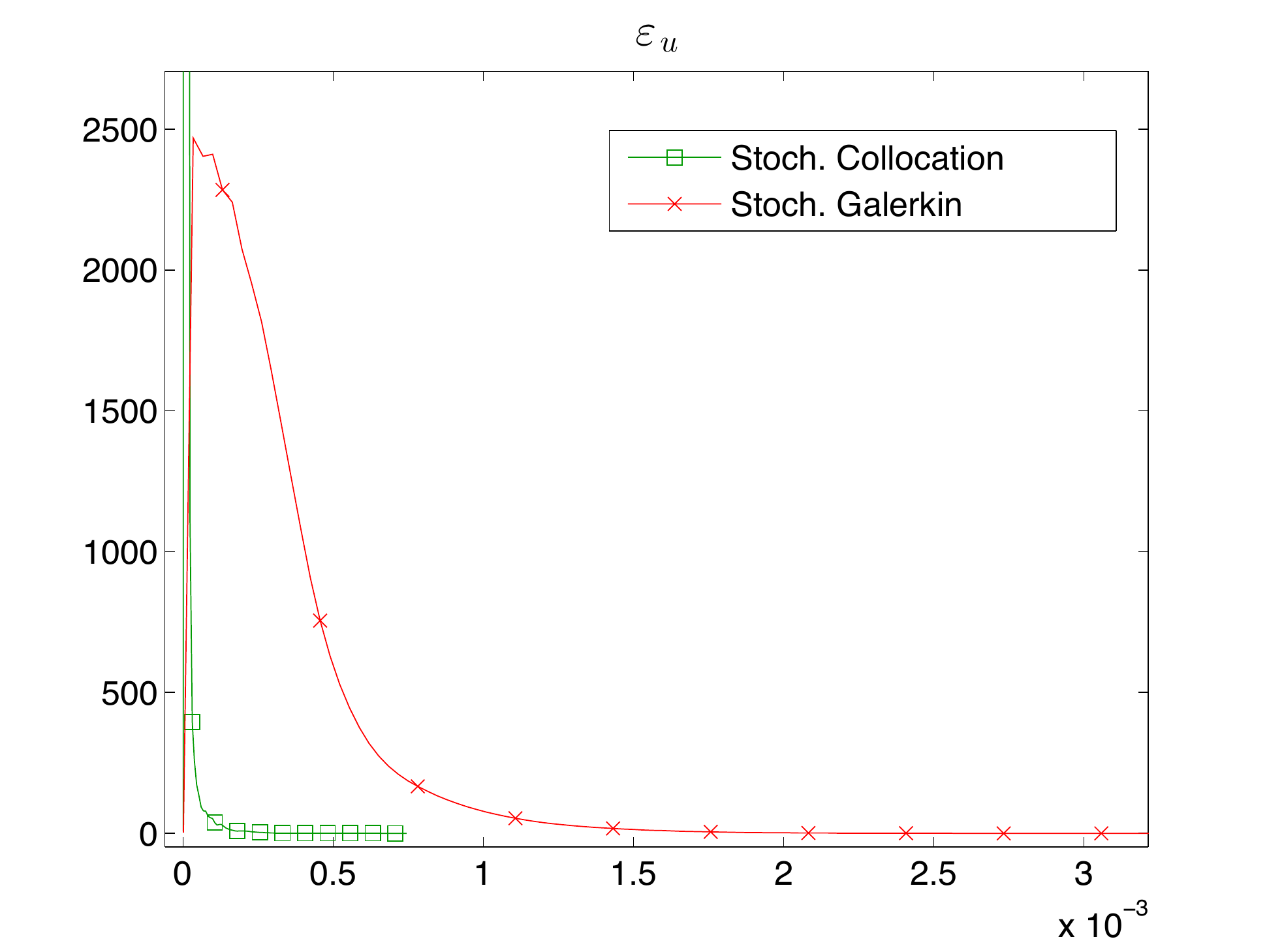}
\includegraphics[width=6.45cm]{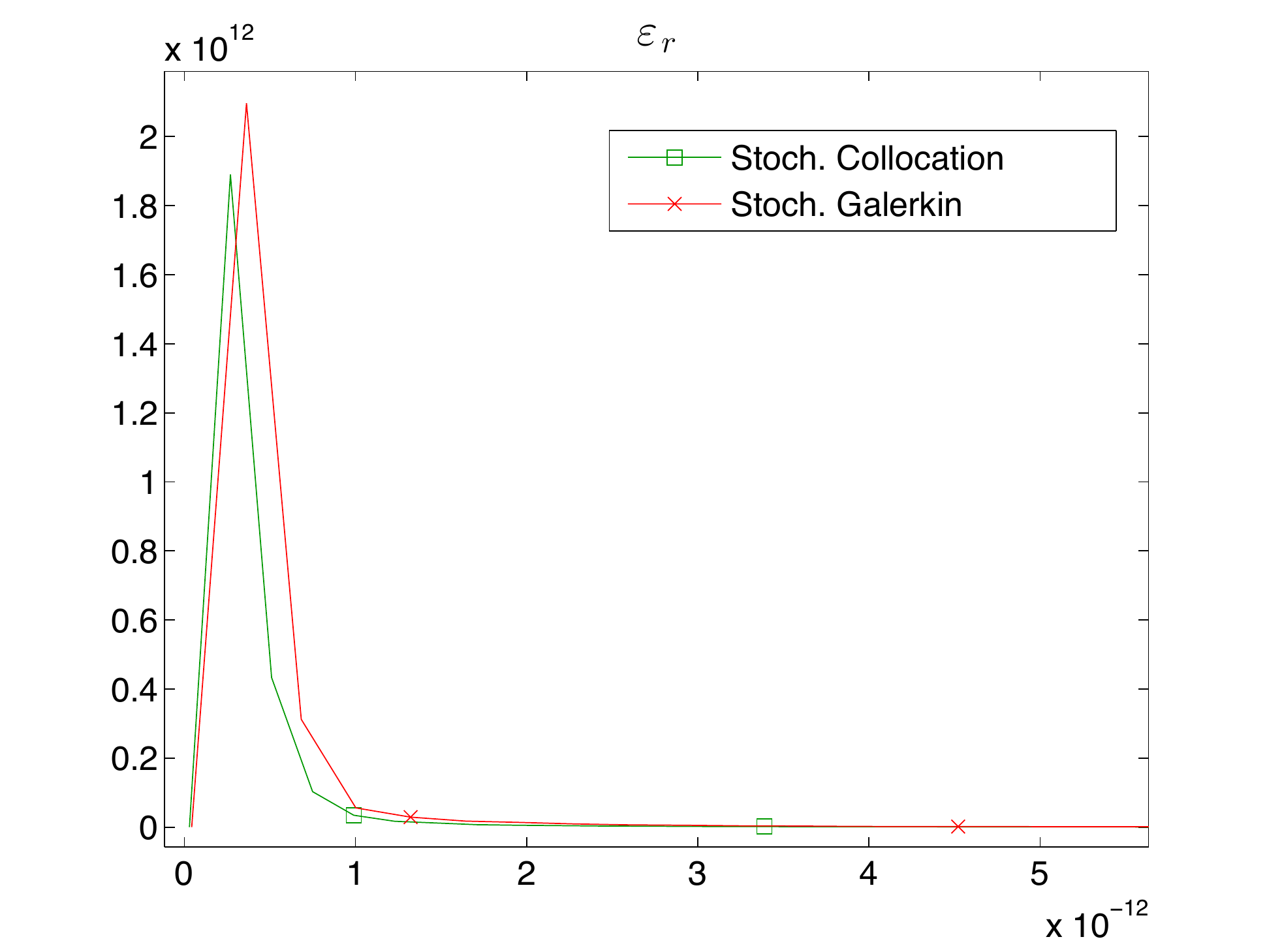}
\includegraphics[width=6.45cm]{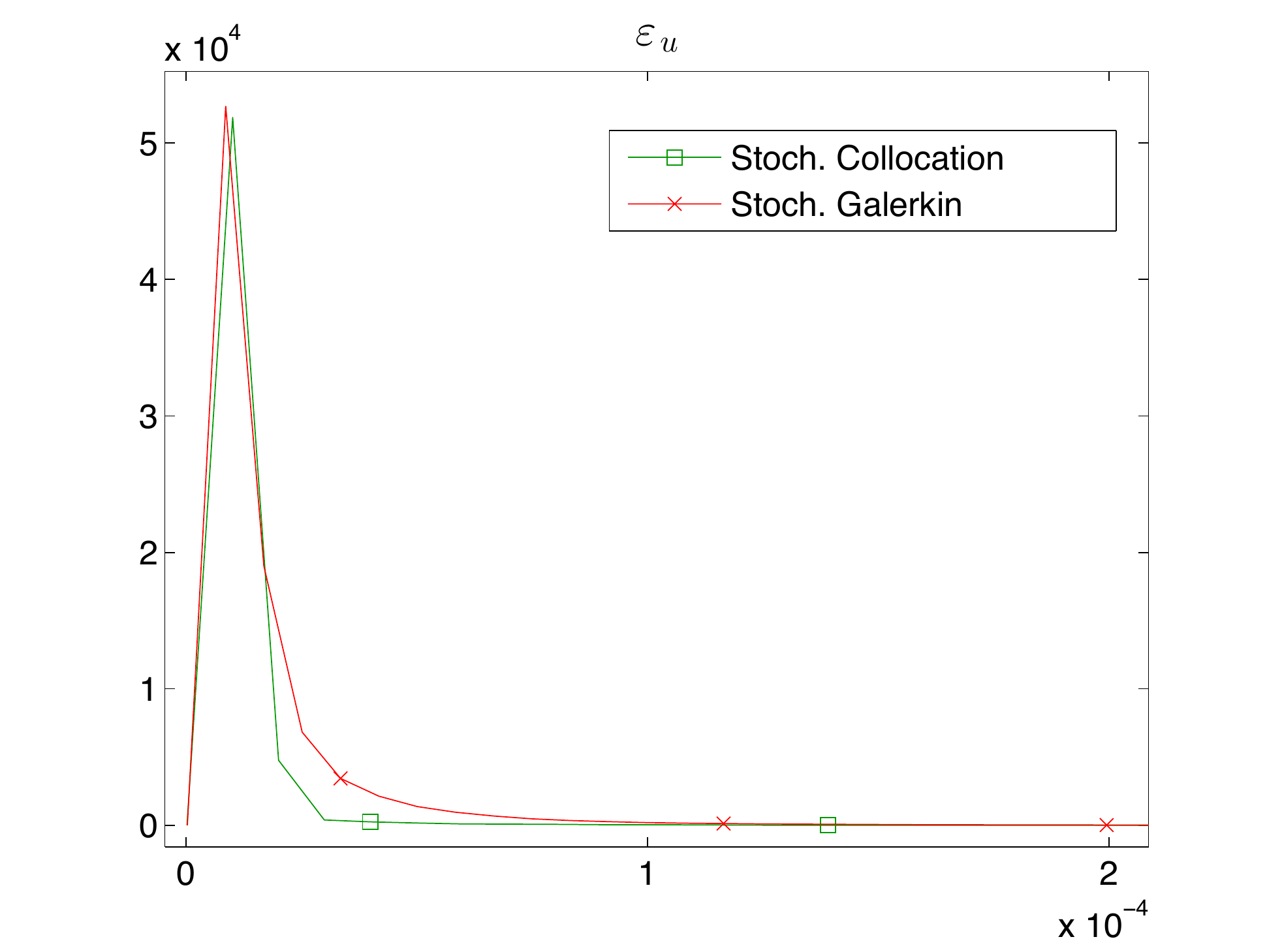}
\includegraphics[width=6.45cm]{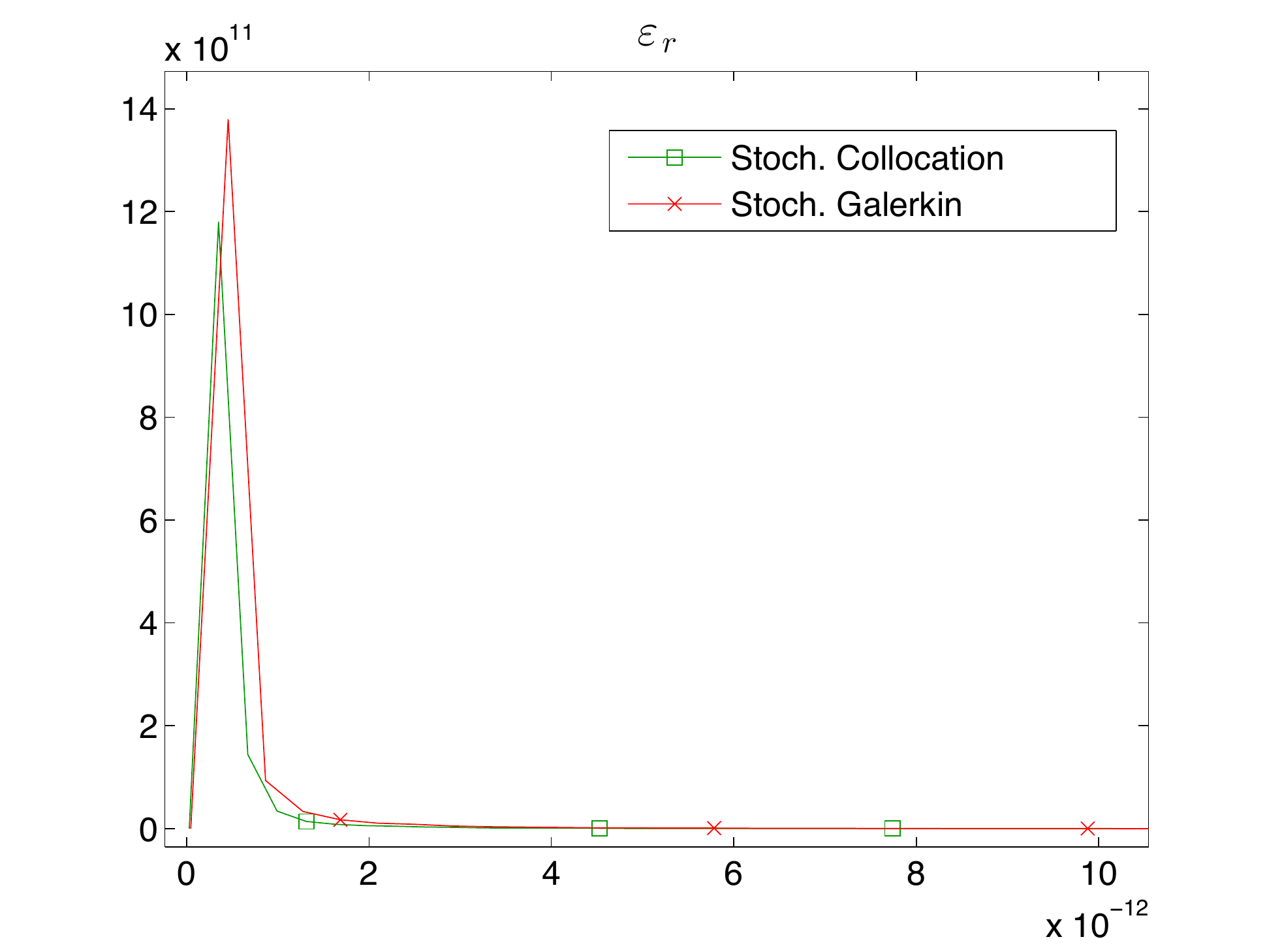}
\end{center}
\caption{Plots of the pdf estimate of the $\ell^{2}$-norms of the relative
eigenvector error (\ref{eq:eps_u}) (left) and the residual (\ref{eq:true_res})
(right) corresponding to the smallest eigenvalue of the Timoshenko beam with
$CoV=25\%$ obtained using stochastic Rayleigh quotient RQ$^{(0)}$ (top), and
after stochastic inverse iteration $1$ (middle) and $2$ (bottom).} 
\label{fig:TB-it-CoV25} 
\end{figure}

\begin{figure}[ptbh]
\begin{center}
\includegraphics[width=6.45cm]{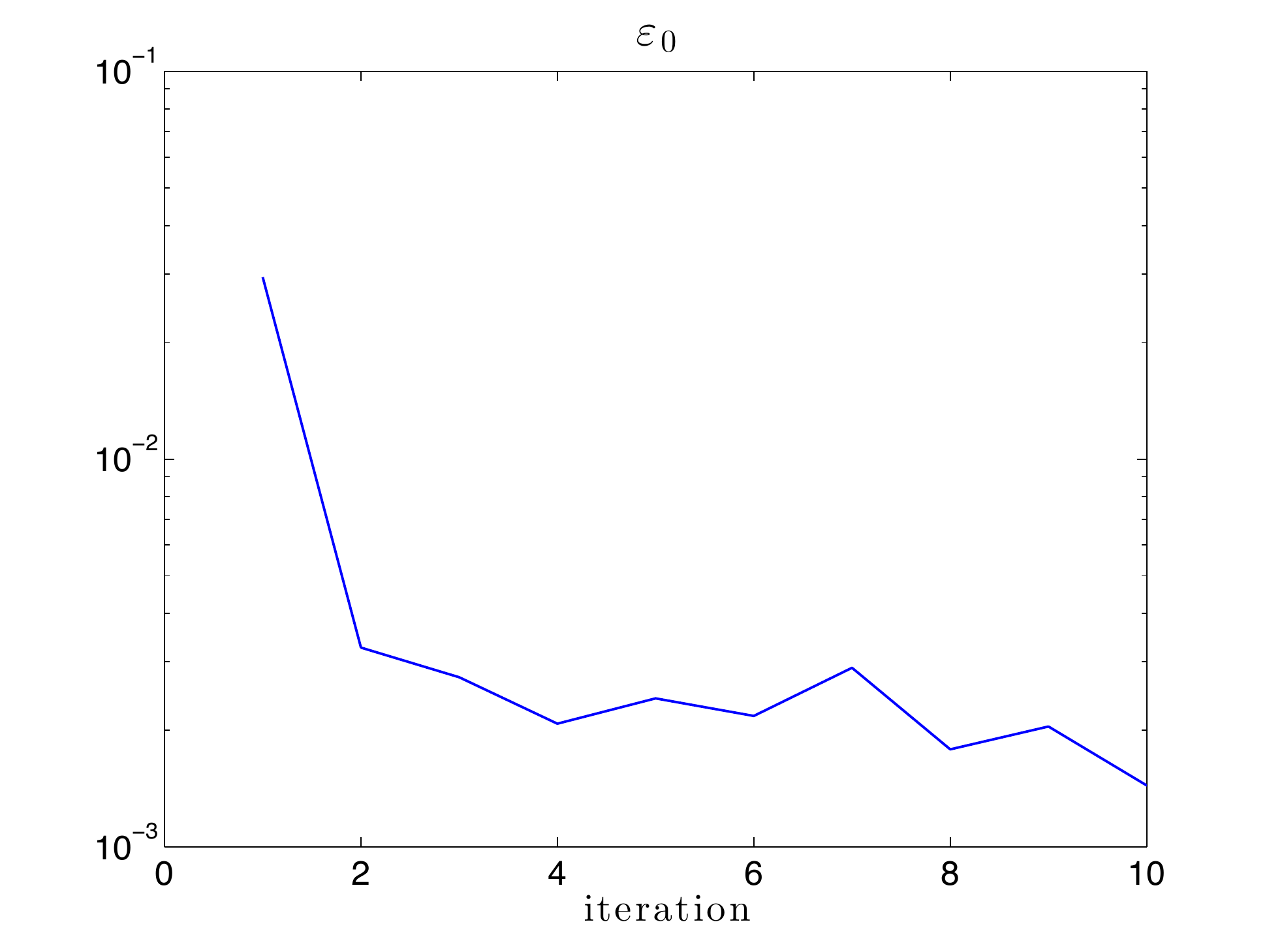}
\includegraphics[width=6.45cm]{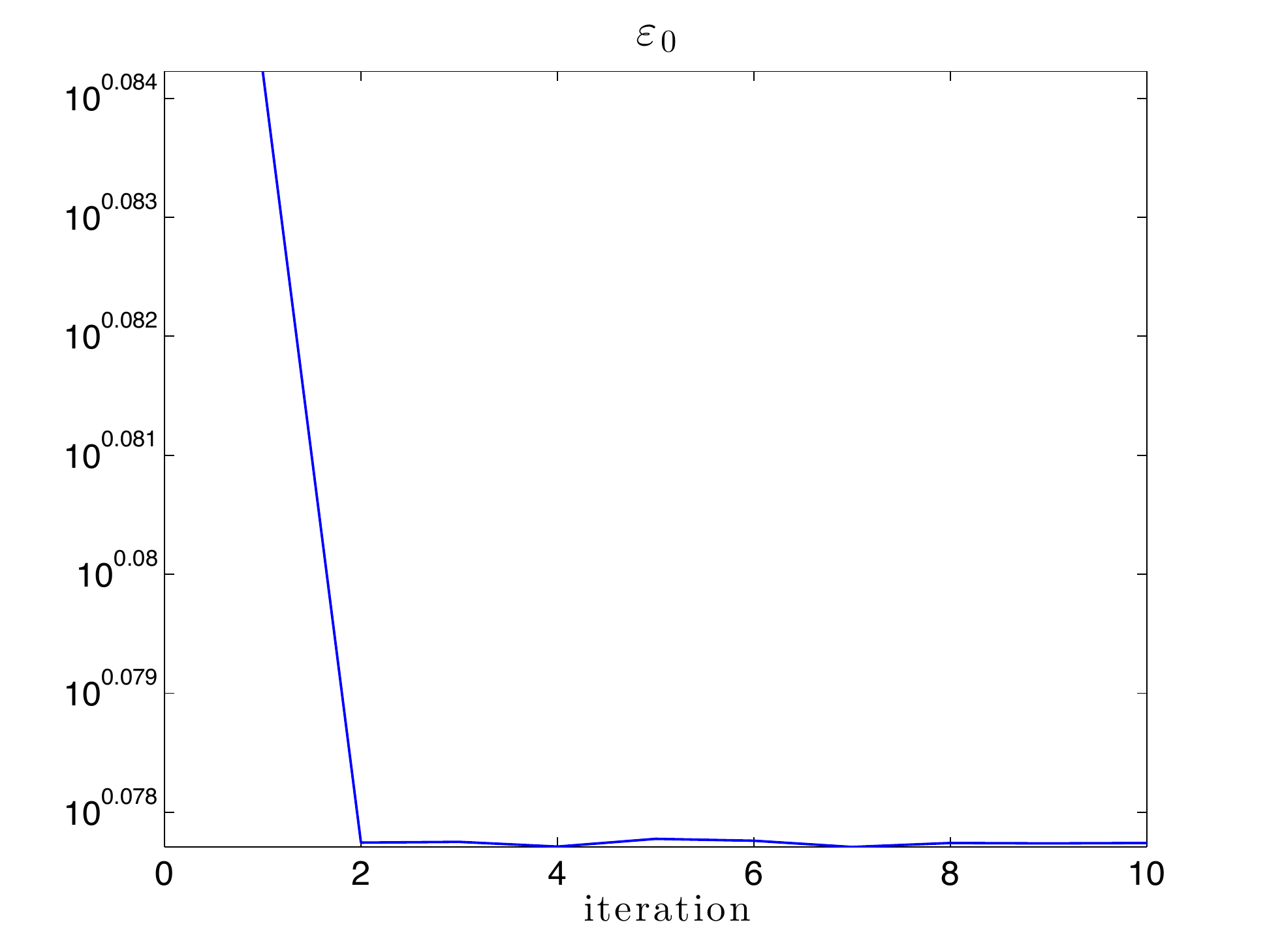}
\includegraphics[width=6.45cm]{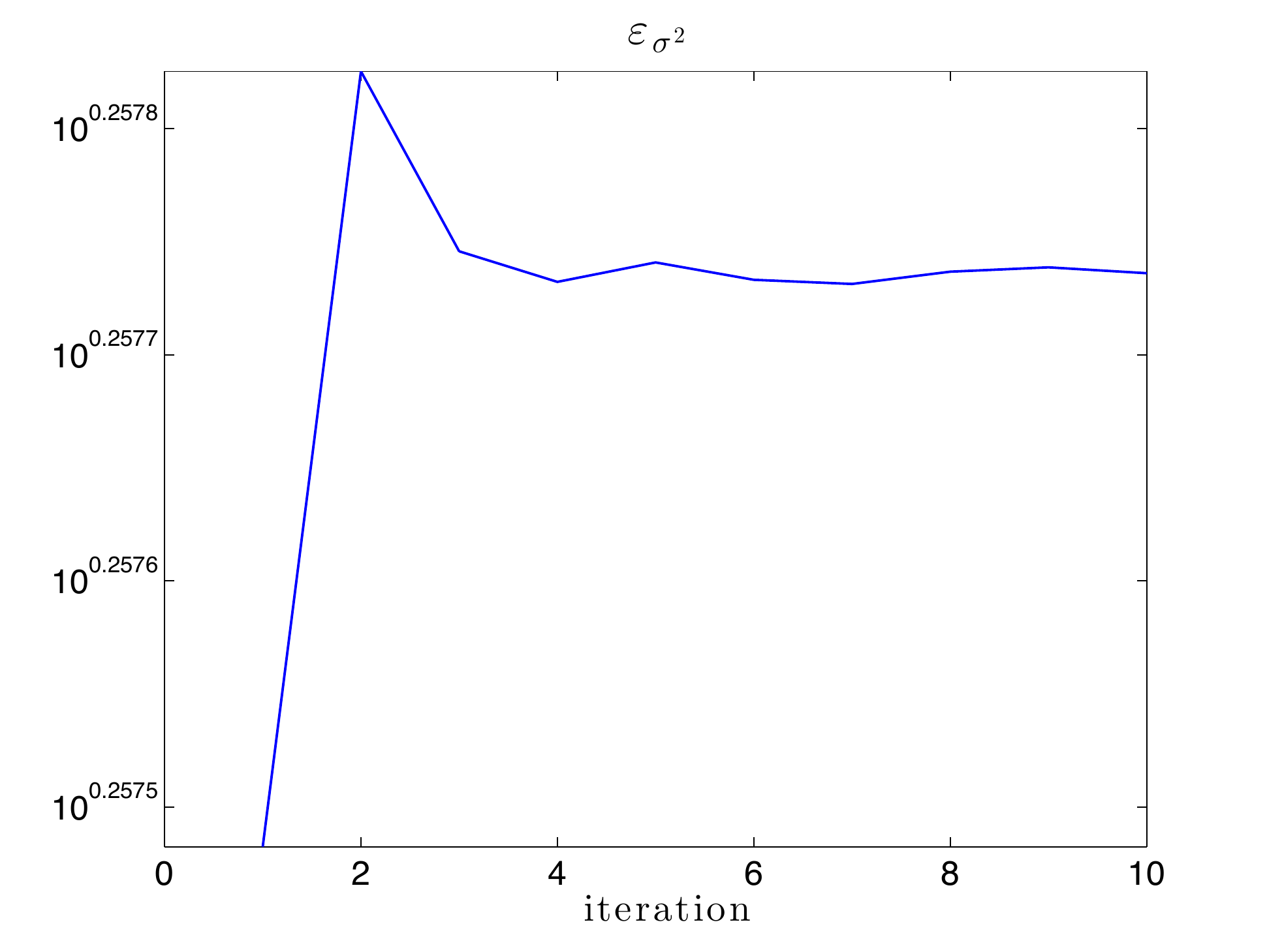}
\includegraphics[width=6.45cm]{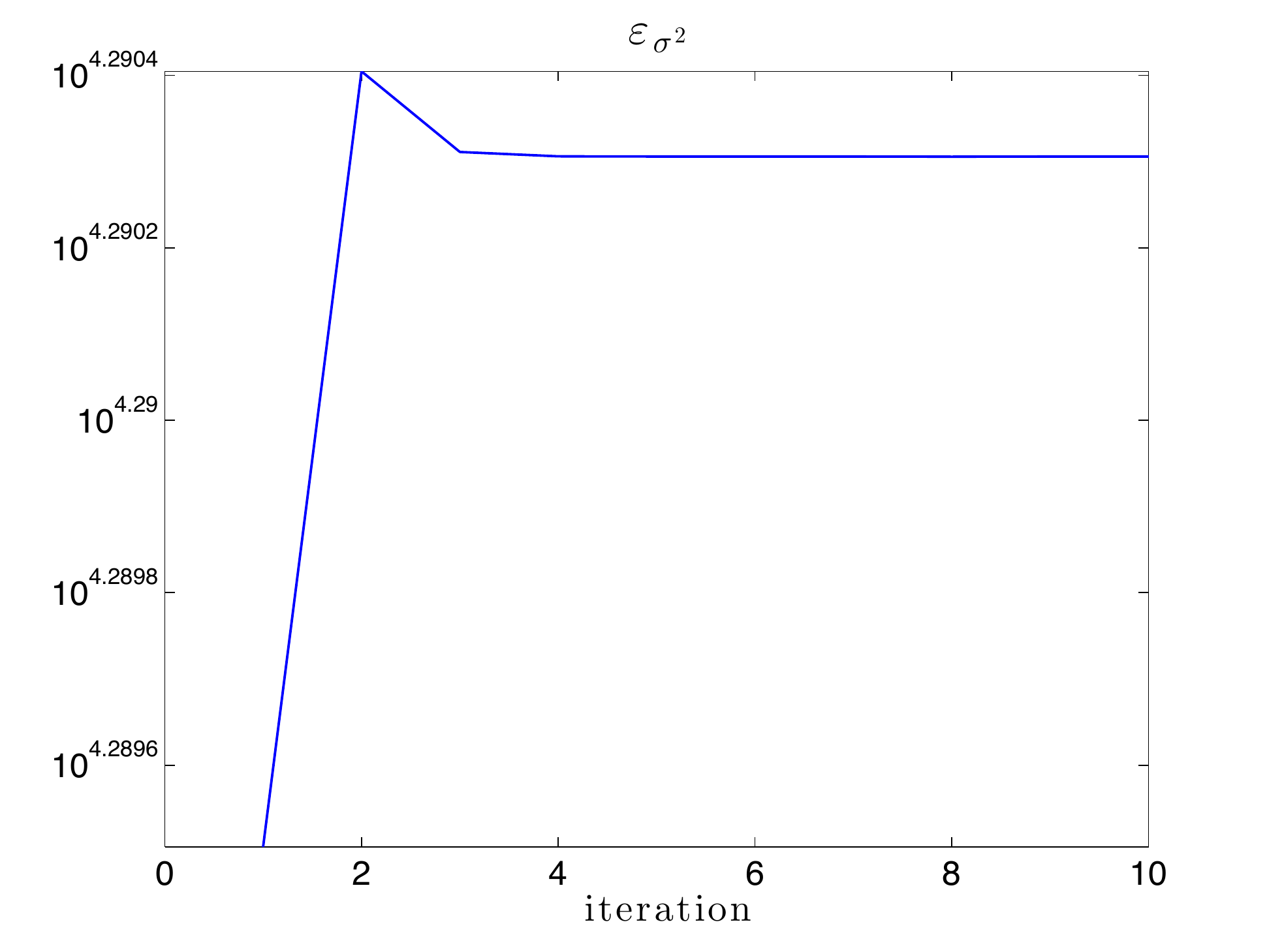}
\includegraphics[width=6.45cm]{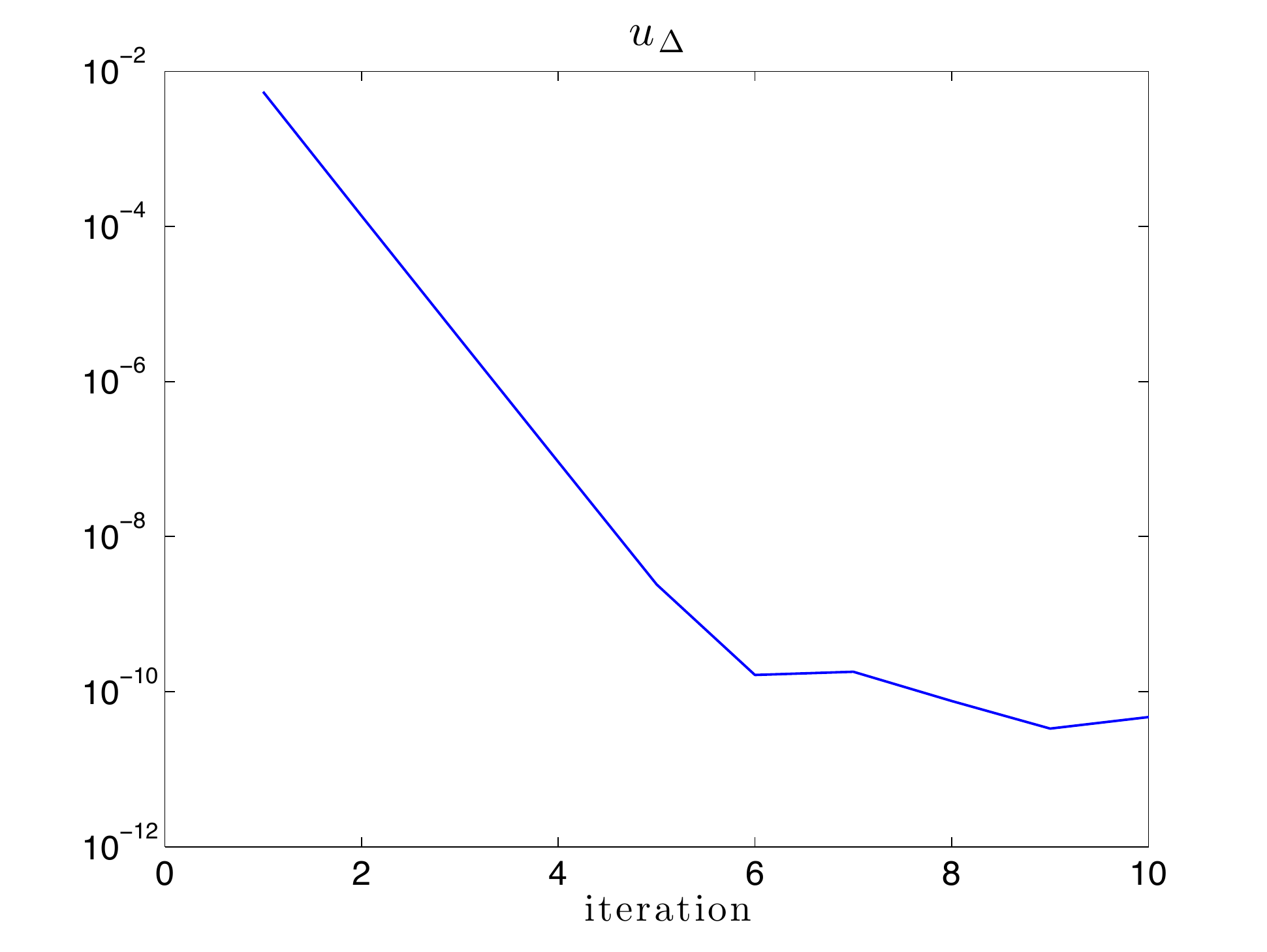}
\includegraphics[width=6.45cm]{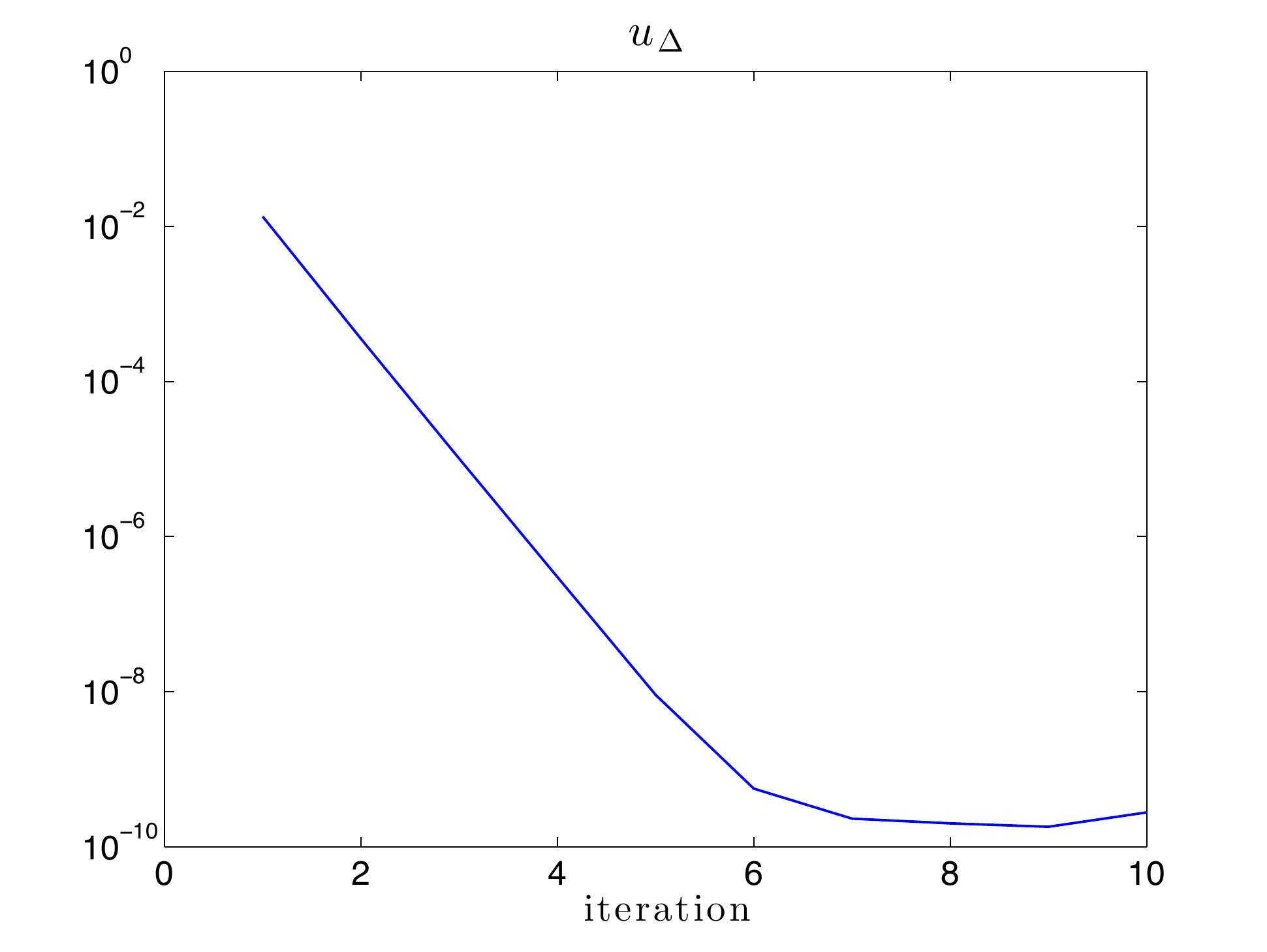}
\end{center}
\caption{Convergence history of convergence indicators~(\ref{eq:eps})
and~(\ref{eq:u_Delta}) corresponding to the smallest eigenvalue of the
Timoshenko beam with $CoV=10\%$ (left panels) and $CoV=25\%$ (right panels). } 
\label{fig:TB-e1_indicators} 
\end{figure}

Analogous computations for eigenvector errors and eigenproblem residuals are
summarized in Figures~\ref{fig:TB-it-CoV10}--\ref{fig:TB-it-CoV25}. These
figures show estimates of pdfs for the eigenvector error~(\ref{eq:eps_u}) and
the residual distribution~(\ref{eq:true_res}) for the eigenvalue/eigenvector
pair, corresponding to the smallest eigenvalue of the Timoshenko beam. The
trends for convergence in both figures (corresponding to $CoV=10\%$ and
$25\%$) are similar.

These figures provide insight into the maximal values of the errors~(\ref{eq:eps_u}) obtained from samples of the discrete eigenvectors.
For example, in the display in the upper left of Figure~\ref{fig:TB-it-CoV10}, the support
of the pdf for RQ$^{(0)}$ (obtained from the mean solution) is essentially the interval 
$[0,0.02]$, which shows that the eigenvector error~$\varepsilon_{u}$ from RQ$^{(0)}$
is of order at most $2\%$.
The analogous result for $CoV=25\%$ is $6\%$
(upper left of Figure~\ref{fig:TB-it-CoV25}), so that RQ$^{(0)}$ is less accurate for 
the larger value of $CoV$. 
Nevertheless, it can be
seen from Figure~\ref{fig:TB-it-CoV25} that even with $CoV=25\%$, the
eigenvector approximation error~$\varepsilon_{u}$ is less than $0.15\%$ after
one step of inverse iteration and after the second step $\varepsilon_{u}$ is
less than~$0.01\%$ and the error essentially coincides with the eigenvector
error from stochastic collocation. In other words, the convergence of~SII is
also indicated by the \textquotedblleft leftward\textquotedblright\ movement
of the pdfs corresponding to~$\varepsilon_{u}$. The pdf estimates of the
residuals are very small after one inverse iteration. We also found that when
the residual indicators (\ref{eq:eps}) stop decreasing and the differences
(\ref{eq:u_Delta}) become small, the sample true residuals~(\ref{eq:true_res})
also become small. Figure \ref{fig:TB-e1_indicators} shows the behavior of the
indicators (\ref{eq:eps})--(\ref{eq:u_Delta}).

Next, we consider the computation of multiple extreme eigenvalues. For the
stochastic Galerkin method, this entails construction of the coefficients of
$n_{s}>1$ eigenvalue fields in~(\ref{eq:RQ}). The stochastic collocation
method computes~$n_{s}$ extreme eigenvalues for each sample point and then
uses these to construct the random fields associated with each of them. Monte
Carlo proceeds in an analogous way.

\begin{figure}[ptbh]
\begin{center}
\includegraphics[width=6.45cm]{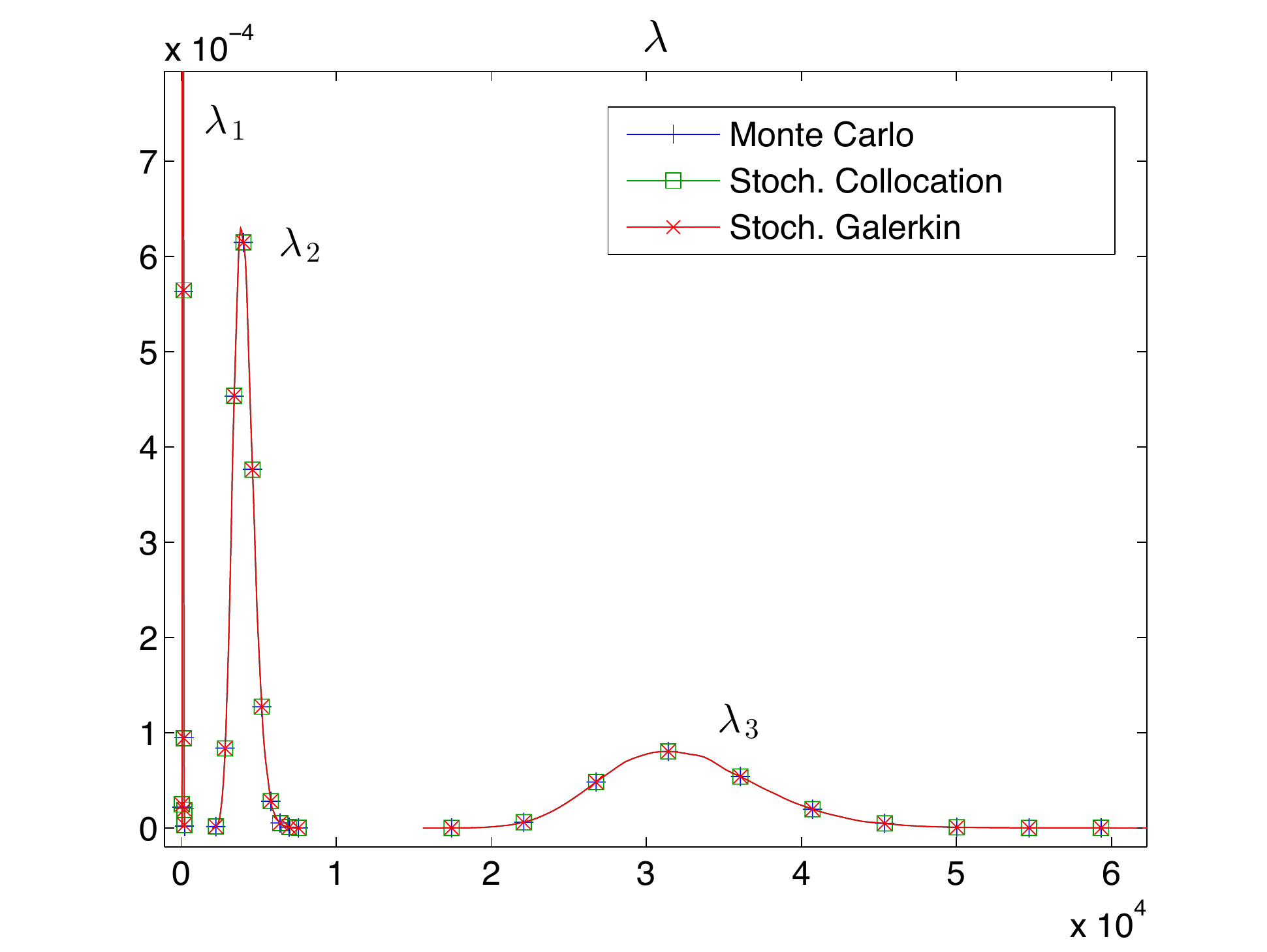}
\includegraphics[width=6.45cm]{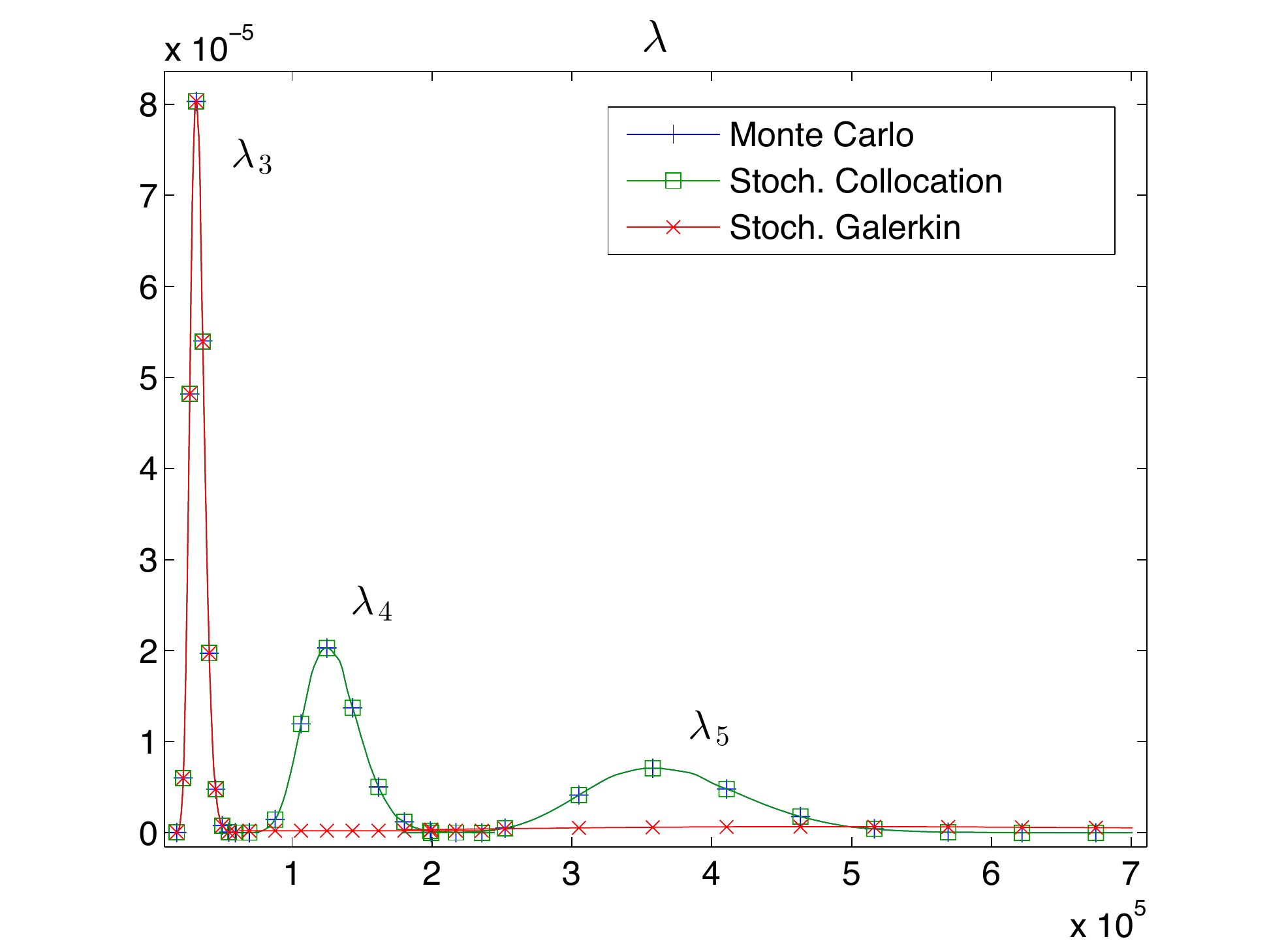}
\includegraphics[width=6.45cm]{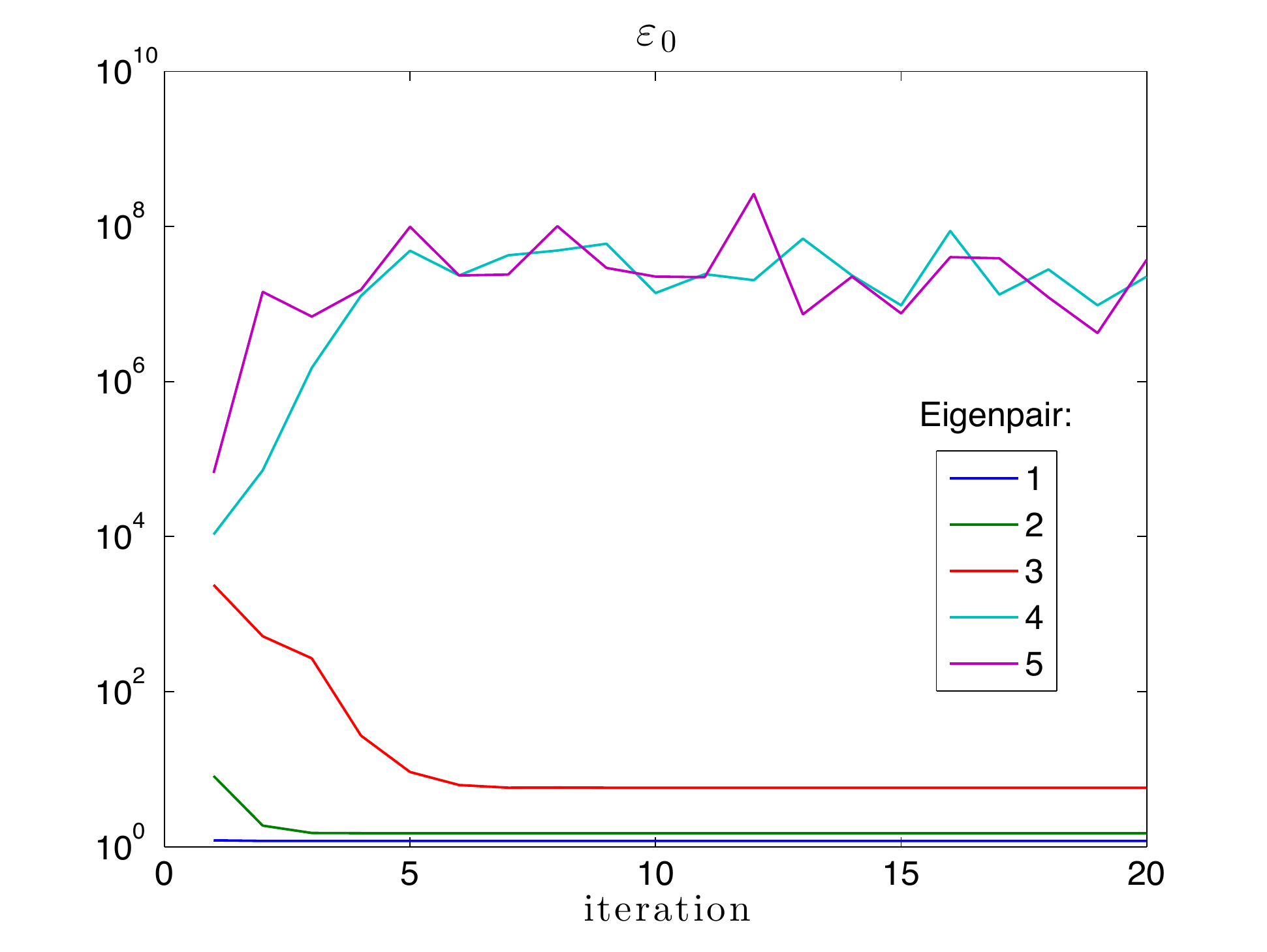}
\includegraphics[width=6.45cm]{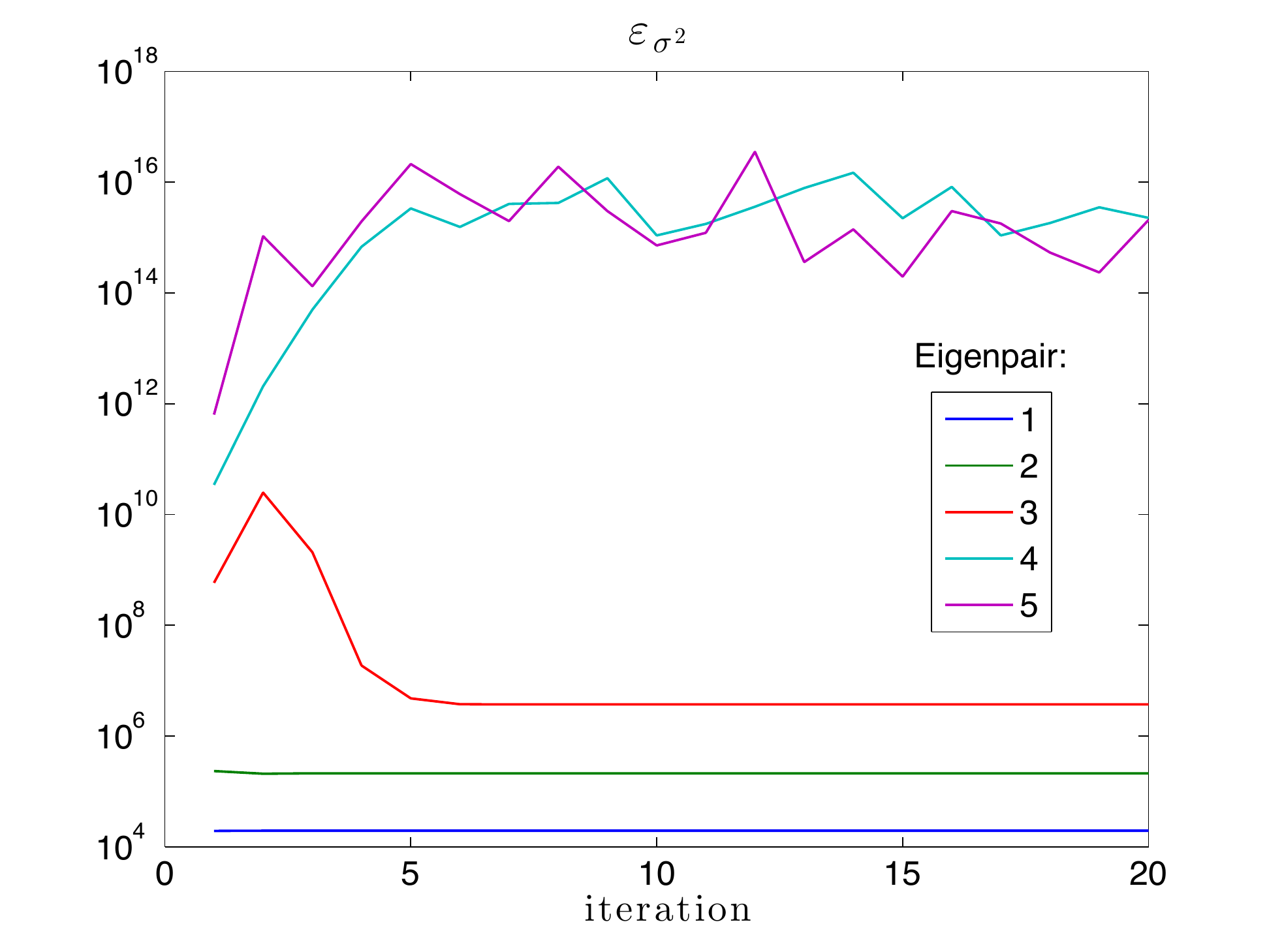}
\end{center}
\caption{Top: pdf estimates of the eigenvalue distribution corresponding to
eigenvalues $\lambda_{1},\lambda_{2},\lambda_{3}$ (left) $\lambda_{3}%
,\lambda_{4},\lambda_{5}$ (right). Bottom: convergence history of the two
indicators $\varepsilon_{0}$ and $\varepsilon_{\sigma}^{2}$ from
(\ref{eq:eps}) obtained using inverse subspace iteration, for the five
smallest eigenvalues of the Timoshenko beam with $CoV=25\%$.} 
\label{fig:TB-min5} 
\end{figure}

\begin{figure}[ptbh]
\begin{center}
\includegraphics[width=6.45cm]{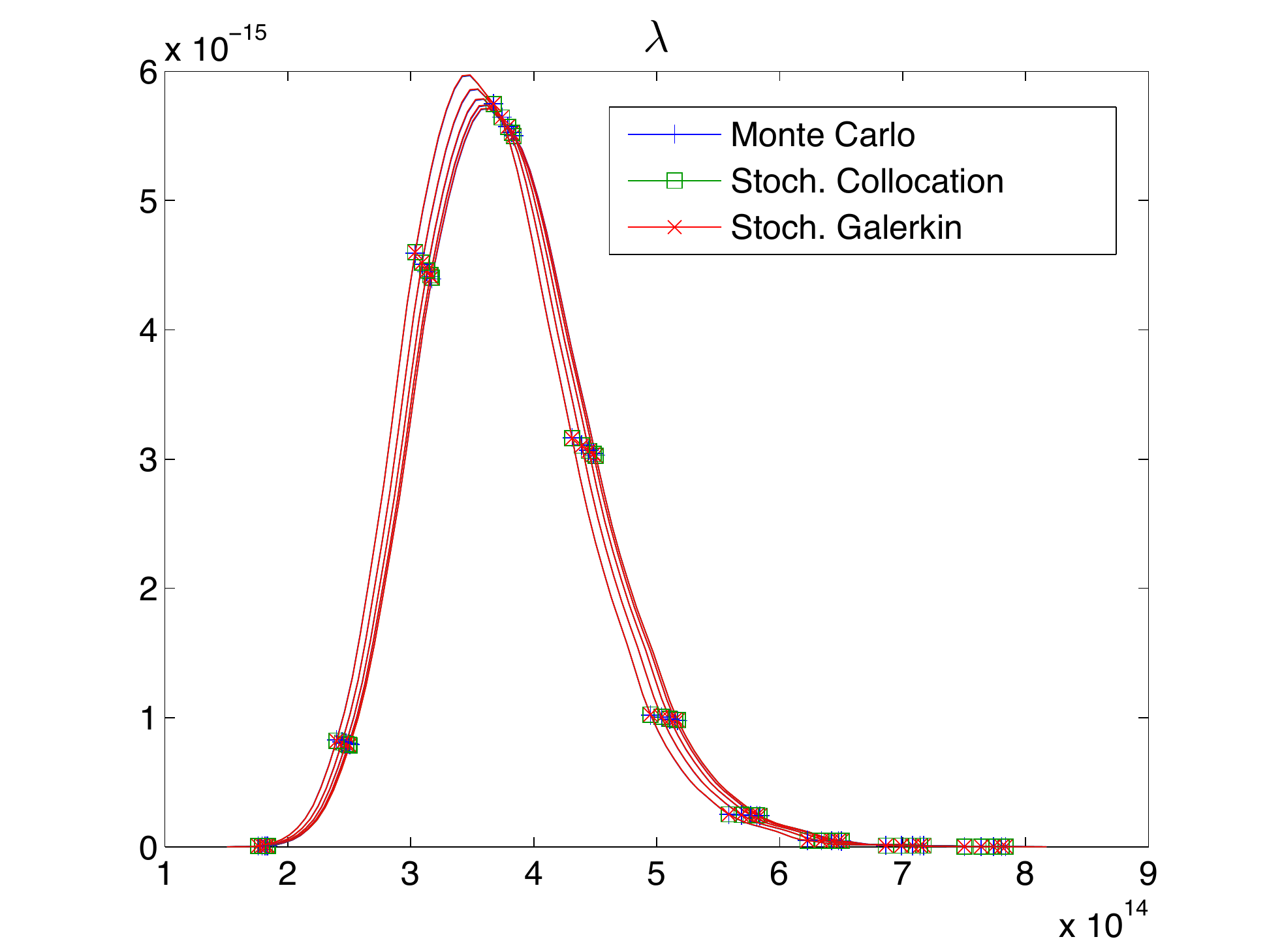}
\includegraphics[width=6.45cm]{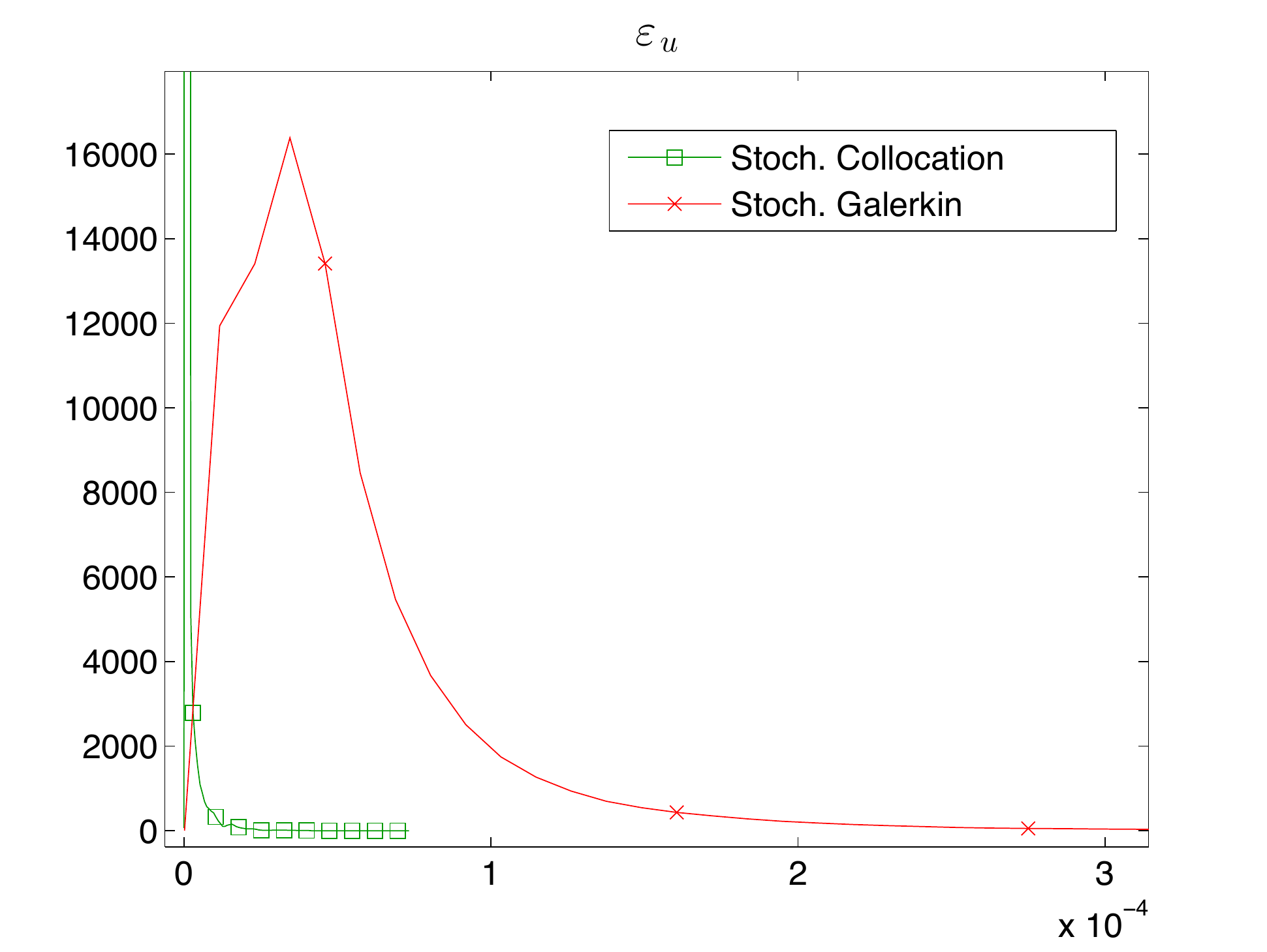}
\includegraphics[width=6.45cm]{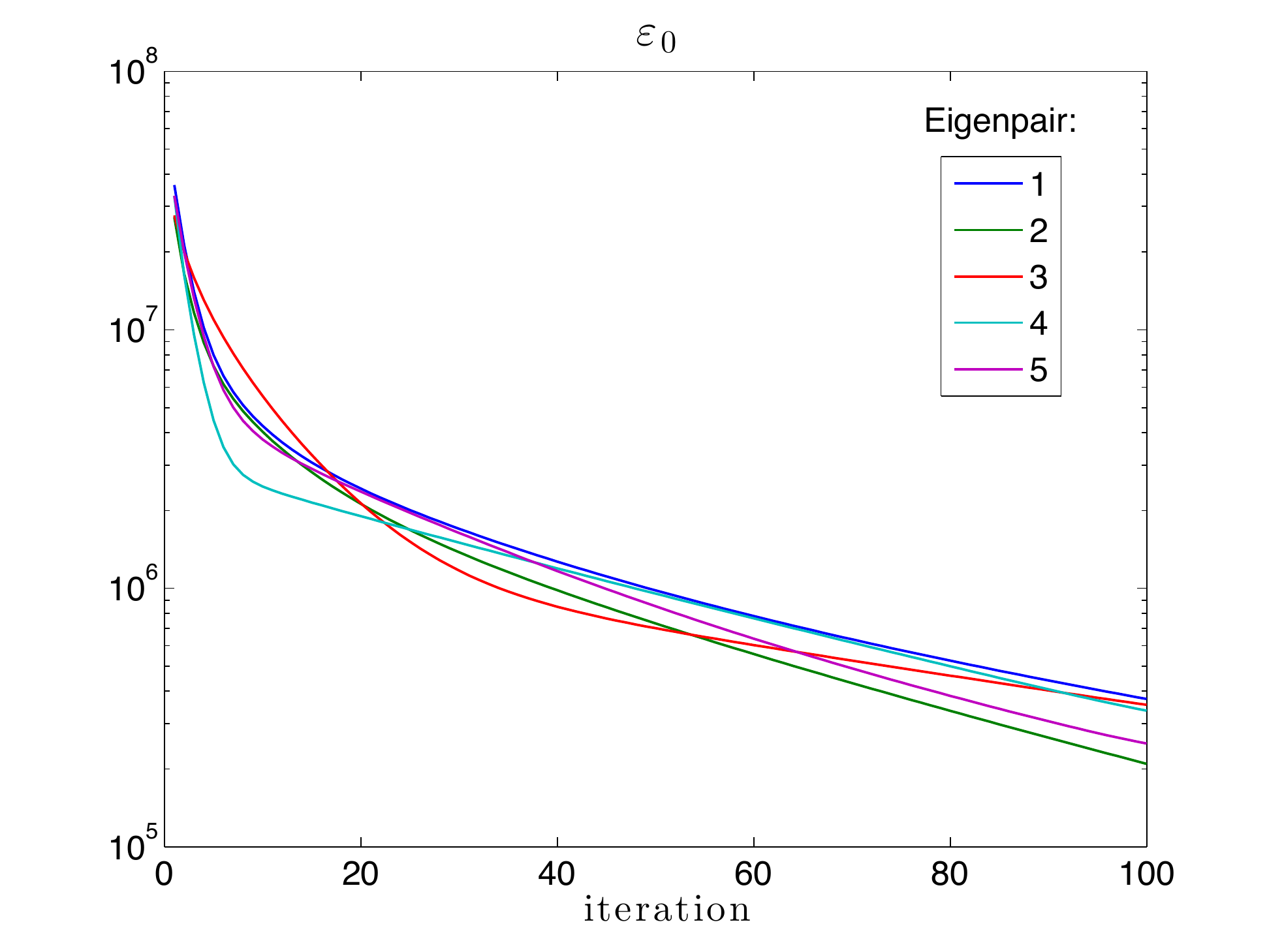}
\includegraphics[width=6.45cm]{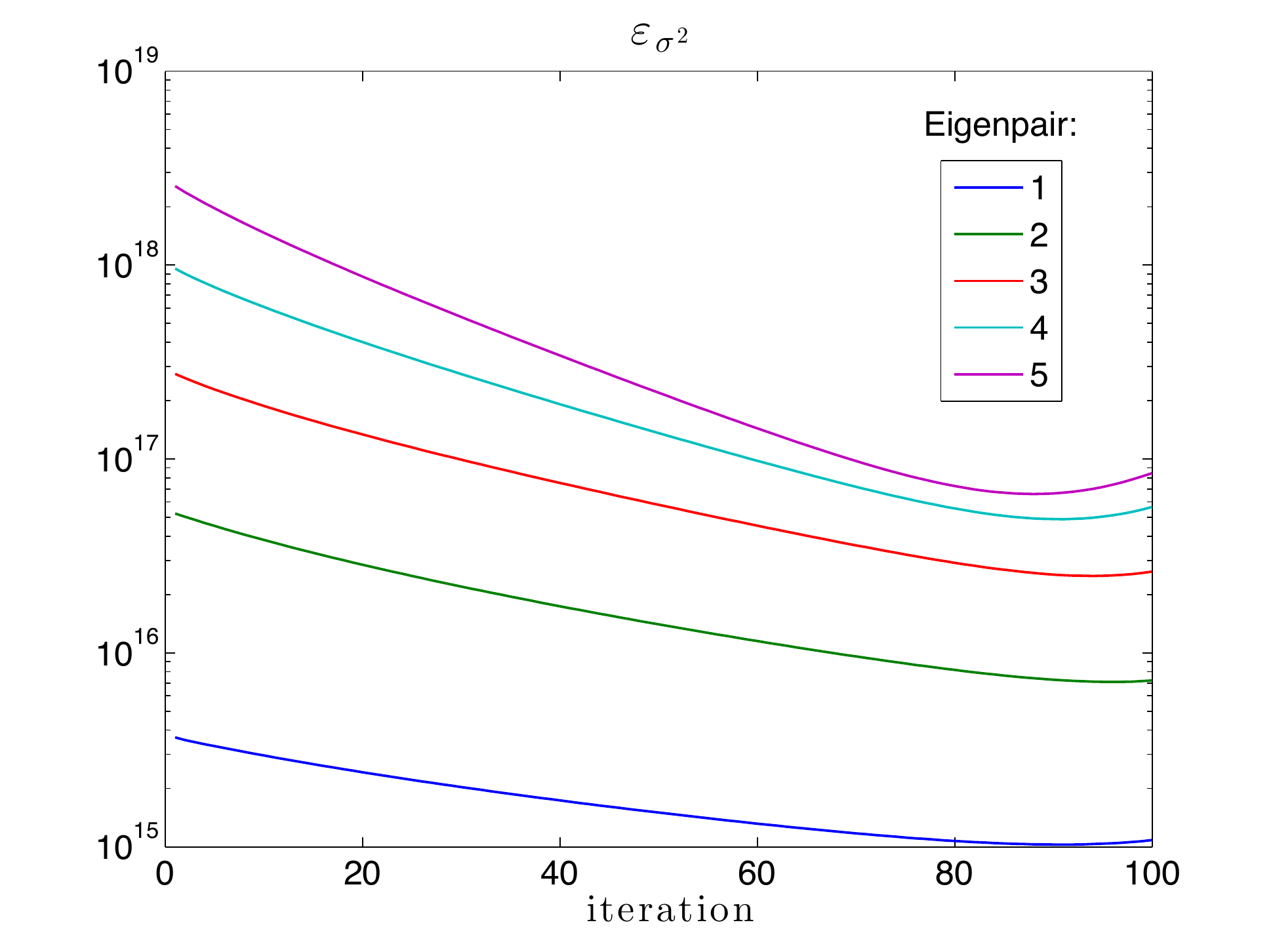}
\end{center}
\caption{Top: pdf estimates of the five largest eigenvalues of the Timoshenko
beam with $CoV=25\%$ (left), and $\ell^{2}$-norm of the relative eigenvector
error (\ref{eq:eps_u}) (right). Bottom: convergence history of the two
indicators $\varepsilon_{0}$ and $\varepsilon_{\sigma^{2}}$ 
from~(\ref{eq:eps}).} 
\label{fig:TB-max5} 
\end{figure}

The performance of the methods for computing the five smallest eigenvalues of
the Timoshenko beam with $CoV=25\%$ is shown in Figure \ref{fig:TB-min5}.
Stochastic Galerkin was able to identify the three smallest eigenvalues
$\lambda_{1},\lambda_{2},\lambda_{3}$, but it failed to identify eigenvalues
$\lambda_{4},\lambda_{5}$. 
(Results were similar for larger values of polynomial degree, $p=4$ and $5$.) 
Stochastic collocation and Monte Carlo were able to
find all five eigenvalues. Note that the error indicators $\varepsilon_{0}$
and $\varepsilon_{\sigma}^{2}$ from~(\ref{eq:eps}), shown in the bottom of the
figure, become flat for the converged eigenvalues but not for those that are
not found. 
Performance results for the five \emph{largest} eigenvalues are shown in
Figure~\ref{fig:TB-max5}. The Galerkin method was robust in this case: for
each of the five eigenvalues, the pdf estimates obtained by all three
computational methods overlap, and the $\ell^{2}$-norm of the relative
eigenvector error~(\ref{eq:eps_u}) corresponding to the fifth maximal
eigenvalue is small. The error indicator~$\varepsilon_{\sigma^{2}}$
from~(\ref{eq:eps}) behaves somewhat inconsistently in this case: after
initial decrease it can be seen that $\varepsilon_{\sigma^{2}}$ increases
slightly after approximately $85$~iterations.

\begin{figure}[ptbh]
\begin{center}
\includegraphics[width=6.45cm]{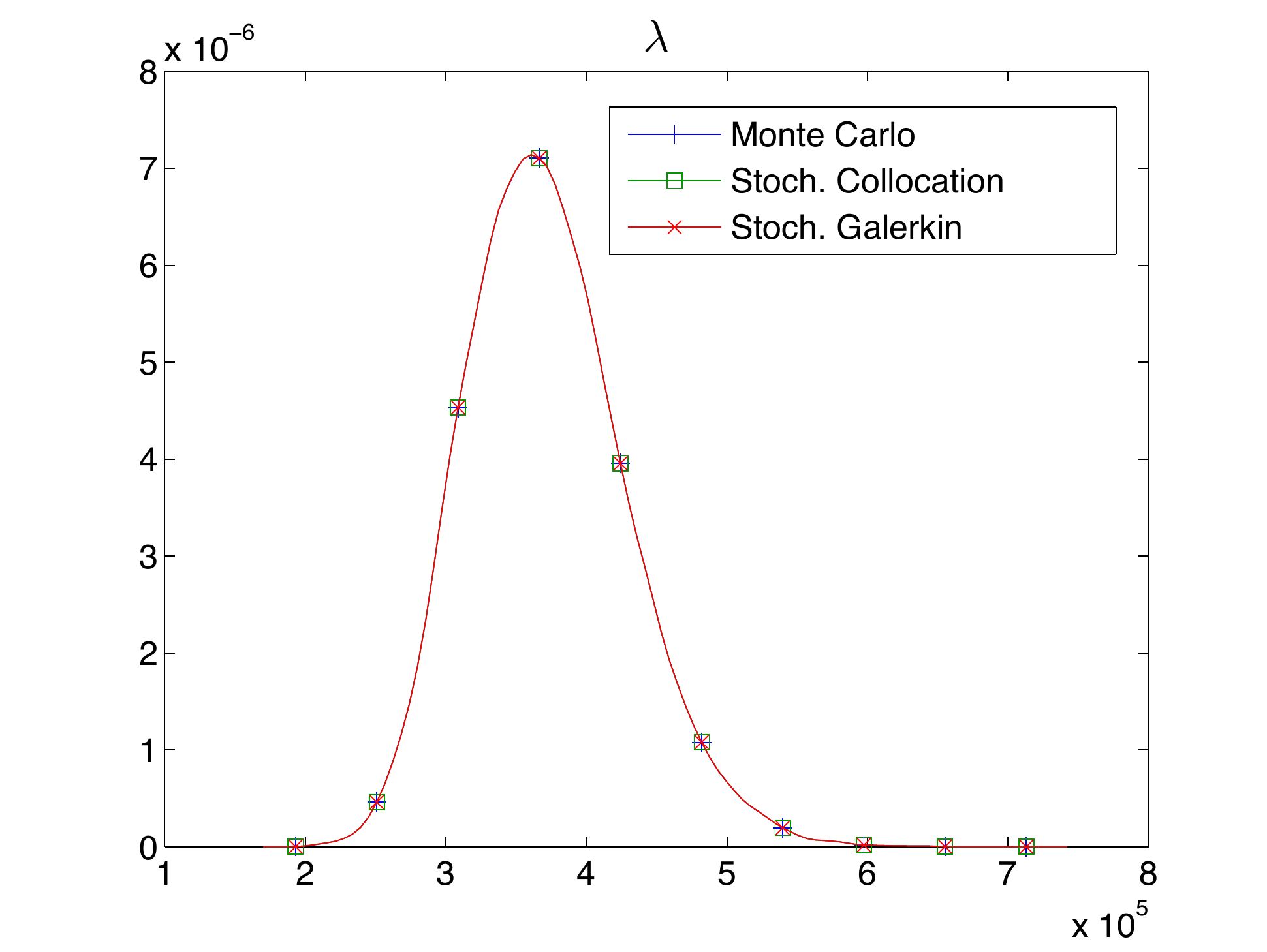}
\includegraphics[width=6.45cm]{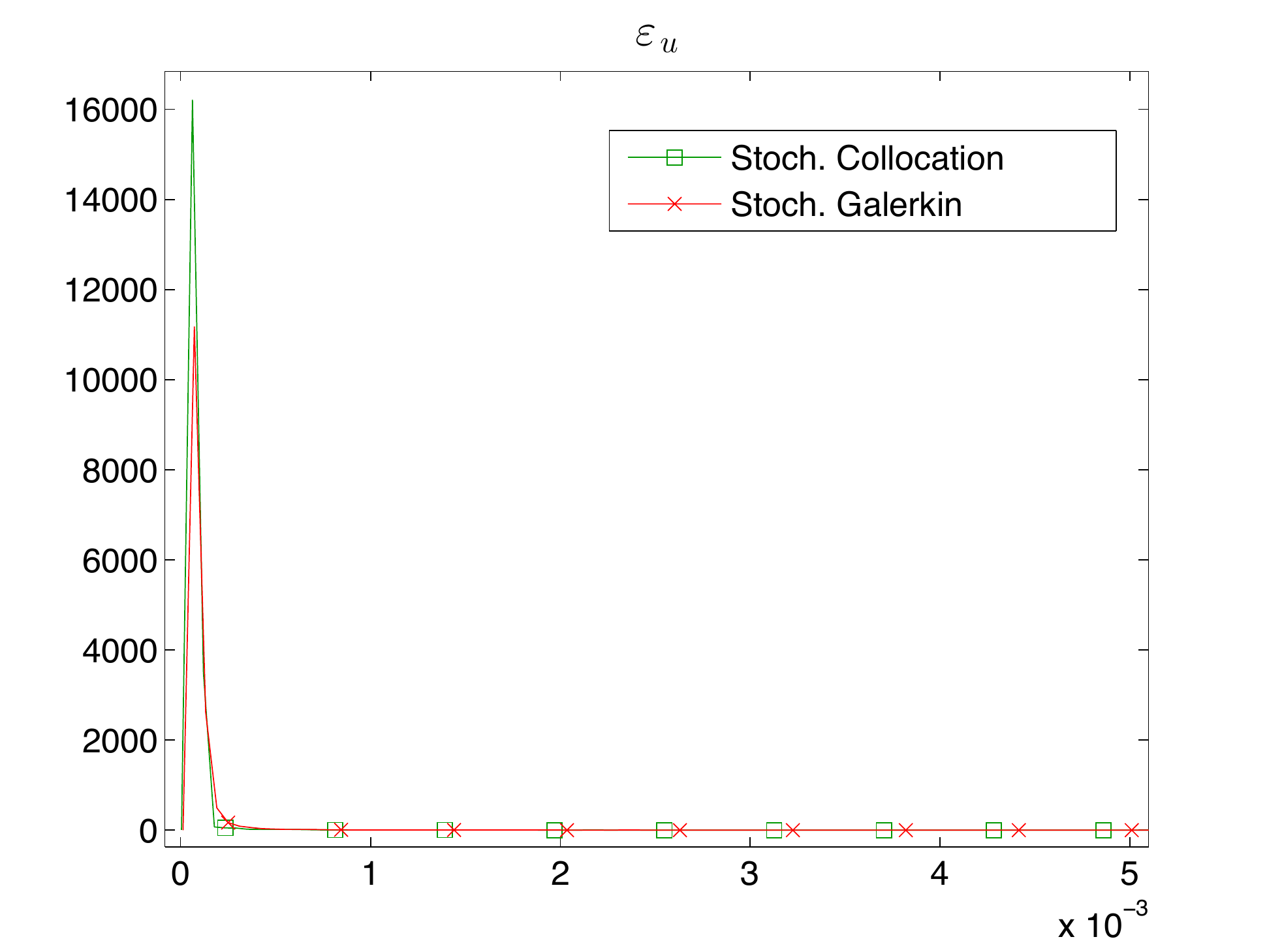}
\includegraphics[width=6.45cm]{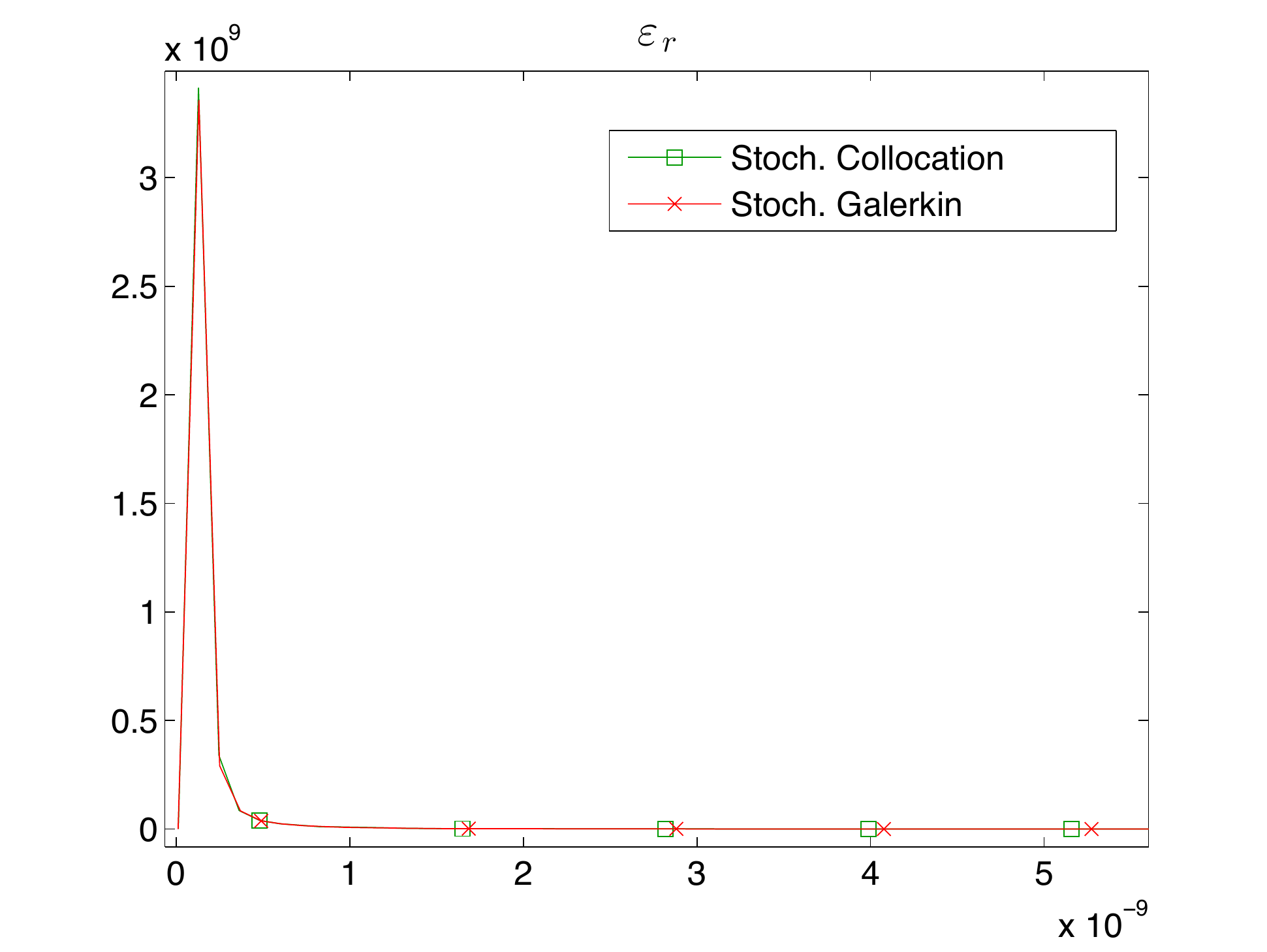}
\includegraphics[width=6.45cm]{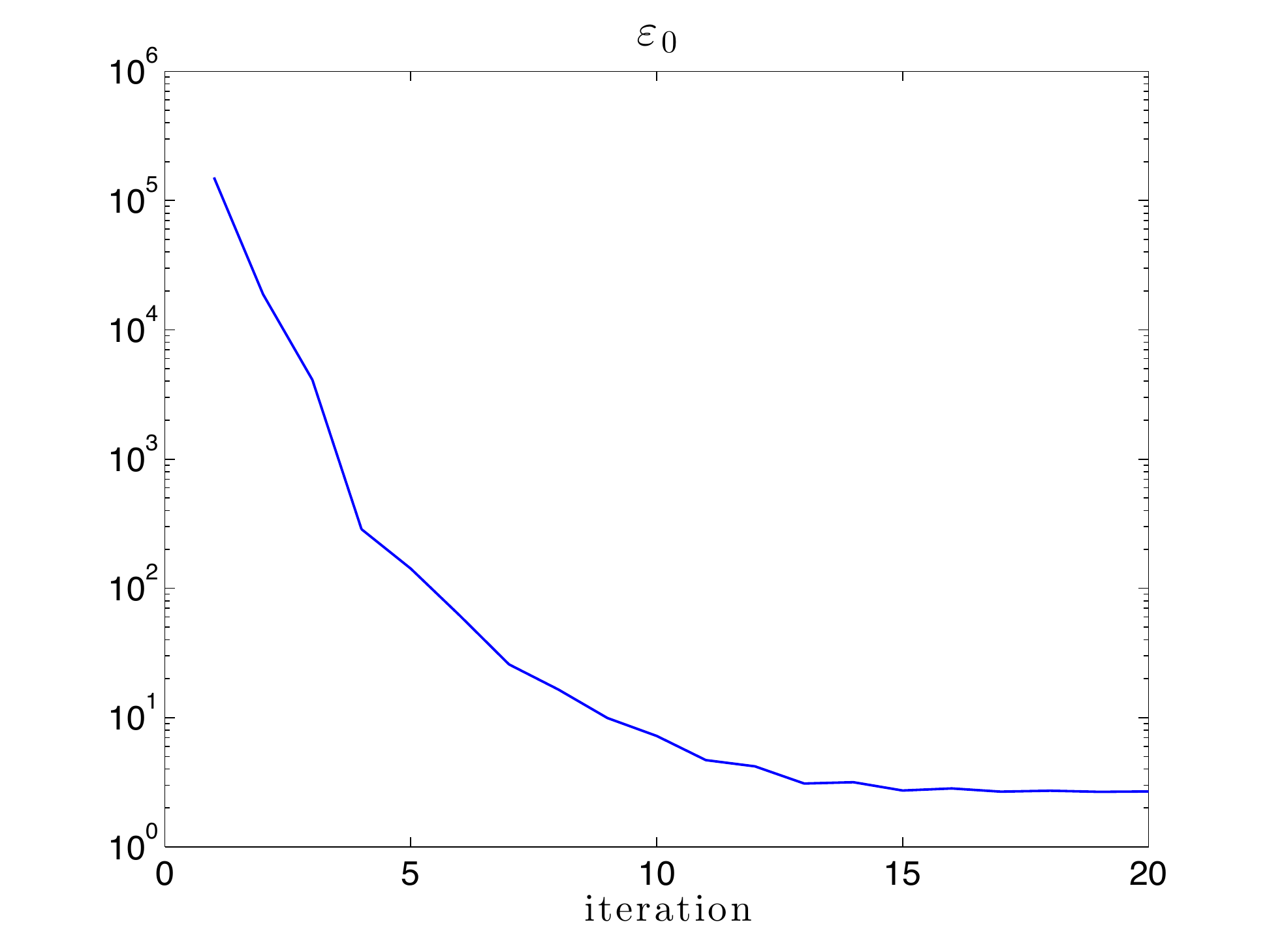}
\includegraphics[width=6.45cm]{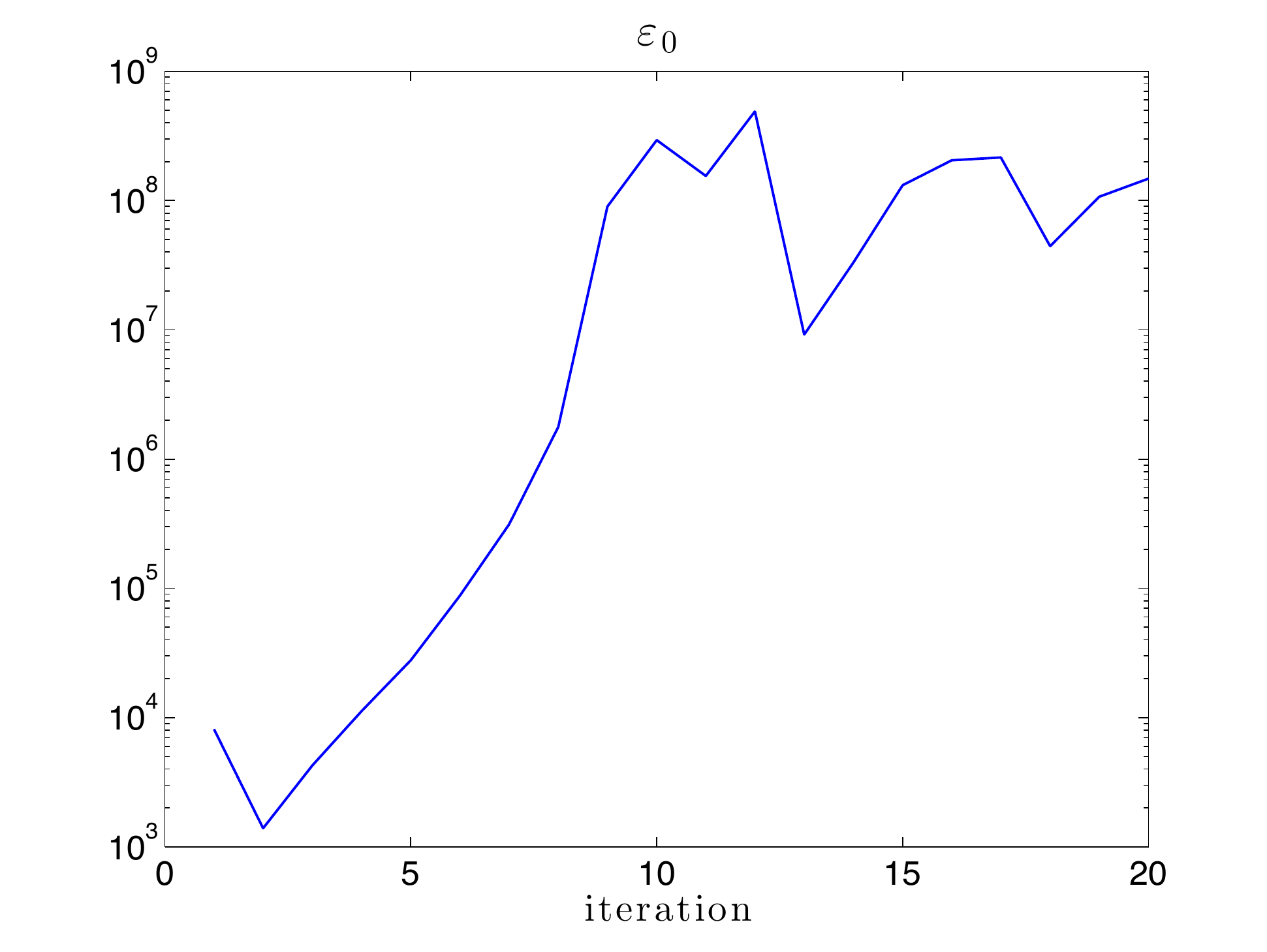}
\includegraphics[width=6.45cm]{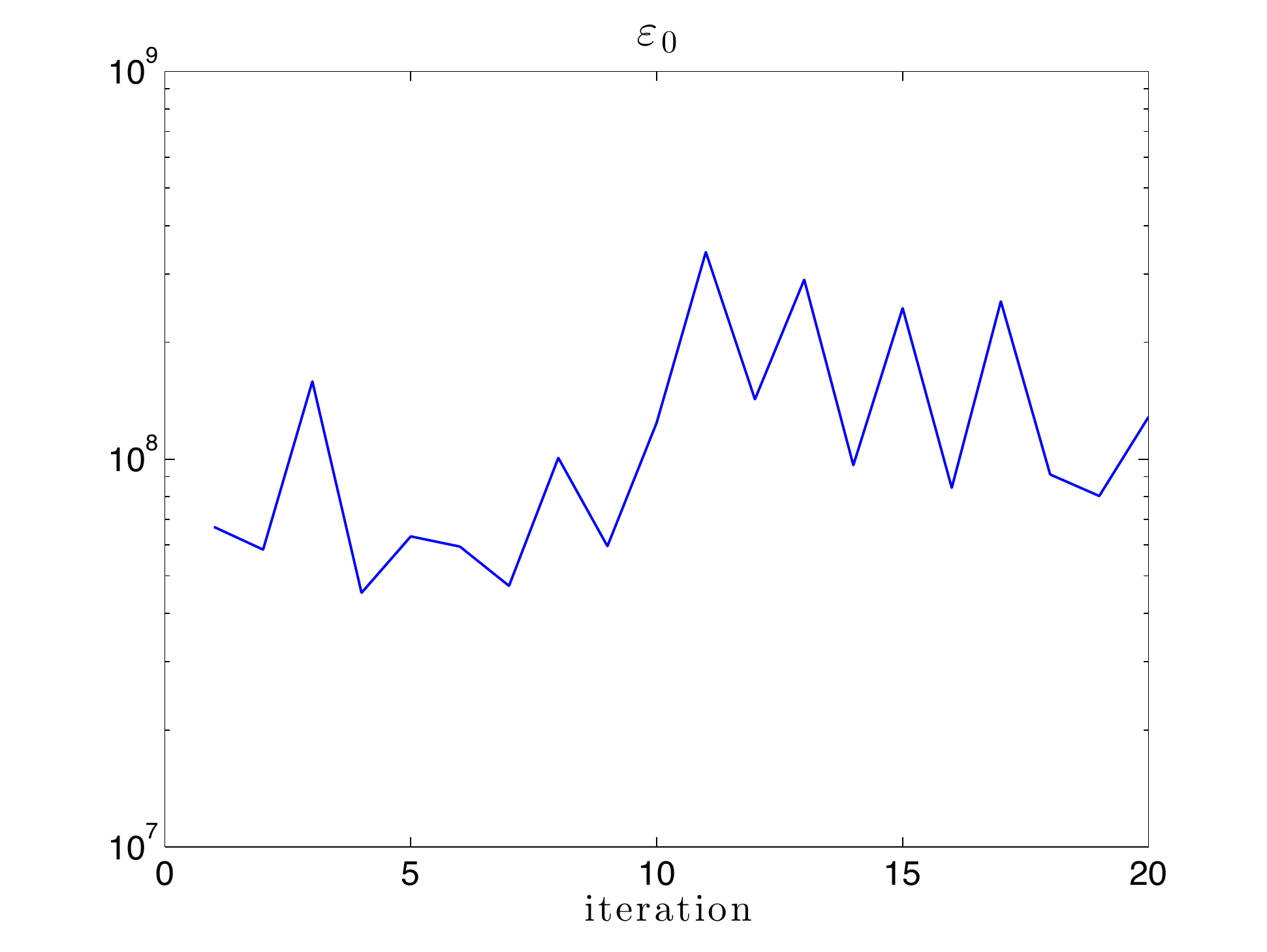}
\end{center}
\caption{Behavior of stochastic inverse iteration with shifts, for the fifth
smallest eigenvalue of the Timoshenko beam, with $CoV=25\%$. Top: pdf
estimates of the eigenvalue distribution (left) and $\ell^{2}$-norm of the
relative eigenvector error~(\ref{eq:eps_u}) (right), obtained using shift
$\rho=4.1\times10^{5}$. Middle: pdf estimate of the true residual
(\ref{eq:true_res}) (left) and convergence history of the indicator
$\varepsilon_{0}$ from (\ref{eq:eps}) (right), also with $\rho=4.1\times
10^{5}$. Bottom: stochastic inverse iteration fails to converge with shift
$\rho=3.9\times10^{5}$ (left) or $\rho=4.3\times10^{5}$ (right), as
illustrated by the convergence history of $\varepsilon_{0}$.} 
\label{fig:TB-5th-q}
\end{figure}

We explored several approaches to enhance the robustness of stochastic
subspace iteration for identifying interior eigenvalues. One possibility is to
use a shift. We tested inverse iteration with a shift to find the fifth
smallest eigenvalue of the Timoshenko beam with $CoV=25\%$. The corresponding
eigenvalue of the mean problem is $\overline{\lambda}_{5}=3.7548\times10^{5}$.
The top four panels in Figure~\ref{fig:TB-5th-q} show plots of the pdf
estimates of the eigenvalue distribution, the $\ell^{2}$-norm of the relative
eigenvector error~(\ref{eq:eps_u}), the true residual~(\ref{eq:true_res}), and
the convergence history of the indicator$~\varepsilon_{0}$ from~(\ref{eq:eps})
with the shift$~\rho=4.1\times10^{5}$. It can be seen that for the estimates
of the pdfs of the eigenvalue, the relative eigenvector errors, and the true
residual of the stochastic inverse iteration, the methods are in agreement.
However, we also found that convergence depends on the choice of the shift
$\rho$. Setting the shift far from the eigenvalue of interest or too close to
it worsens the convergence rate and the method might even fail to converge.
For this eigenvalue, the best convergence occurs with the shift set close to
either $\rho=3.5\times10^{5}$ or $\rho=4.1\times10^{5}$, but with shift set to
$\rho=3.9\times10^{5}$\ or $\rho=4.3\times10^{5}$ the method fails to
converge. Similar behavior was also reported in~\cite{Verhoosel-2006-ISR}. We
note that the mean of the sixth smallest eigenvalue is $\overline{\lambda}%
_{6}=8.9196\times10^{5}$, that is, the means of the fifth and sixth
eigenvalues are well separated, see also Figure~\ref{fig:TB-eig-A_0}.

\begin{figure}[ptbh]
\begin{center}
\includegraphics[width=6.45cm]{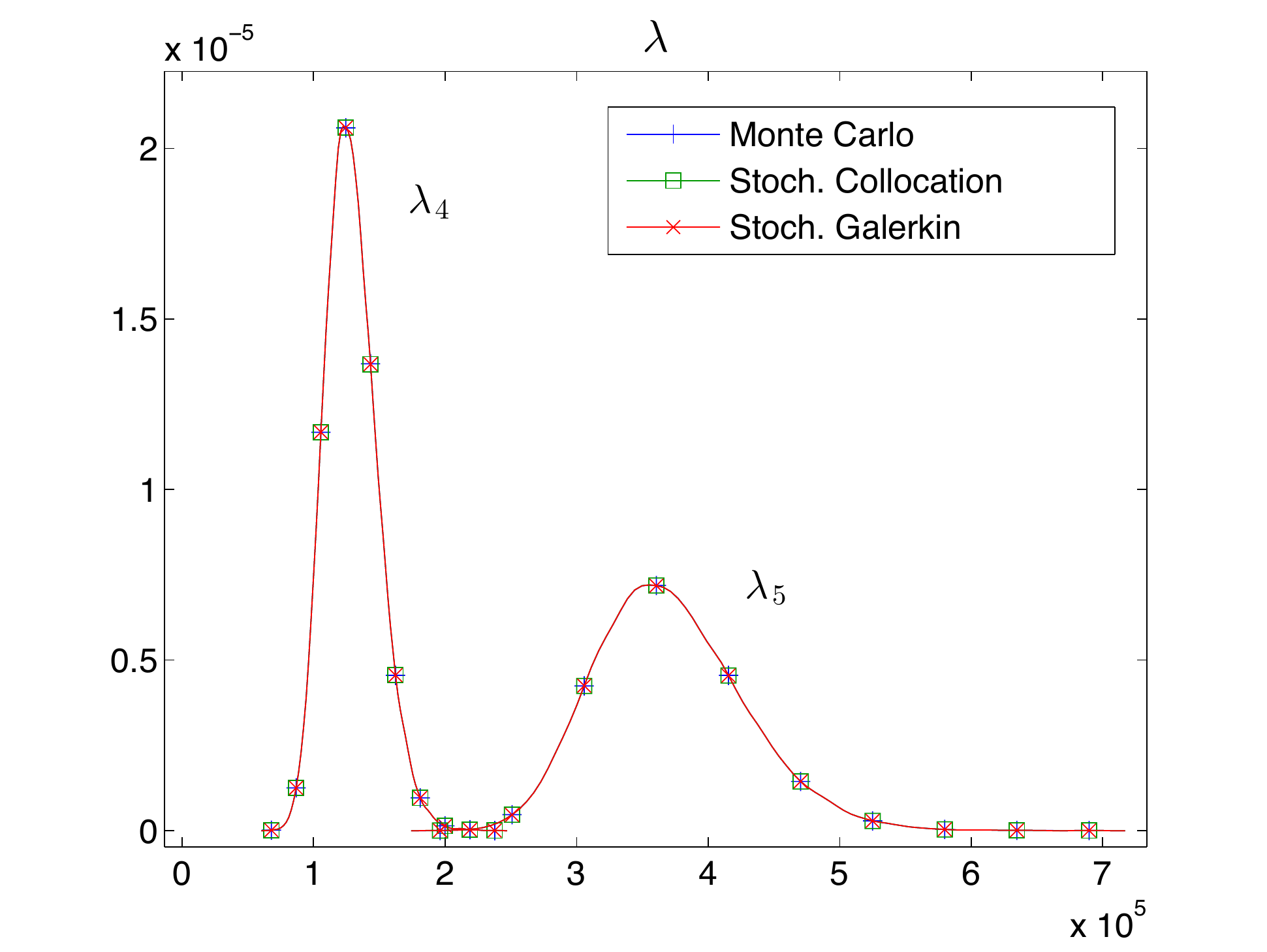}
\includegraphics[width=6.45cm]{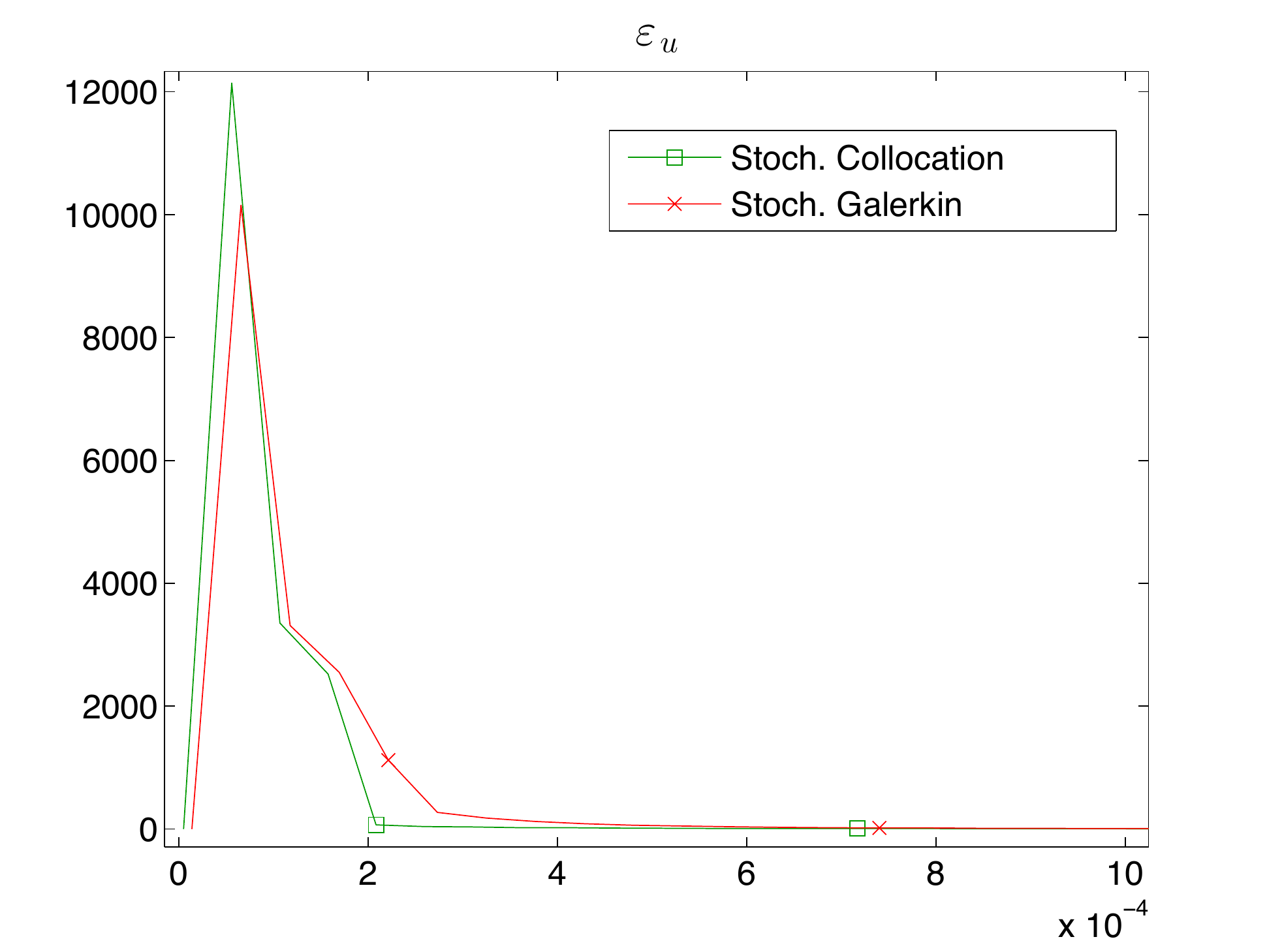}
\includegraphics[width=6.45cm]{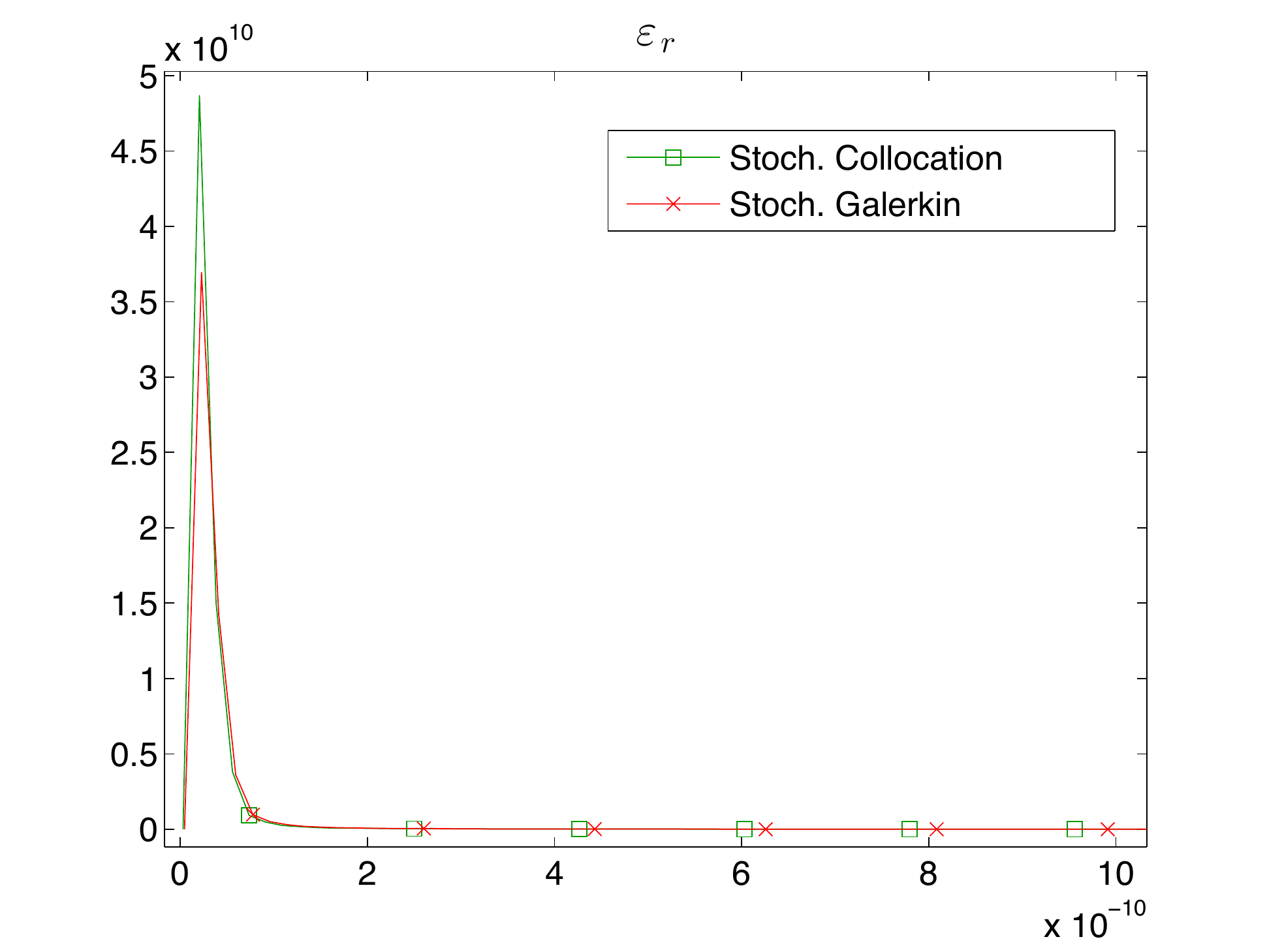}
\includegraphics[width=6.45cm]{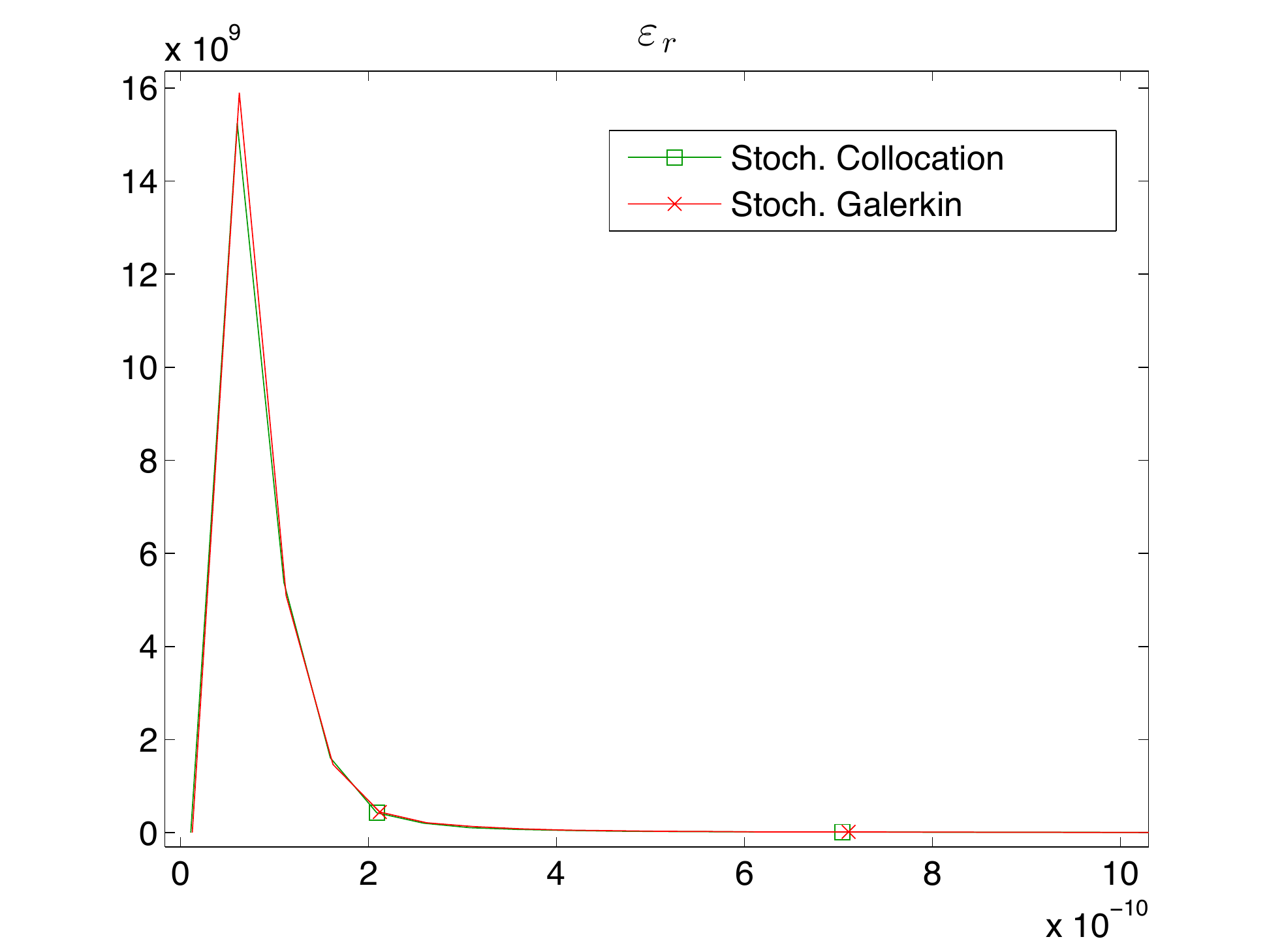}
\includegraphics[width=6.45cm]{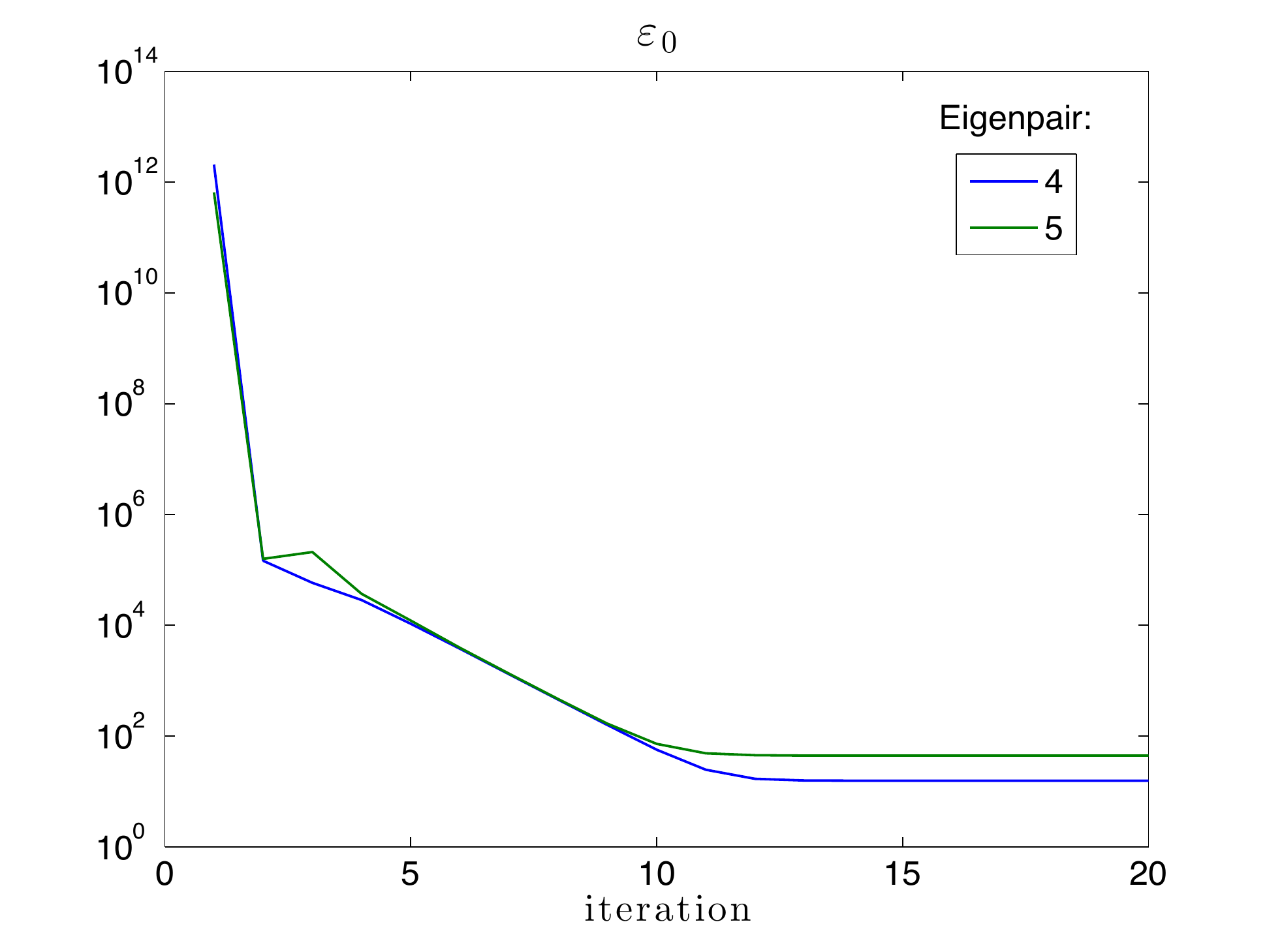}
\includegraphics[width=6.45cm]{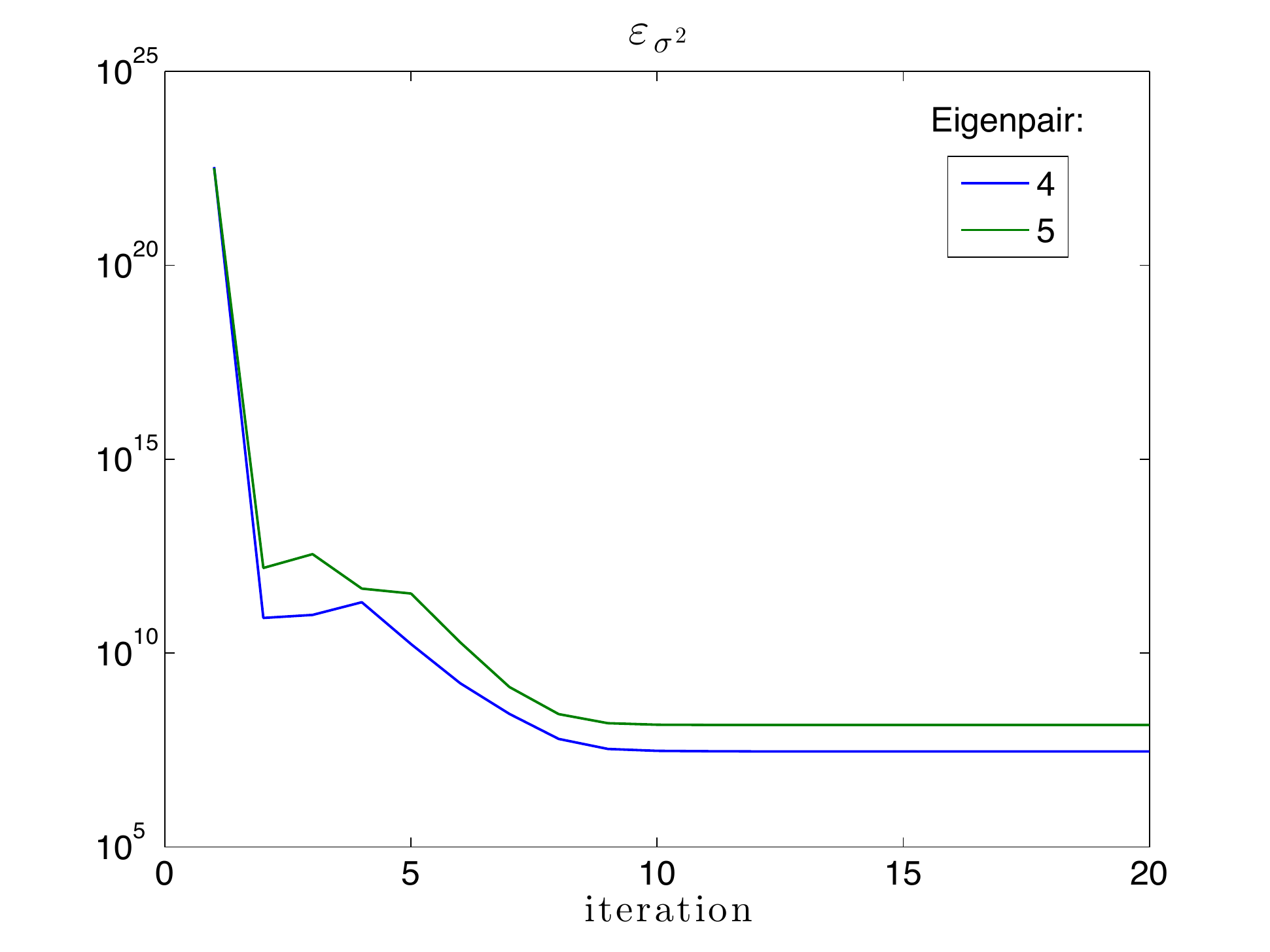}
\end{center}
\caption{Top: pdf estimates of the eigenvalue distribution $\lambda_{4}$ and
$\lambda_{5}$ (left), and of the $\ell^{2}$-norm of the relative eigenvector
error~(\ref{eq:eps_u}) corresponding to $\lambda_{5}$ (right). Middle: pdf of
the $\ell^{2}$-norm of the true residual~(\ref{eq:true_res}) corresponding to
eigenvalues~$\lambda_{4}$ (left), and~$\lambda_{5}$ (right). Bottom:
convergence history of the two indicators~$\varepsilon_{0}$ and~$\varepsilon
_{\sigma}^{2}$ from~(\ref{eq:eps}) corresponding to eigenvalues~$\lambda_{4}$
and~$\lambda_{5}$\ of the Timoshenko beam with $CoV=25\%$ obtained using
inverse subspace iteration and deflation~(\ref{eq:deflation-2}) of the three
smallest eigenvalues $\lambda_{1}$, $\lambda_{2}$, and $\lambda_{3}$.}%
\label{fig:TB-defl}%
\end{figure}

\begin{figure}[ptbh]
\begin{center}
\includegraphics[width=6.45cm]{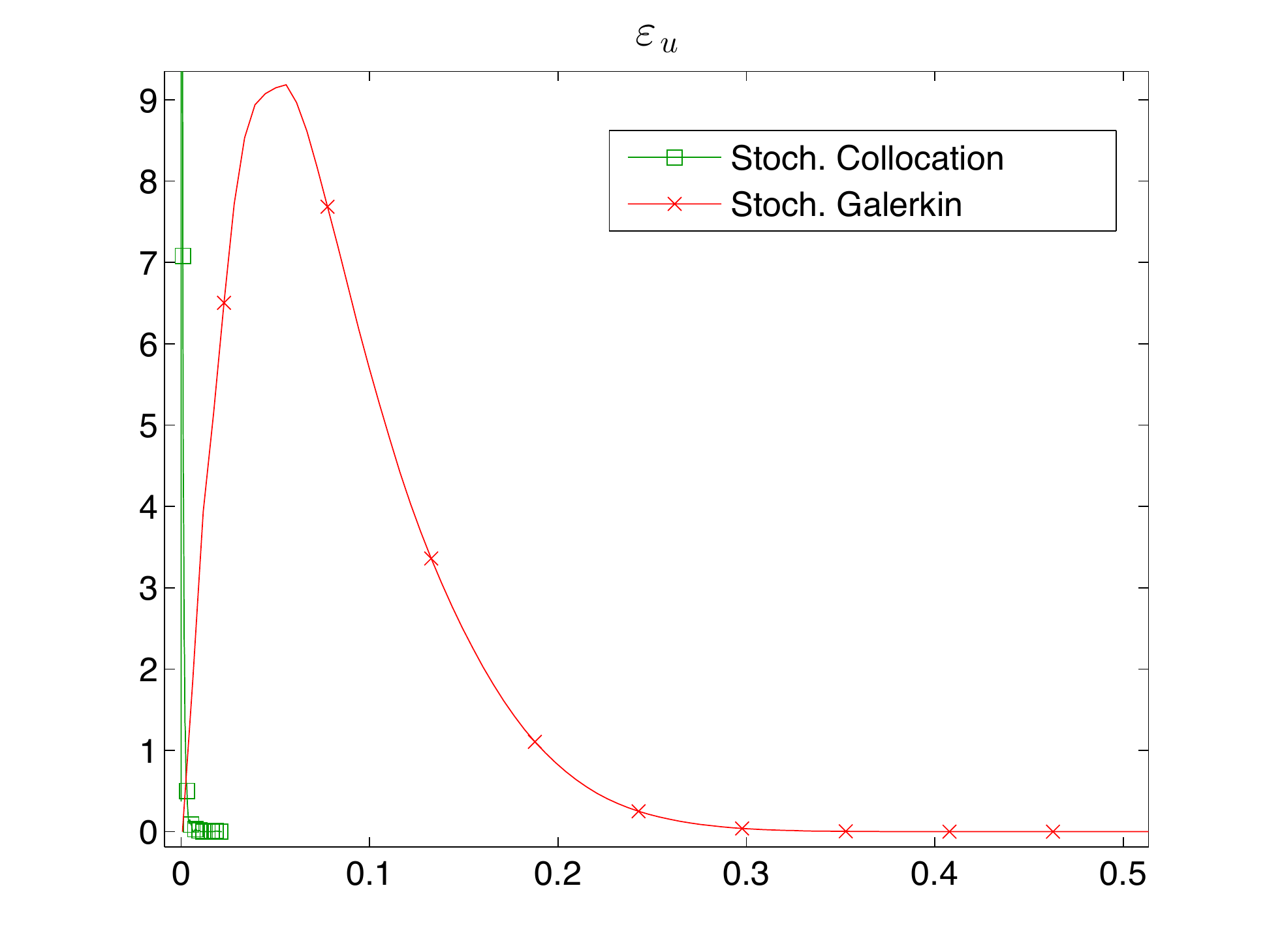}
\includegraphics[width=6.45cm]{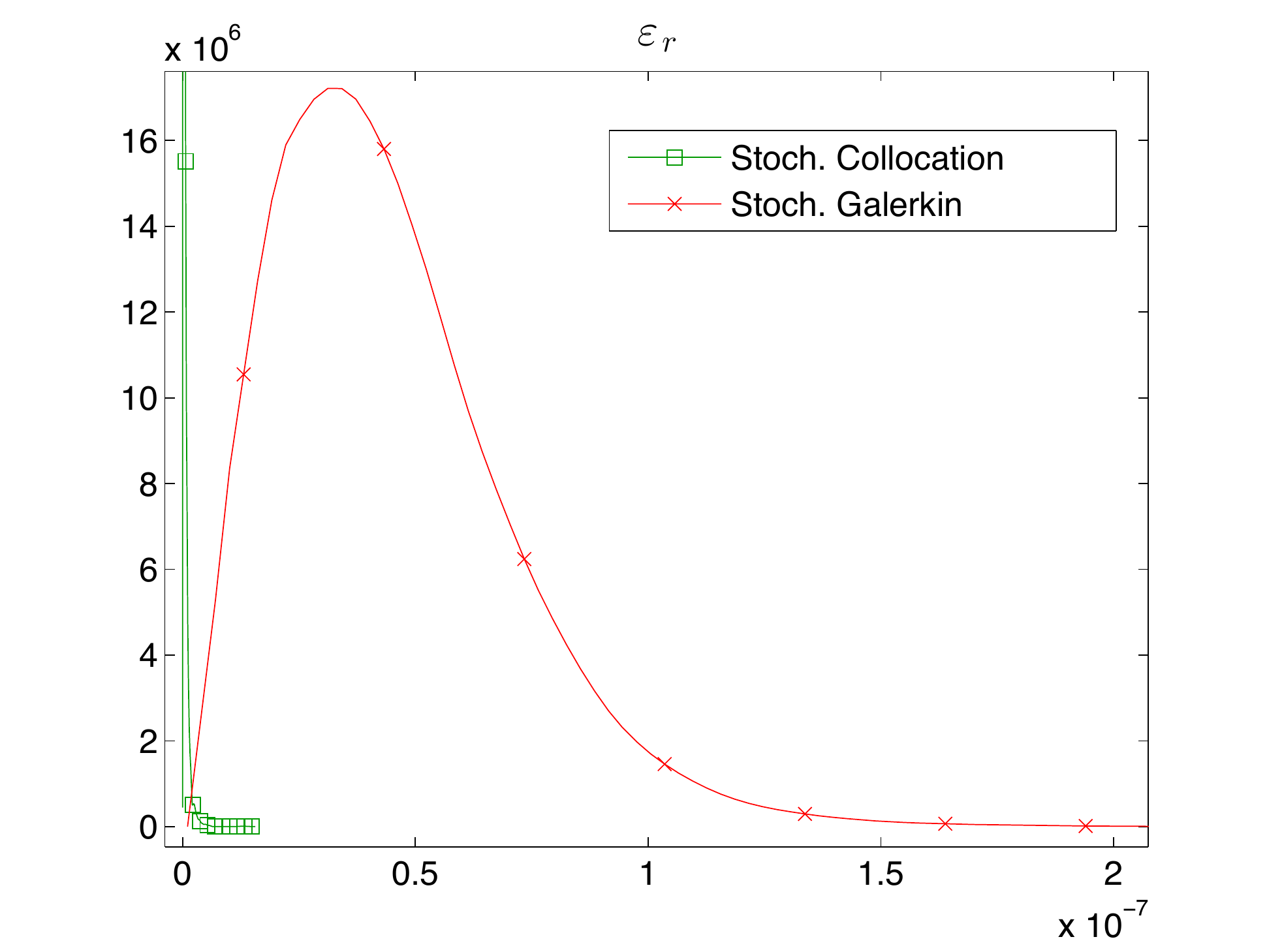}
\includegraphics[width=6.45cm]{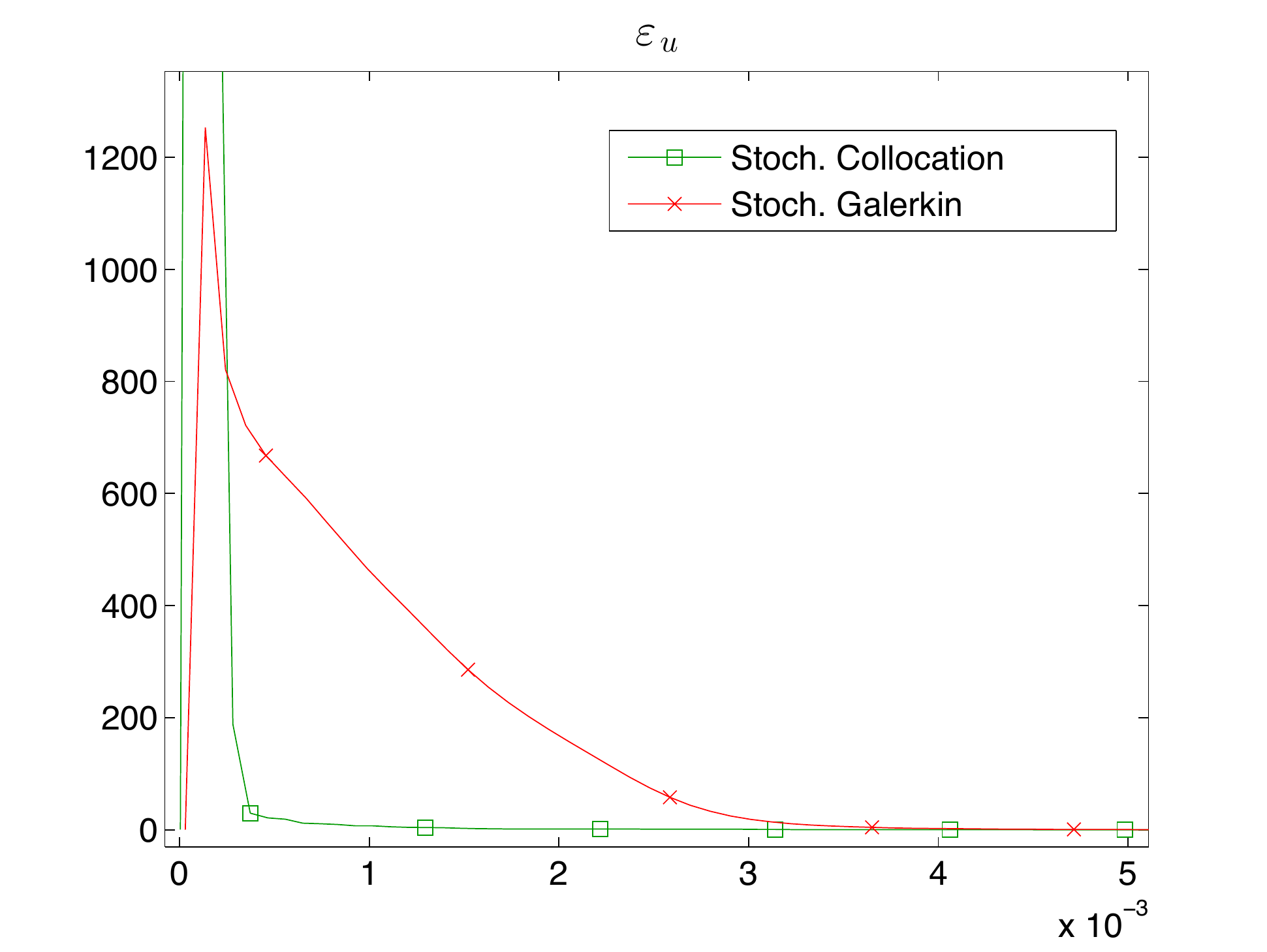}
\includegraphics[width=6.45cm]{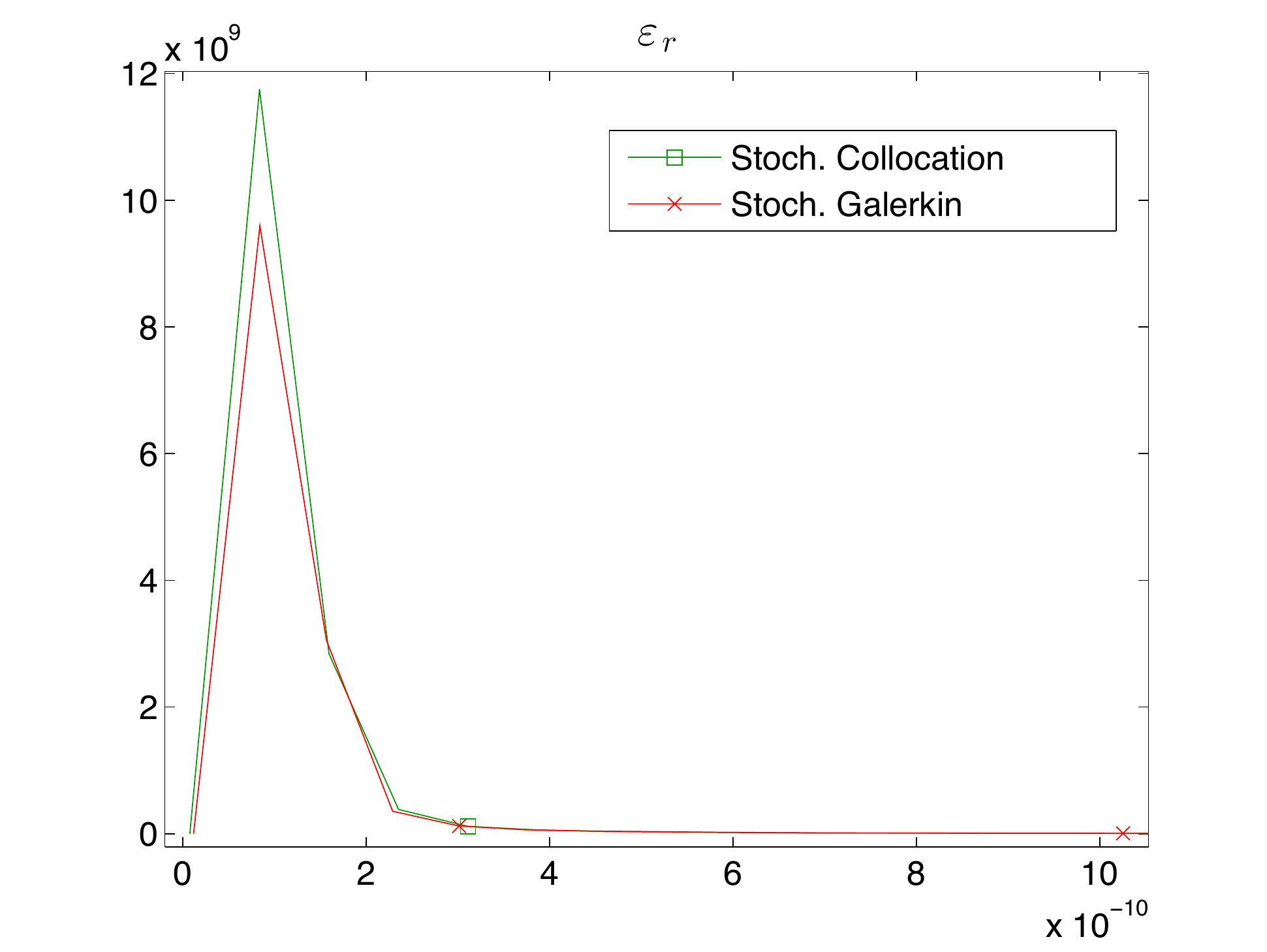}
\includegraphics[width=6.45cm]{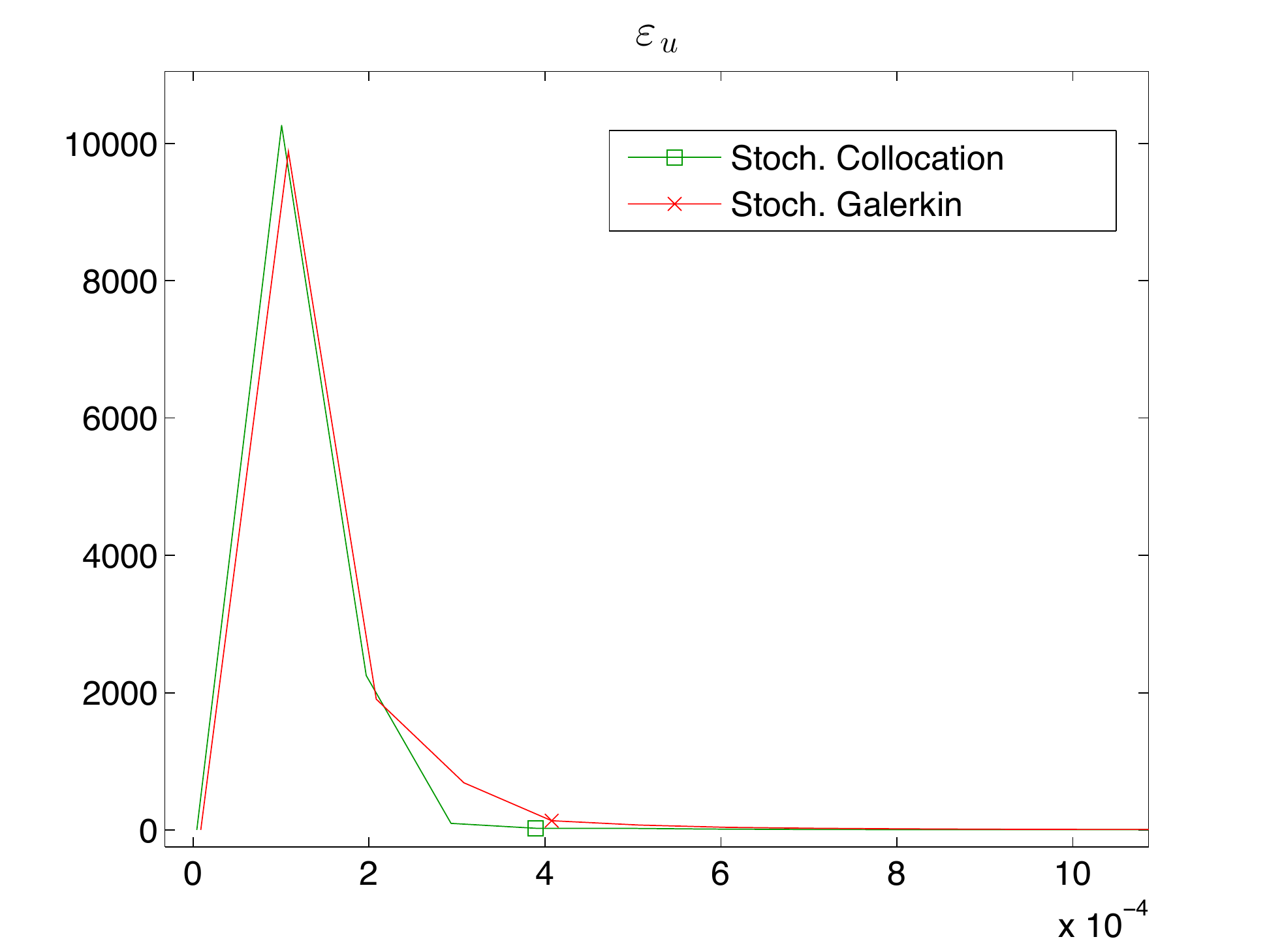}
\includegraphics[width=6.45cm]{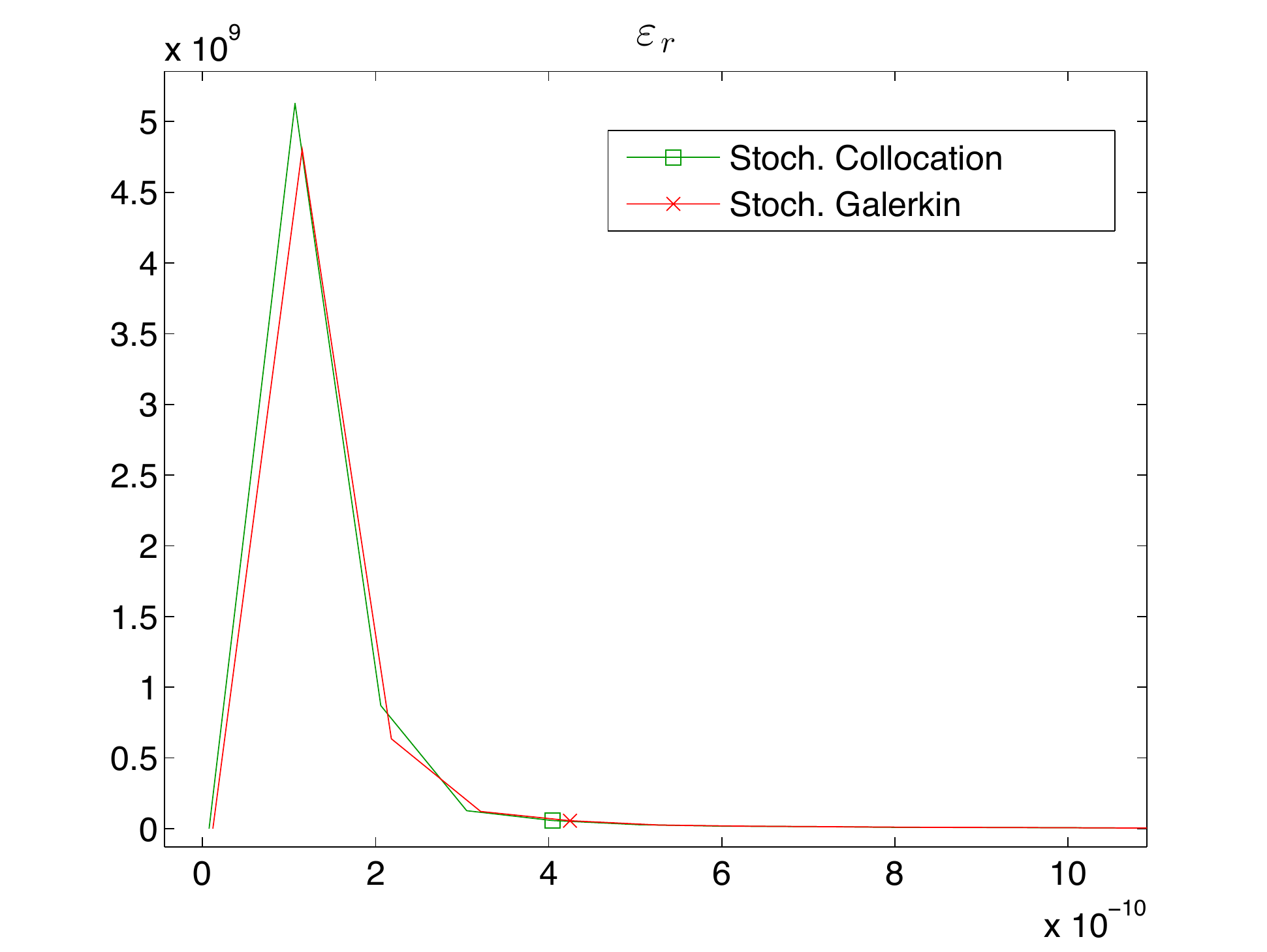}
\end{center}
\caption{Plots of the pdf estimate of the $\ell^{2}$-norms of the relative
eigenvector error~(\ref{eq:eps_u}) and the true residual~(\ref{eq:true_res})
corresponding to the fifth smallest eigenvalue of the Timoshenko beam with
$CoV=25\%$ obtained using inverse iteration with deflation of the four
smallest eigenvalues in iteration $0$ (top), $5$ (middle) and $10$ (bottom).}%
\label{fig:TB-it} 
\end{figure}

An approach that we found to be more robust was to use
deflation of the mean matrix. 
Suppose we are interested in some interior eigenvalues in the lower side of
the spectrum, for example $\lambda_{4}$ and $\lambda_{5}$, which we were
unable to identify in a previous attempt (Figure \ref{fig:TB-min5}). To
address this, as suggested in (\ref{eq:deflation-2}) we can deflate the mean
matrix $A_{0}$ using the mean eigenvectors corresponding to $\lambda_{1}$,
$\lambda_{2}$\ and $\lambda_{3}$. Figure~\ref{fig:TB-defl} shows that in this
case, Algorithm \ref{alg:SISI} was able to identify the fourth and fifth
smallest eigenvalues, and the relative eigenvector errors~(\ref{eq:eps_u})
almost coincide. We note that the results in Figures~\ref{fig:TB-defl} (and
also in Figure \ref{fig:TB-it}) were obtained using the deflated mean matrix
also in stochastic collocation and Monte Carlo methods.

One significant advantage of stochastic inverse subspace iteration over Monte
Carlo and stochastic collocation is that it allows termination of the
iteration at any step, and thus the coefficients of the expansions
(\ref{eq:gPC-lambda-u}) can be found only approximately.
Figure~\ref{fig:TB-it} shows the $\ell^{2}$-norms of the relative eigenvector
error~(\ref{eq:eps_u}) and the pdf estimates of the true
residual~(\ref{eq:true_res}) corresponding to the fifth smallest eigenvalue of
the Timoshenko beam with $CoV=25\%$, obtained using inverse iteration with
deflation of the four smallest eigenvalues in iteration $0$, $5$, and $10$.
For example, the initial mean of the relative eigenvector error$~\varepsilon
_{u}$ from ~(\ref{eq:eps_u}) is centered around $10\%$, after $5$ iterations
it is reduced to less than $0.5\%$, and after $10$ iterations the results of
stochastic inverse iteration and stochastic collocation essentially agree, and
the difference from Monte Carlo represented by$~\varepsilon_{u}$ is less than
$0.05\%$.

\subsection{Example 2: Mindlin plate}

\begin{figure}[ptbh]
\begin{center}
\includegraphics[width=8.5cm]{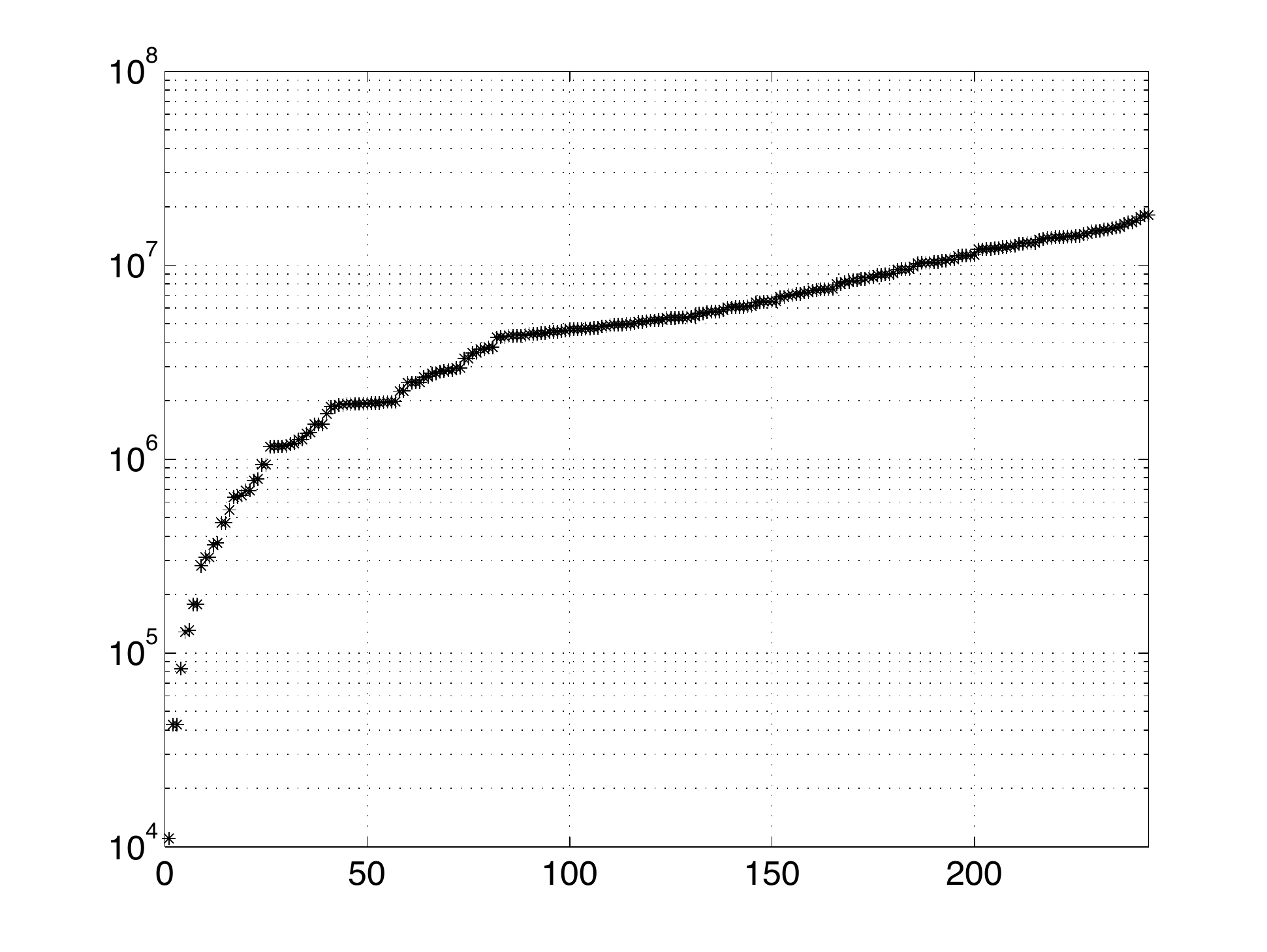}
\end{center}
\caption{Eigenvalues of the matrix $A_{0}$ corresponding to the Mindlin
plate.} 
\label{fig:MP-eig-A_0} 
\end{figure}

For the second example, we analyzed vibrations of a square, fully simply
supported Mindlin plate. For this problem we used $3\times10^{4}$ Monte Carlo
samples. The physical parameters were set according to~\cite[Section
12.5]{Ferreira-2009-MCF-1Ed} as follows: the mean Young's modulus of the
lognormal random field was $E_{0}=10,920$, Poisson's ratio $\nu=0.30$, length
of a side $L_{\text{plate}}=1$, thickness $0.1$, $\kappa=5/6$, and density
$\rho=1$. The plate was discretized using $10\times10$ bilinear (Q4) finite
elements with $243$ physical degrees of freedom. The condition number of the
mean matrix$~A_{0}$ from (\ref{eq:A_ell})\ is$~1.6436\times10^{3}$, the norm
$\left\Vert A_{0}\right\Vert _{2}=1.8153\times10^{7}$, and the eigenvalues
of$~A_{0}$\ are displayed in Figure~\ref{fig:MP-eig-A_0}.
Coefficient of variation of the Young's modulus was set to $CoV=25\%$, and the
spatial correlation length $L_{corr}=L_{\text{plate}}/4$. This is a
two-dimensional problem which means that there are repeated eigenvalues: for
example, the four smallest eigenvalues of the mean problem are $\overline
{\lambda}_{1}=1.1044\times10^{4}$, $\overline{\lambda}_{2}=\overline{\lambda
}_{3}=4.2720\times10^{4}$, and $\overline{\lambda}_{4}=8.3014\times10^{4}$.

\begin{figure}[ptbh]
\begin{center}
\includegraphics[width=6.45cm]{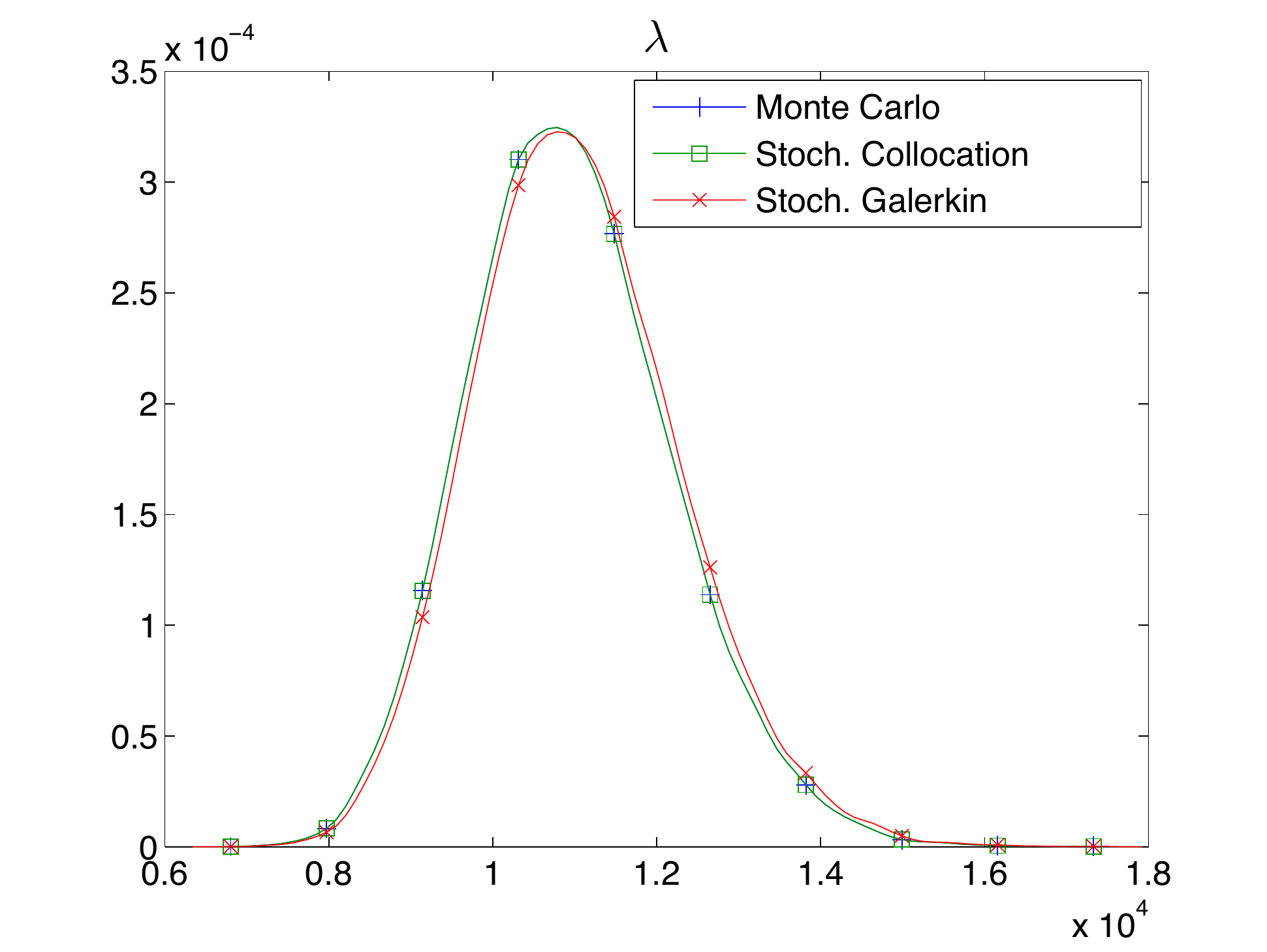}
\includegraphics[width=6.45cm]{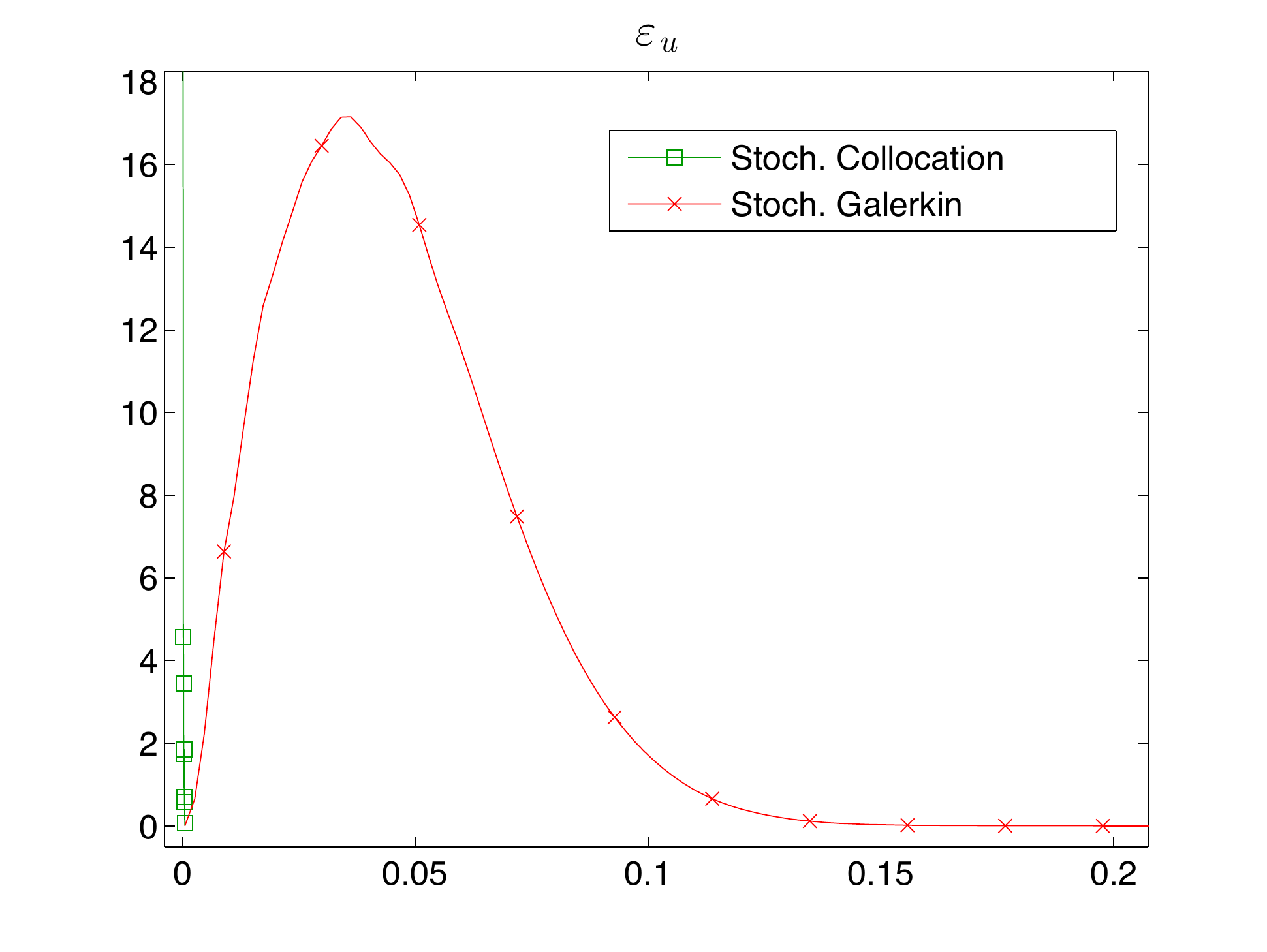}
\includegraphics[width=6.45cm]{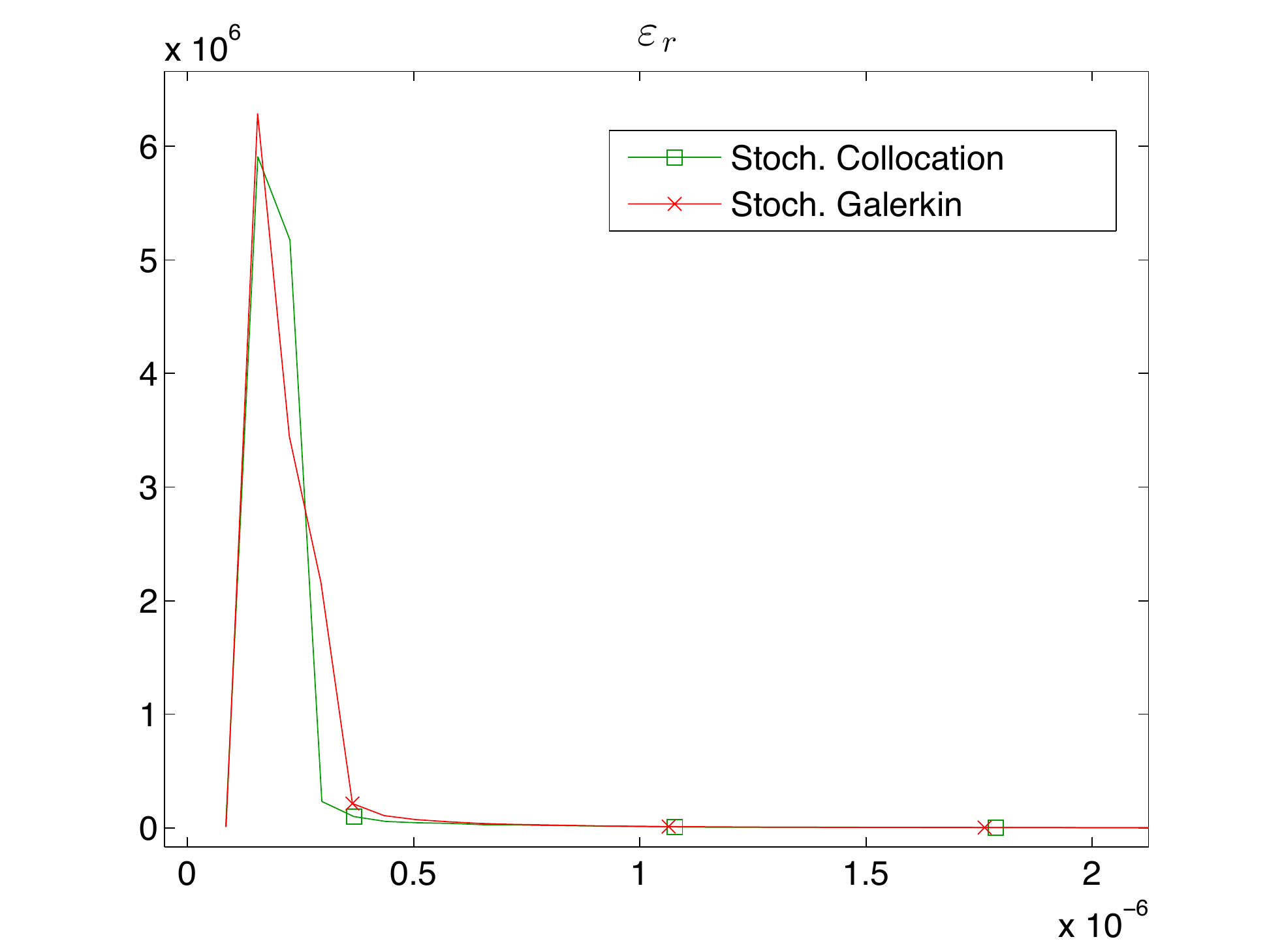}
\includegraphics[width=6.45cm]{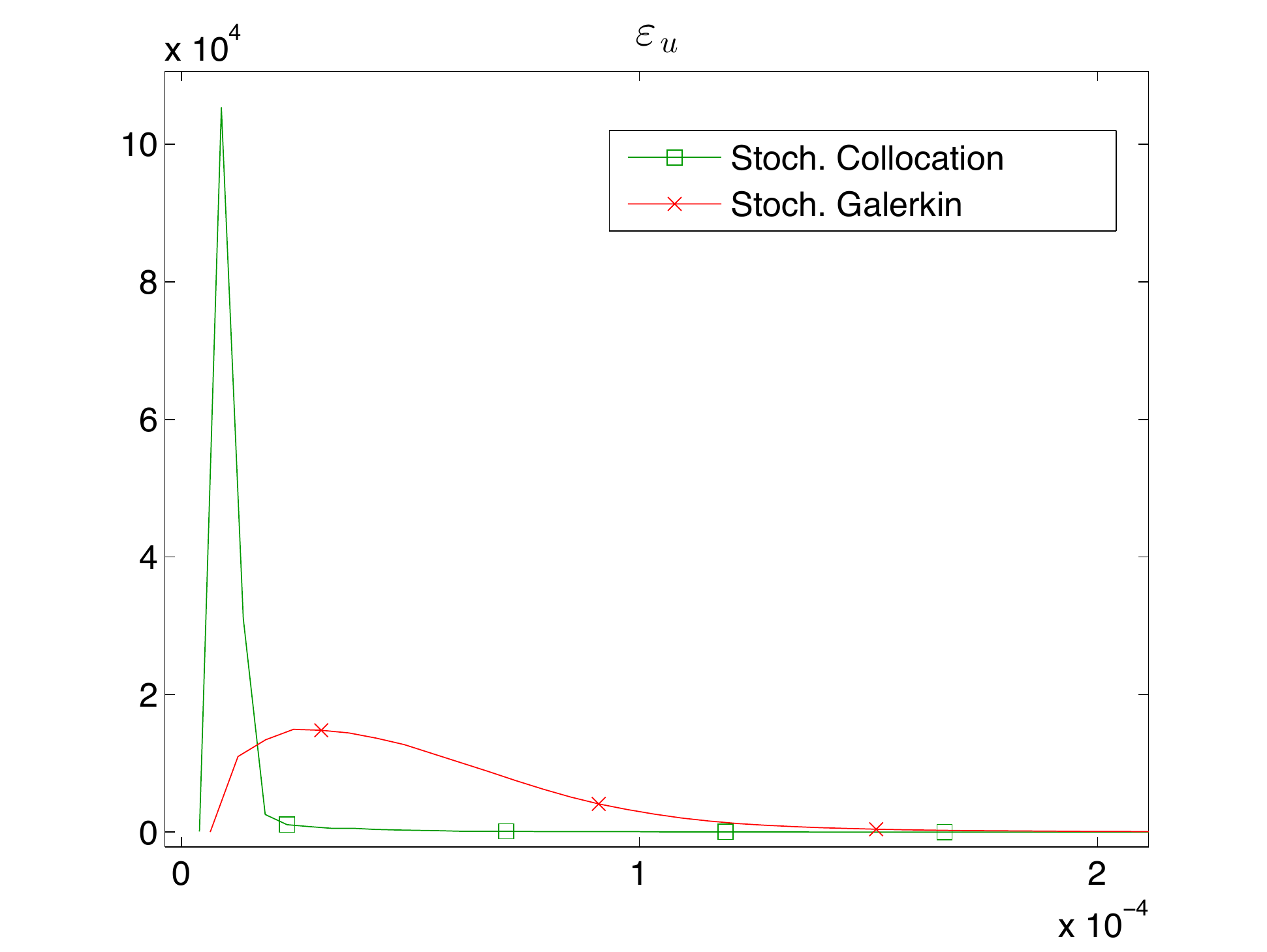}
\end{center}
\caption{Top: pdf estimates of the eigenvalue distribution (left) and the
$\ell^{2}$-norm of the relative eigenvector error~(\ref{eq:eps_u}) (right)
corresponding to the minimal eigenvalue~$\lambda_{1}$ of the Mindlin plate
with $CoV=25\%$ obtained directly using stochastic Rayleigh
quotient~RQ$^{\left(  0\right)  }$. Bottom: pdf estimates of the true
residual~(\ref{eq:true_res}) (left) and of the relative eigenvector
error~(\ref{eq:eps_u}) (right) after five steps of stochastic inverse
iteration.} 
\label{fig:MP-RQ} 
\end{figure}

As before, we first examined the performance of stochastic inverse iteration
and stochastic collocation to identify the smallest eigenvalue. The results
are in Figure~\ref{fig:MP-RQ} and Table~\ref{tab:MP-gPC} presents a comparison
of the first 10 coefficients of the gPC expansion of the smallest eigenvalue
obtained using RQ$^{\left(  0\right)  }$, one and five~steps of stochastic
inverse iteration, and stochastic collocation. Monte Carlo simulation gave
sample mean $1.0952\times10^{4}$ and standard deviation $1.2224\times10^{3}$,
i.e., $CoV\approx11\%$. 
As before, RQ$^{(0)}$ alone provides a close estimate of the eigenvalue
expansion~(\ref{eq:gPC-lambda-u}), and the results of stochastic inverse
iteration and stochastic collocation essentially agree.

\begin{table}[ptbh]
\caption{The first ten coefficients of the gPC expansion of the smallest
eigenvalue of the Mindlin plate with $CoV=25\%$ using $0$, $1$ or $5$ steps of
stochastic inverse iteration, or using stochastic collocation. Here $d$ is the
polynomial degree and $k$ is the index of basis function in expansion
(\ref{eq:gPC-lambda-u}). } 
\label{tab:MP-gPC}
\begin{center} 
\begin{tabular}
[c]{|c|c|c|c|c|c|}\hline
$d$ & $k$ & \multicolumn{1}{|c}{RQ$^{(0)}$} & \multicolumn{1}{|c|}{SII$^{(1)}%
$} & \multicolumn{1}{|c|}{SII$^{(5)}$} & \multicolumn{1}{c|}{SC}\\\hline
0 & 0 & 11,044.1637 & 10,960.2185 & 10,954.8258 & 10,954.8256\\\hline
\multirow{3}{*}{1} & 1 & 1227.0431 & 1217.0164 & 1216.3543 & 1216.3548\\
& 2 & 0 & 0 & 0 & 0\\
& 3 & 0 & 0 & 0 & 0\\\hline
\multirow{6}{*}{2} & 4 & 103.4832 & 97.6116 & 97.5149 & 97.5167\\
& 5 & 0 & 0 & 0 & 0\\
& 6 & 0 & 0 & 0 & 0\\
& 7 & 40.8972 & -15.9359 & -19.7006 & -19.6992\\
& 8 & 0 & 0 & 0 & 0\\
& 9 & 40.8972 & -15.9359 & -19.7006 & -19.6992\\\hline
\end{tabular}
\end{center}
\end{table}

\begin{figure}[ptbh]
\begin{center}
\includegraphics[width=6.45cm]{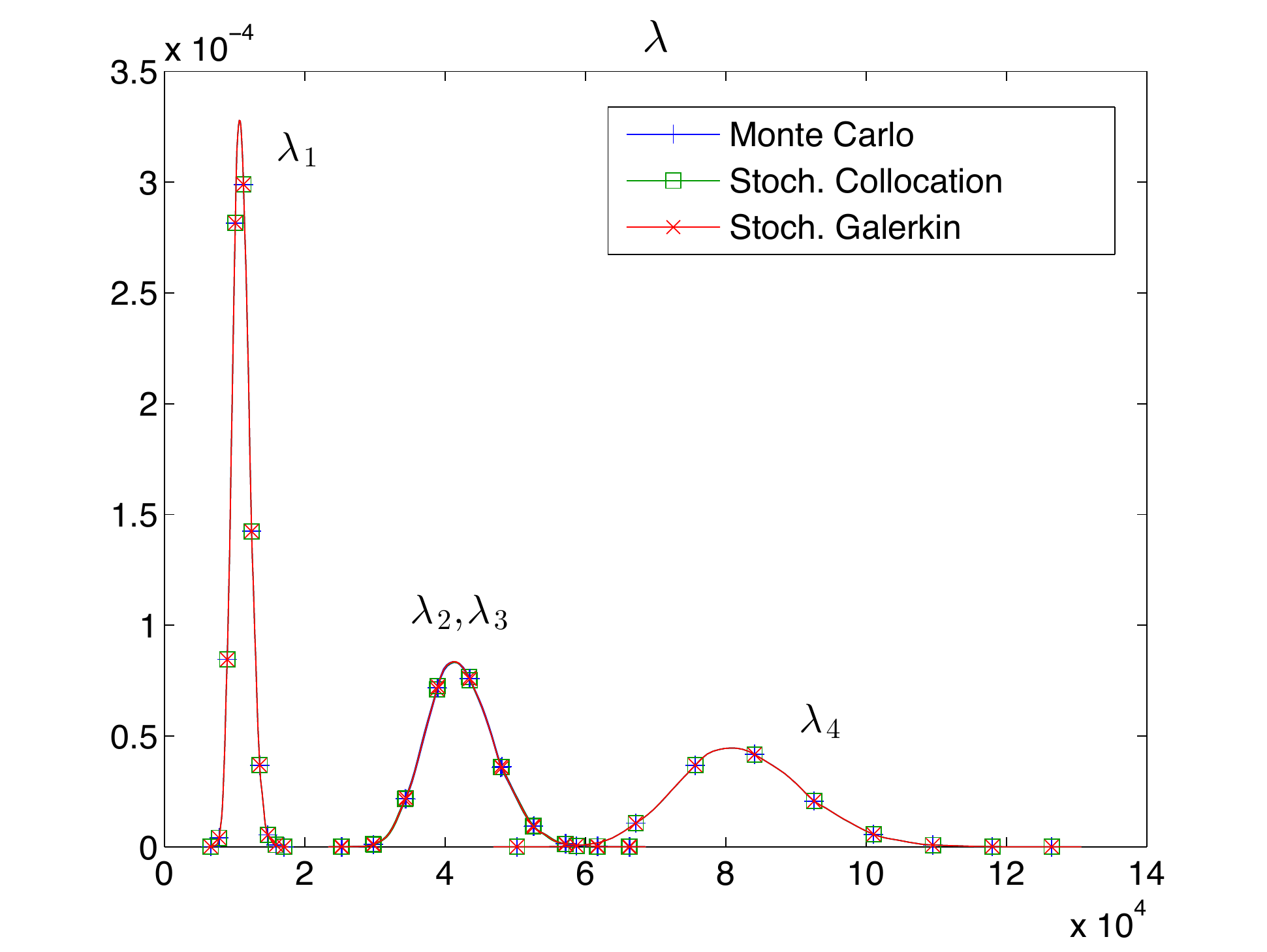}
\includegraphics[width=6.45cm]{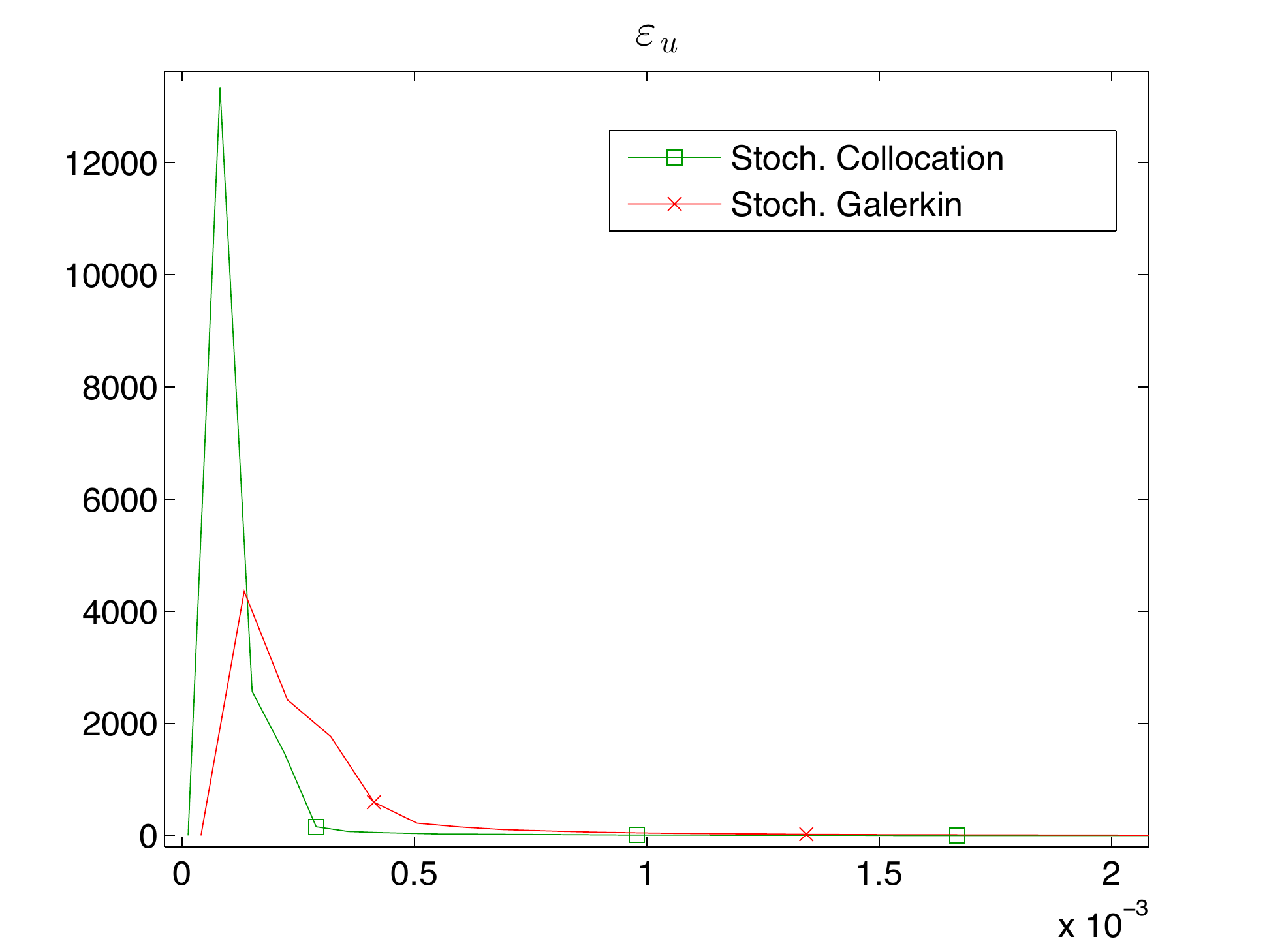}
\includegraphics[width=6.45cm]{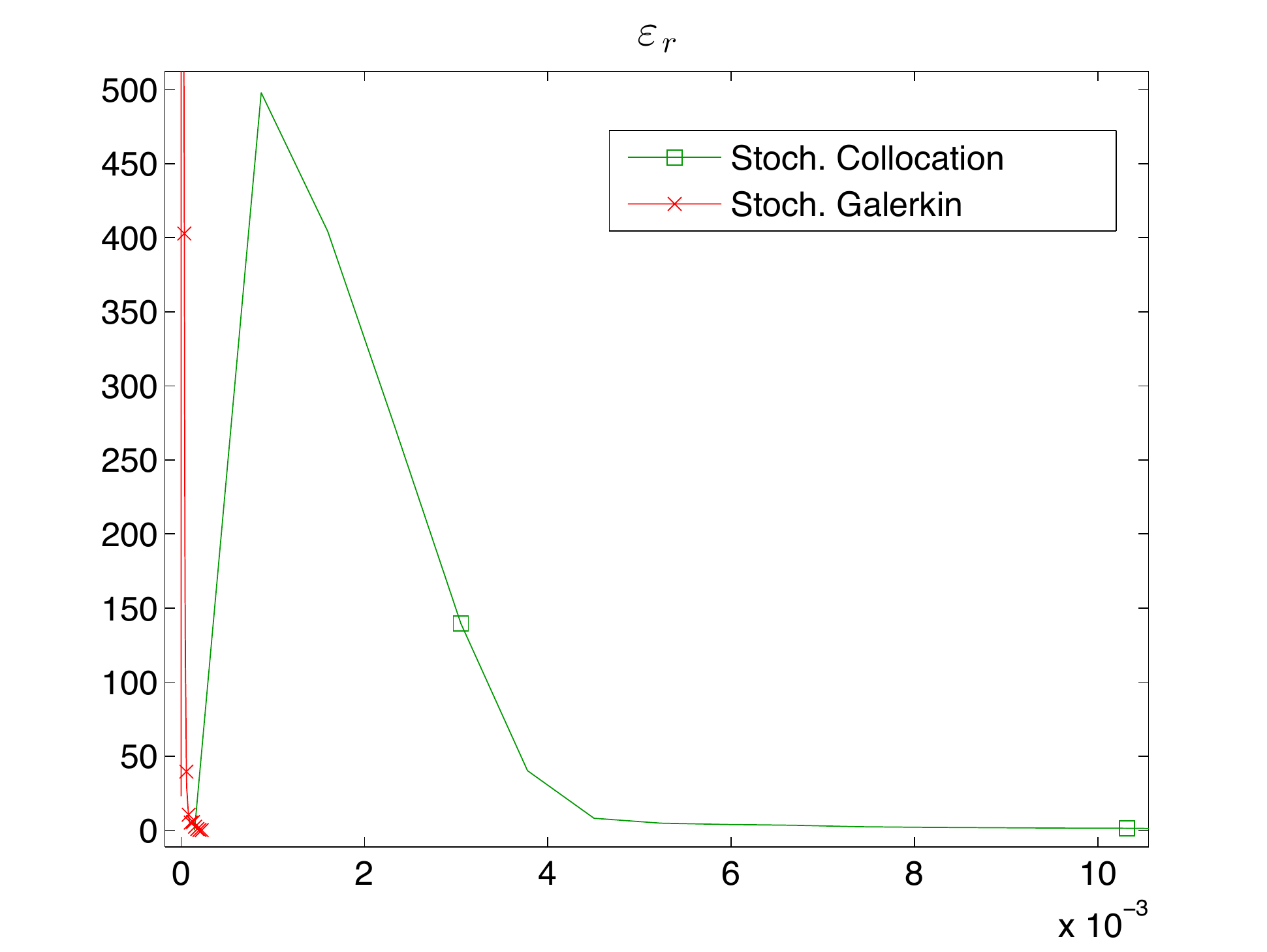}
\includegraphics[width=6.45cm]{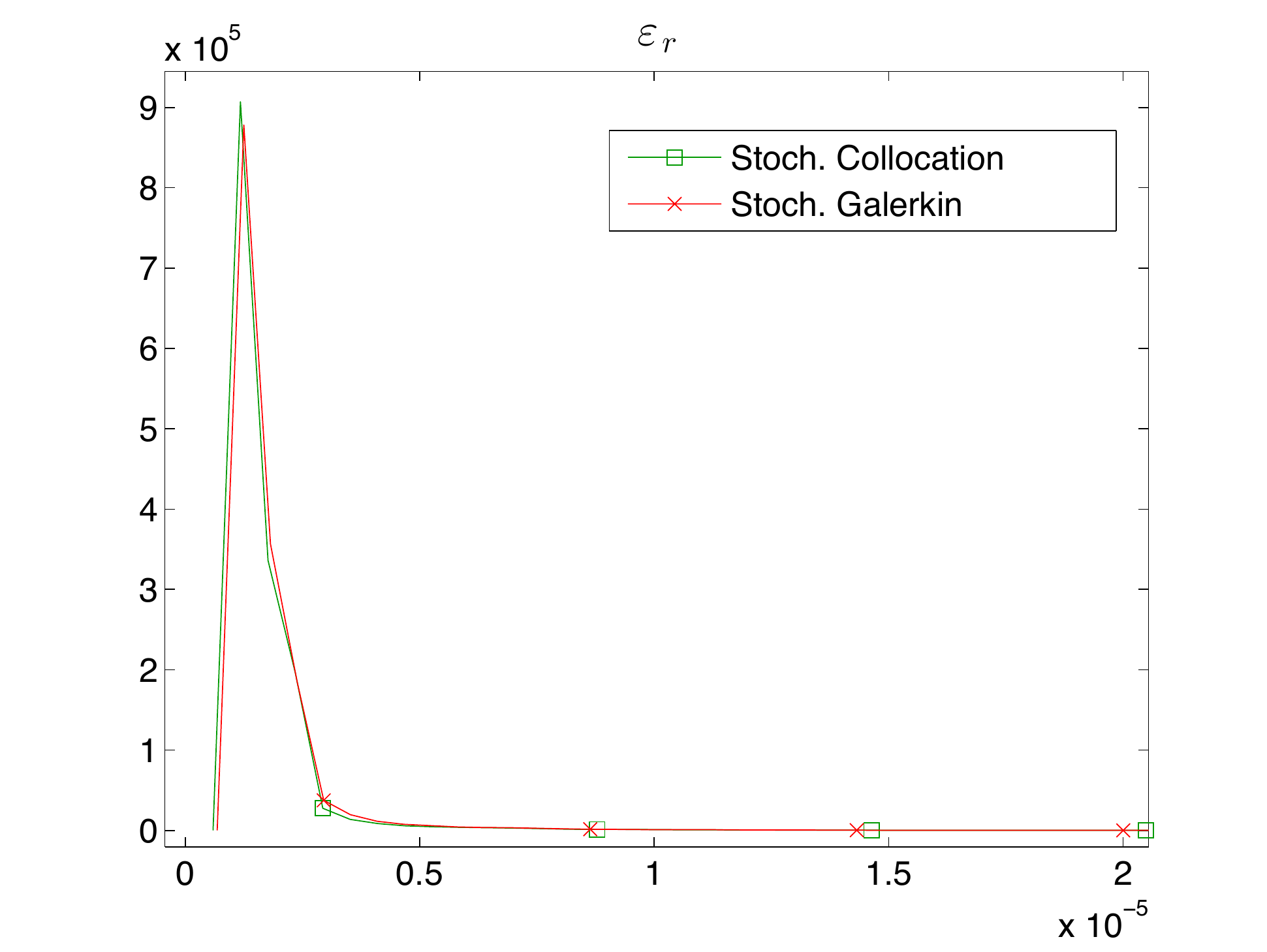}
\includegraphics[width=6.45cm]{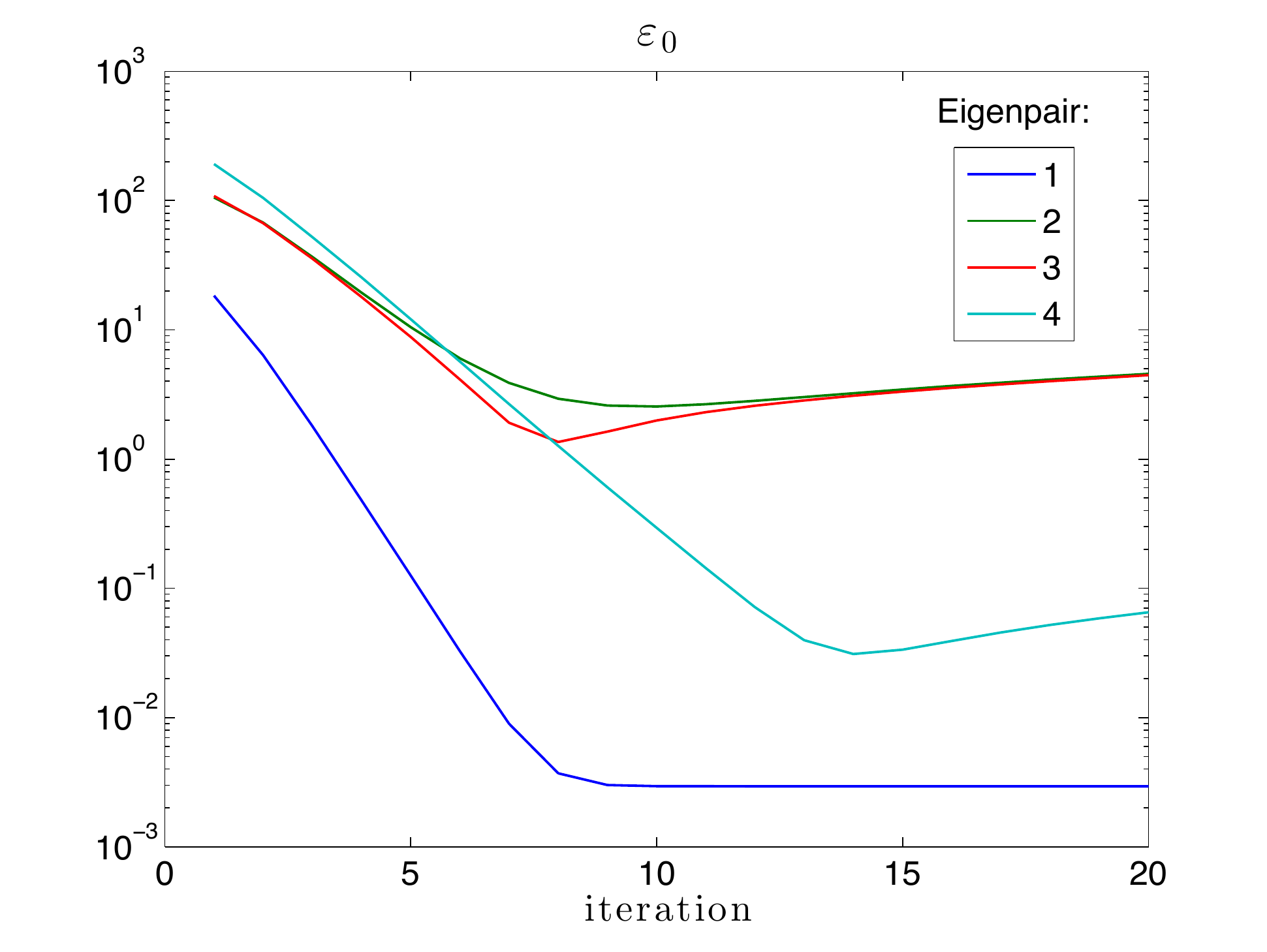}
\includegraphics[width=6.45cm]{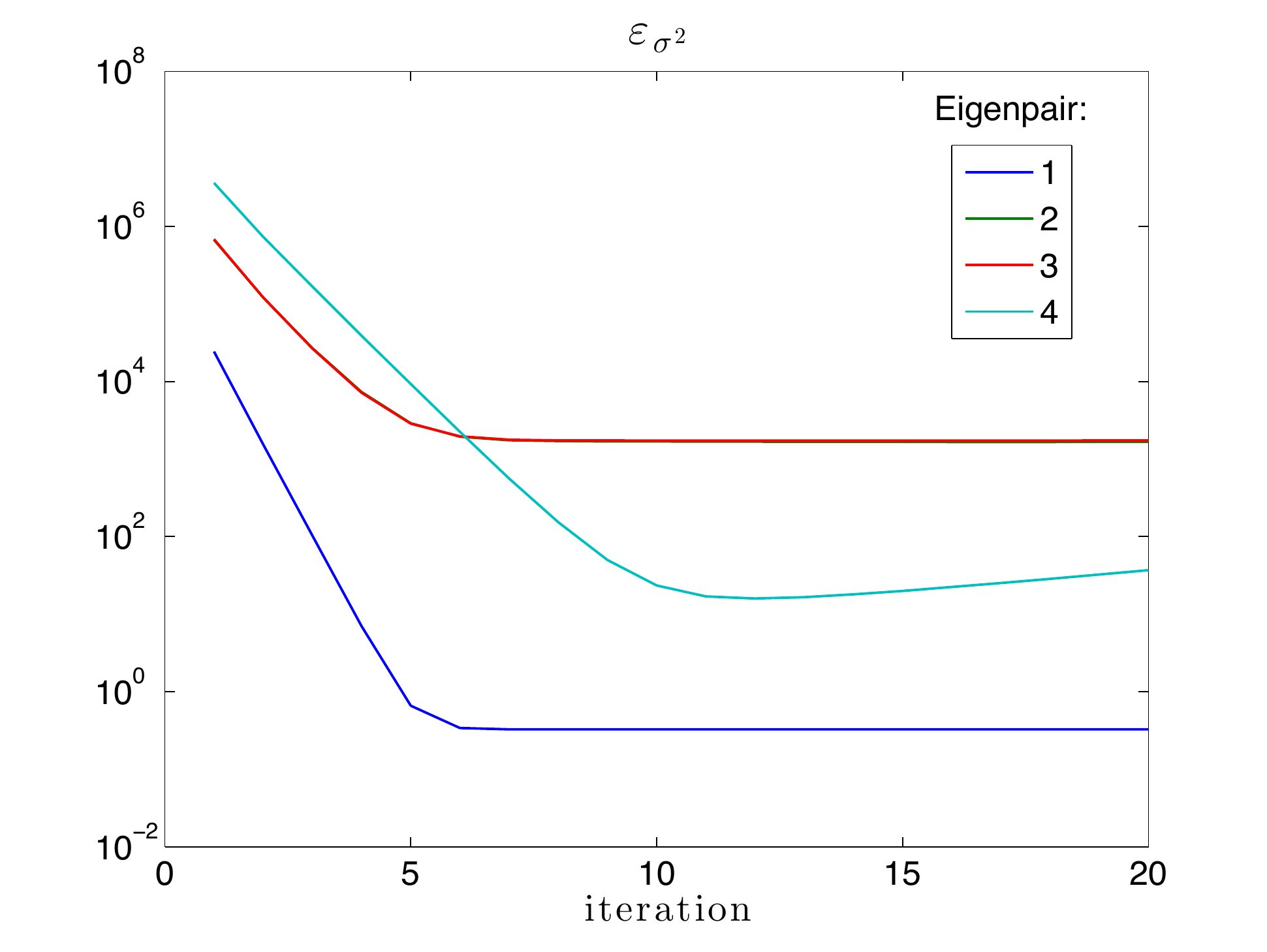}
\end{center}
\caption{Top: pdf estimates of the four smallest eigenvalues (left) and of the
$\ell^{2}$-norm of the relative eigenvector error~(\ref{eq:eps_u})
corresponding to the fourth smallest eigenvalue of the Mindlin plate (right).
Middle: pdf estimates of the true residual~(\ref{eq:true_res}) corresponding
to the second (left) and fourth (right) eigenvector. Bottom: convergence
history of the two indicators $\varepsilon_{0}$ and $\varepsilon_{\sigma^{2}}$
from~(\ref{eq:eps}).} 
\label{fig:MP-min4} 
\end{figure}

Next, we used stochastic inverse subspace iteration to identify the four
smallest eigenvalues. The results are in Figure~\ref{fig:MP-min4}. 
It can be seen that the distributions of all four eigenvalues match and, in
particular, the distributions of the repeated eigenvalues$~\lambda_{2}$
and$~\lambda_{3}$ overlap. However, it also appears that stochastic
collocation exhibits some difficulties detecting the subspace corresponding to
$\lambda_{2}$, whereas the distribution of$~\varepsilon_{r}$ corresponding
to$~\lambda_{4}$ suggests that in this case stochastic inverse subspace
iteration and stochastic collocation methods are in excellent agreement.

\section{Conclusion}

\label{sec:conclusion}We studied random eigenvalue problems in the context of
spectral stochastic finite element methods. We formulated the algorithm of
stochastic inverse subspace iteration and compared its performance in terms of
accuracy with stochastic collocation and Monte Carlo simulation. While overall
the experiments indicate that in terms of accuracy all three methods are quite
comparable, we also highlighted some differences in their methodology. In the
stochastic inverse subspace iteration we formulate and solve a global
stochastic Galerkin problem in order to find the coefficients of the
gPC\ expansion of the eigenvectors. The coefficients of the eigenvalue
expansion are computed from a stochastic version of the Rayleigh quotient. In
fact, we found that
the coefficients of the eigenvector expansion corresponding to the underlying
mean-value problem, with the coefficients of the higher order terms set to
zero, provide a good estimate of the probability distribution 
of the corresponding eigenvalue.
From our experiments it also appears that the stochastic inverse subspace
iteration is not robust when the nature of the eigenvalues is very different, for
example, due to a badly conditioned problem. Moreover, the performance of
inverse iteration for interior eigenvalues seems to be sensitive to the choice
of the shift. Nevertheless, we were able to successfully resolve both issues
by deflation of the eigenvalues of the mean matrix. The algorithm also
performs well in cases when the spectrum is clustered and even for repeated
eigenvalues. However, a unique description of stochastic subspaces
corresponding to repeated eigenvalues, which would allow a comparison of
different bases, is more delicate~\cite{Ghosh-2012-ISA} and is not addressed
here.

We briefly comment on computational cost.
Stochastic inverse subspace iteration is a computational intensive algorithm 
because  it requires repeated solves with the global stochastic Galerkin matrix. 
However, our main focus here was on the methodology, and we view 
a cost analysis to be beyond the scope of this project.
In our experiments,
we used direct solves for the global stochastic Galerkin system, and 
for deterministic eigenvalue problems required by 
the sampling (collocation and Monte Carlo) methods, we used 
the \texttt{Matlab} function \texttt{eig}. 
Many other strategies can be brought to this discussion, for example preconditioned
Krylov subspace methods, e.g., \cite{Sousedik-2014-THP,Sousedik-2014-HSC},
to approximately solve the Galerkin systems, and state-of-the art iterative 
eigenvalue solvers for the sampling methods.
Moreover, the solution of Galerkin systems is also a topic 
of ongoing study~\cite{Matthies-2013-SSS}.
Finally, we note that an appealing feature of the Galerkin approach is
that it allows solution of the random eigenvalue problem only approximately,
performing zero (in case of the stochastic Rayleigh quotient) or only a few
steps of the stochastic iteration, unlike the Monte Carlo and the stochastic
collocation methods which are based on sampling.

\vspace{.1in}
{\bf Acknowledgement.}  We thank the two reviewers and the Associate Editor 
for careful readings of the manuscript and many helpful suggestions.

\bibliographystyle{siam}
\bibliography{eig.bib}

\end{document}